\documentclass[11pt]{amsart}
\usepackage{amsthm}
\usepackage{enumitem}
\usepackage{moreenum}
\usepackage{amssymb,verbatim}
\usepackage{upgreek}
\usepackage{xspace,cmap}
\usepackage{graphicx}
\usepackage{csquotes}
\usepackage{stmaryrd}
\usepackage[colorlinks,citecolor=blue,urlcolor=blue,linkcolor=blue]{hyperref}

\renewcommand{\restriction}{\mathbin\upharpoonright}

\newtheorem*{theorem*}{Theorem}

\newtheorem*{question*}{Question}
\newtheorem{theorem}{Theorem}[section]

\newtheorem{prop}[theorem]{Proposition}
\newtheorem{claim}{Claim}[theorem]
\newtheorem{subclaim}{Subclaim}[claim]
\newtheorem{subsubclaim}{Subsubclaim}[subclaim]
\newtheorem{lemma}[theorem]{Lemma}
\newtheorem{cor}[theorem]{Corollary}

\newtheorem{question}{Question}
\newtheorem{fact}[theorem]{Fact}

\theoremstyle{definition}
\newtheorem*{definition*}{Definition}

\newtheorem{definition}[theorem]{Definition}
\newtheorem{notation}[theorem]{Notation}

\newtheorem{setup}{Setup}

\newtheorem*{MAIN}{Main Theorem}

\theoremstyle{remark}
\newtheorem{remark}[theorem]{Remark}

\newcommand\HOD{\textnormal{HOD}}
\newcommand\OD{\textnormal{OD}}

\newcommand\cat[1]{{}^\curvearrowright #1}

\def\br{\blacktriangleright}

\DeclareMathOperator{\supp}{supp}
\DeclareMathOperator{\crit}{crit}

\DeclareMathOperator{\id}{id}
\DeclareMathOperator{\Suc}{Suc}
\DeclareMathOperator{\Reg}{Reg}
\DeclareMathOperator{\reg}{reg}

    \newcommand{\one}{\mathop{1\hskip-3pt {\rm l}}\nolimits}

\makeatletter
\newcommand{\tpitchfork}{%
  \vbox{
    \baselineskip\z@skip
    \lineskip-.52ex
    \lineskiplimit\maxdimen
    \m@th
    \ialign{##\crcr\hidewidth\smash{$-$}\hidewidth\crcr$\pitchfork$\crcr}
  }%
}
\makeatother

\def\s{\subseteq}
\def\forces{\Vdash}

\DeclareMathOperator{\otp}{otp}

\DeclareMathOperator{\acc}{acc}

\DeclareMathOperator{\ro}{RO^+}

\DeclareMathOperator{\cf}{cf}

\renewcommand{\mid}{\mathrel{|}\allowbreak}

\newcommand{\dom}{\mathop{\mathrm{dom}}\nolimits}

\newcommand{\dP}{\mathbb{P}}
\newcommand{\dQ}{\mathbb{Q}}
\newcommand{\dR}{\mathbb{R}}
\newcommand{\dCS}{\mathbb{CS}}
\newcommand{\dC}{\mathbb{C}}

\title[HOD-supercompactness]{HOD-supercompactness and the optimality of Woodin's HOD dichotomy}
\author[Jakob]{Hannes Jakob}
\address[Jakob]{Department of Mathematics, College of Science, University of North Texas, Denton TX, 76201, USA}
\email{hannes.jakob@unt.edu}
\urladdr{https://hannesjakob.github.io}
\author[Kaplan]{Eyal Kaplan}
\address[Kaplan]{Department of Mathematics, Carnegie Mellon University, Pittsburgh PA, USA}
\email{eyalk@andrew.cmu.edu}
\urladdr{https://eyalkpl.github.io/website}

\author[Poveda]{Alejandro Poveda}
\address[Poveda]{Departamento de Matemáticas, CUNEF Universidad,  Madrid, 28040, Spain.}
\email{alejandro.poveda@cunef.edu}
\urladdr{https://www.alejandropovedaruzafa.com}

\begin{document}
\maketitle

\begin{abstract}
We  answer  a question of Goldberg--Poveda from \cite{GOP} by showing that the HOD hypothesis is consistent with $\kappa$ being the first $\HOD$-supercompact cardinal carrying a club of cardinals that are regular in $\HOD$. This configuration yields a sharp distinction between $\HOD$-supercompactness and extendibility when it comes to covering properties between HOD and the set-theoretic universe. To obtain it, we develop a new forcing construction: a homogeneous Radin-like poset. Finally, we also address a conjecture of Woodin from~\cite{WOEM1} concerning successors of singular cardinals in $\HOD$.
\end{abstract}

\section{Introduction}
Questions of definability have played a central role in set theory since its origins. In a 1946 lecture at Princeton \cite{Godel}, G\"odel introduced the class of \emph{hereditarily ordinal definable sets}, nowadays denoted by $\HOD$. Informally, a set $X$ is ordinal definable if it is the unique object satisfying some set-theoretic property that can be described using finitely many ordinal parameters. A set $X$ belongs to $\HOD$ if $X$ is ordinal definable, every element of $X$ is ordinal definable, every element of every element of $X$ is ordinal definable, and so on; see \S\ref{sec:preliminaries} for the formal definition. Subsequent work of Myhill--Scott \cite{Myhill}, Vop\v{e}nka \cite{Vopenka}, and others began to reveal the surprisingly rich structure of the class of ordinal definable sets.

\smallskip

The class $\HOD$ is {transitive}, contains all ordinals, and satisfies the usual axioms of set theory; namely, the \emph{Zermelo--Fraenkel axioms together with the axiom of choice} (ZFC). In modern set-theoretic jargon, $\HOD$ is an \emph{inner model} of the set-theoretic universe $V$. Inner models provide canonical approximations to the universe, and a central theme in set theory is to understand how much of the structure of $V$ they can capture. \emph{G\"odel's constructible universe} $L$ (\cite{GodelConstructible}) initiated this program, and subsequent work of Jensen, Mitchell, Steel, Neeman, Woodin, and many others (\cite{Steel, Zeeman}) have produced increasingly sophisticated canonical inner models capable of accommodating stronger and stronger large-cardinal axioms.

\smallskip

The position of $\HOD$ within this picture is particularly intriguing. Unlike the fine-structural inner models just mentioned, $\HOD$ is defined directly from the ambient universe and is highly sensitive to forcing. Nevertheless, it contains every canonical inner model constructed to date and can accommodate virtually all  large cardinals (see e.g., \cite{BagPov}). Thus the comparison between $\HOD$ and $V$ provides a natural testing ground for a basic problem in the foundations of mathematics; namely,  \emph{to what extent does the definable part of the universe reflect the structure of the universe itself?}

\smallskip

\emph{Large-cardinal} axioms postulate the existence of infinite cardinals whose existence cannot be established from ZFC. They form a remarkable hierarchy of principles which calibrates the consistency strength of a broad range of mathematical theories and has played a decisive role in the modern theory of independence \cite{Kan,Koel}. Large cardinals have also had consequences for more familiar mathematical objects --- for instance, for regularity properties of sets of real numbers \cite{Solmod,DPT,PT} and for structural questions in algebra \cite{MagShe,Eklof,CortesPoveda}. At sufficiently high levels of this hierarchy, their interaction with $\HOD$ becomes unexpectedly rigid.

\smallskip

In recent years, Woodin's groundbreaking work \cite{WOEM1, midrasha} has stimulated considerable activity around the interaction between $\HOD$, large cardinals, compactness phenomena, and forcing. One direction has sought to understand how dramatically the cardinal and large-cardinal structure of $\HOD$ can diverge from that of the ambient universe, using sophisticated forcing techniques; see, among others, \cite{GitikMer,BenUng,BenHay,GOP,PovOmega}. A complementary line of research has uncovered compactness and covering principles that constrain this divergence in the presence of strongly compact cardinals and singular cardinals of uncountable cofinality; see \cite{GoldbergHODStrongCompact,GolPov,GOP,CompactnessMagidor}. These developments form part of a broader renewed interest in definable inner models and their interaction with large cardinals, also visible in the recent influential work on \emph{Nairian models}~\cite{Nairian}.

\smallskip

A fundamental discovery of Woodin is the following dichotomy.

\begin{theorem*}[Woodin's HOD Dichotomy~\cite{WOEM1}]
Suppose that $\delta$ is an extendible cardinal. Then exactly one of the following holds:
\begin{enumerate}
    \item $\HOD$ is a weak extender model for the supercompactness of $\delta$;
    \item every regular cardinal $\lambda\geq\delta$ is $\omega$-strongly measurable in $\HOD$.
\end{enumerate}
\end{theorem*}

Let us briefly explain the first alternative, which will be central to this paper. An inner model $N$ is a \emph{weak extender model for the supercompactness of $\kappa$} if, for every cardinal $\lambda\geq\kappa$, there is a normal fine $\kappa$-complete ultrafilter $U$ on $\mathcal P_\kappa(\lambda)$ such that $\mathcal P_\kappa(\lambda)\cap N\in U$ and $U\cap N\in N$. Thus the supercompactness measures of the ambient universe can be chosen so as to concentrate on $N$ while leaving behind measures internal to $N$. In particular, if $N$ is a weak extender model for the supercompactness of $\kappa$, then $\kappa$ is itself supercompact in $N$. Additionally, these models enjoy of strong forms of \emph{covering} and correctness properties above $\kappa$ (see \cite{midrasha}). Among others, there is the \emph{weak covering} property which states the following: if $\lambda>\kappa$ is a singular cardinal then it is singular in $N$ and $\lambda^{+N}=\lambda^+$.

\smallskip

The two alternatives of Woodin's dichotomy therefore describe radically different pictures. In the first, $\HOD$ retains a substantial part of the large-cardinal structure of the universe. In the second, $\HOD$ is very far from $V$: all sufficiently large regular cardinal are large cardinals inside $\HOD$.

Woodin's \emph{HOD hypothesis} (see \S\ref{sec:preliminaries}) rules out the second alternative. His broader \emph{HOD conjecture} predicts that the HOD hypothesis should follow from suitable large-cardinal assumptions. In this sense, the HOD hypothesis expresses the expectation that, in the presence of sufficiently strong large cardinals, $\HOD$ should eventually resemble the universe of sets.

An important point is that ordinary supercompactness does not by itself enforce this relationship with $\HOD$. Woodin introduced the notion of \emph{$\HOD$-supercompactness} (see Definition~\ref{def:hodsupercompactness}) precisely to capture the additional coherence that is required. 
Every extendible cardinal is $\HOD$-supercompact (\cite{WOEM1}), while the consistency strength of $\HOD$-supercompactness is only that of a supercompact cardinal. The distinction between supercompactness and $\HOD$-supercompactness is nevertheless genuine. The first observation in this direction is due to Sargsyan~\cite{Sar08}. More recently, Goldberg and Poveda showed that it is consistent for $\kappa$ to be supercompact and for the HOD hypothesis to hold, while $\HOD$ is nevertheless not a weak extender model for the supercompactness of $\kappa$. Thus ordinary supercompactness, even together with the HOD hypothesis, does not suffice to recover the first alternative of Woodin's dichotomy.
For $\HOD$-supercompact cardinals the situation is completely different: if $\kappa$ is $\HOD$-supercompact and the $\HOD$ hypothesis holds, then $\HOD$ is a weak extender model for the supercompactness of $\kappa$ (Theorem~\ref{thm:HODsupercompactWoodin}).

 Thus, despite having only the consistency strength of supercompactness, $\HOD$-supercompactness recovers much of the extendible-like rigidity of $\HOD$ \emph{above} the large cardinal. This raises the question at the heart of this paper:
\begin{question*}
How much extendible-like rigidity of $\HOD$ does a $\HOD$-supercom\-pact cardinal $\kappa$ impose below itself?    
\end{question*}

Recent work of Goldberg--Poveda \cite{GOP} shows that the situation below an ordinary supercompact cardinal can be remarkably wild. On the one hand, if the HOD hypothesis holds and $\kappa$ is supercompact, then there are unboundedly many singular cardinals $\lambda<\kappa$ which remain singular in $\HOD$ and satisfy $\lambda^{+\HOD}=\lambda^+$.\footnote{This degree of resemblance was also observed by Woodin.} Thus some degree of correctness below $\kappa$ is unavoidable. On the other hand, this cannot in general be strengthened from ``unboundedly many'' to ``club many''. Indeed, it is consistent with the HOD hypothesis that the first supercompact cardinal $\kappa$ carries a club $C\subseteq\kappa$ such that every $\lambda\in C$ is regular in $\HOD$ and $\lambda^{+\HOD}<\lambda^+$. 

\begin{question}[{\cite{GOP}}]
Is it consistent with the $\HOD$ hypothesis that $\kappa$ is $\HOD$-supercompact and there is a club $C\s\kappa$ consisting of $\HOD$-regulars?   
\end{question}

 The question is particularly natural because the corresponding configuration is incompatible with extendibility, whereas $\HOD$-supercompactness reproduces many of the consequences of extendibility for the structure of $\HOD$. Our main theorem gives a positive partial answer to this question.
\begin{MAIN}
Assume the existence of a supercompact cardinal. Then there is a cardinal-preserving generic extension in which the HOD hypothesis holds and the first $\HOD$-supercompact cardinal $\kappa$ carries a club $C\subseteq\kappa$ consisting of cardinals which are measurable in $\HOD$.\footnote{Actually, by starting from a stronger large cardinal assumptions, we can upgrade this to encompass other large cardinal properties in $\HOD$; such as supercompactness in $\HOD$.}
\end{MAIN}

 The proof of the Main Theorem requires a new forcing construction: a \emph{homogeneous Radin-like forcing}. At an intuitive level, our aim is to reproduce one of the characteristic effects of Radin forcing --- the introduction of a club of former regular cardinals --- while retaining enough homogeneity to control $\HOD$ of the final generic extension and enough large-cardinal structure to preserve $\HOD$-supercompactness.  The earlier construction from~\cite{GOP}, based on Cummings--Friedman--Golshani's work  \cite{CumFriGol}, does not seem to yield the configuration described in the Main Theorem. The basic issue is the lack of homogeneity  and the effects that it has on HOD-supercompactness.
 
 Our construction achieves the desired configuration by melting two forcing techniques: The first originates in Gitik's classical work on adding closed unbounded sets consisting of former regular  cardinals \cite{GitikOnClosed};
the second is the nonstationary-support iteration of Prikry-type forcings developed by Ben-Neria and Unger~\cite{BenUng}, which provides the degree of homogeneity needed to control the HOD of a forcing extension. The similarities and differences between our methods and theirs are explained at length in \S\ref{sec:comparison with Gitik}.

Our expectation is that this novel homogenous Radin-like forcing will find other applications in the study of HOD beyond those pursued here.

\smallskip

Separately, the paper also addresses a  conjecture of Woodin concerning successors of singular cardinals in $\HOD$. In~\cite[p.~299]{WOEM1}, he  observed that if $\gamma$ is a singular strong limit cardinal and $\gamma^{+\HOD}<\gamma^+$, then either $\cf(\gamma^{+\HOD})=\omega$ or $\cf(\gamma^{+\HOD})=\cf(\gamma)$, and conjectured that the first possibility should always occur.  We prove this in the negative:

\begin{theorem*}
Assume the existence of a $\kappa^+$-supercompact cardinal $\kappa$ with $2^\kappa=\kappa^+$. Then there is a generic extension in which $\kappa$ is a singular strong limit cardinal, $\cf(\kappa^{+\HOD})=\omega_1$, and $\kappa^{+\HOD}<\kappa^+$. Moreover, $\cf^\HOD(\kappa)=\omega_1$.
\end{theorem*}

The paper is organized as follows. In Section~\ref{sec:preliminaries} we review the relevant background on $\HOD$, weak extender models, large cardinals, and iterations of Prikry-type forcings.  In Section~\ref{sec: the iteration} we define the global iteration and establish its preservation properties. Section~\ref{sec: basic building block} develops the local forcing notions used in the construction extending Gitik's machinery to the Mitchell orders required for our application. Section~\ref{sec: homogeneity} contains the homogeneity analysis. In Section~\ref{sec:PreservingHODsupercompacts} we prove the preservation of $\HOD$-supercompactness and establish the Main Theorem. Finally, Section~\ref{sec:SuccessorsOfSingulars} gives the negative answer to Woodin's conjecture concerning successors of singular cardinals in $\HOD$.

\section{Preliminaries}\label{sec:preliminaries}

\subsection*{Notations and conventions}
We follow standard set-theoretic notation. If $M$ is an inner model and $\lambda$ is a cardinal, we write $\lambda^{+M}$ for the successor of $\lambda$ as computed in $M$.  We will also use expressions such as ``$\lambda$ is $M$-regular'' to mean that $M\models``\lambda$ is regular''. When no model is indicated, all cardinal-theoretic notions are computed in the ambient universe.

Given a forcing poset $\mathbb{P}$, we shall adopt the convention that $p\leq q$ means that \emph{$p$ is stronger than $q$}. We denote the trivial condition of a forcing $\mathbb{P}$ by $\one_{\mathbb{P}}$, or simply by $\one$ when $\mathbb{P}$ is clear from context.  We usually suppress the subscript in the forcing relation when the relevant forcing is unambiguous.

Given a forcing poset $\mathbb{P}$, a $\mathbb{P}$-name $\tau$, and a formula $\varphi(x)$ in the forcing language of $\mathbb{P}$, we say that \emph{$\tau$ is forced to satisfy $\varphi$} whenever
$\one_{\mathbb{P}}\forces_{\mathbb{P}}\varphi(\tau)$. When the context allows it, we will also write $\forces \varphi(\tau)$ rather than the more informative $\one_{\mathbb{P}}\forces_{\mathbb{P}}\varphi(\tau)$.

\subsection{HOD}
The primary object of study in this manuscript is Gödel's universe of Hereditary Ordinal Definable sets, $\HOD$ (see \cite{Godel}):
\begin{definition}
    Let $X$  be a set. The class of  \emph{Ordinal Definable sets relative to $X$}, denoted $\OD_X$,  consists of all sets $A$ for which there is a formula $\phi(z, x_0,\dots, x_{m-1},y_0\dots, y_{n-1})$ in the language of set theory, sets $a_0,\dots, a_{m-1}$ in  $X$ and ordinals $\beta, \alpha_0,\dots, \alpha_{n-1}$ such that $A$ can be presented as $$A=\{z\in V_\beta \mid V_\beta \models \varphi(z,a_0,\dots,a_{m-1},\alpha_0,\dots, \alpha_{n-1})\}.$$
    The class of \emph{Hereditarily Ordinal Definable sets relative to $X$}, denoted $\HOD_X$, is the class of all $A$ such that $\mathrm{trcl}(\{A\})\s \OD_X$.

   We shall write $\OD$ and $\HOD$ for $\OD_\emptyset$ and $\HOD_\emptyset$, respectively.
\end{definition}
$\HOD$ is an \emph{inner model}; namely, a transitive $\in$-structure satisfying the ZFC axioms and containing all the ordinals (through the paper this is what we shall mean when we refer to an \emph{inner model}). It also possesses a canonically definable global well-order: namely, for a set $A\in \HOD$, let $\theta_{A}$ be the least beth fixed-point such that $A$ is ordinal definable in $V_{\theta_{A}}$, and let $\langle n^A, \alpha^A_0,\ldots, \alpha^A_{m_A} \rangle\in V_{\theta_{A}}$ be the lexicographically least finite sequence for which $n^A\in \mathbb{N}$ codes a formula $\varphi$ that defines $A$ in $V_{\theta_A}$ from the ordinals $\langle \alpha^A_0,\ldots, \alpha^A_{m_A} \rangle$. We then order $A<B$ if and only if $$\langle \theta^A, n^A, \alpha^A_0,\ldots, \alpha^A_{m_A} \rangle<_{\text{lex}}  \langle \theta^B,  n^B, \alpha^B_0,\ldots, \alpha^B_{m_B} \rangle.$$\label{wellorderofHOD}
By forcing with a class iteration that codes every set of ordinals  into the continuum function,  McAloon proved the following well-known theorem:

\begin{fact}[McAloon, {\cite{McAloon}}]
  There is a class forcing poset sets yielding a generic extension where the axiom $V=\mathrm{HOD}$ holds.
\end{fact}

A strengthening of $``V=\HOD$" relevant to this paper was introduced by Fuchs--Hamkins--Reitz in \cite{Geology}:

\begin{definition}[Fuchs--Hamkins--Reitz]
The \emph{generic $\HOD$ of $V$}, denoted $\mathrm{gHOD}$, is the intersection of the $\HOD$s of all the generic extensions of $V$. That is, $x\in \mathrm{gHOD}$ if and only if $x\in \HOD^{V[G]}$ for every set forcing $\mathbb{P}$ and $G\s \mathbb{P}$ generic. Equivalently, 
$$\mathrm{gHOD}=\{x\mid \forall \mathbb{P}\, \one\forces_{\mathbb{P}}\check{x}\in\dot{\HOD}\}$$
\end{definition}

\begin{fact}[{\cite[Theorem~15]{Broo2}}]\label{fact:BrookeTaylor}
Suppose the $\mathrm{GCH}$ holds and that $\kappa$ is supercompact. Then, there is a class iteration forcing $V=\mathrm{gHOD}$ and preserving both the $\mathrm{GCH}$ and the supercompactness of $\kappa$.
\end{fact}
We now present a series of  results of Woodin's landmark work \cite{WOEM1}.
\begin{definition}[Woodin, {\cite{WOEM1}}]
   Given a cardinal $\kappa$, an inner model $M$ of ZFC is a \emph{weak extender model for the supercompactness of $\kappa$} (shortly, \emph{WEM for $\kappa$}) if for all \(X\in M\), there is a normal fine \(\kappa\)-complete ultrafilter
\(\mathcal U\) on \(\mathcal{P}_\kappa(X)\cap M\) with \(\mathcal U\cap M\in M\). 
\end{definition}

   $M$ being a WEM for  $\kappa$ implies that  \(\kappa\) is  supercompact in both \(M\) and \(V\). 

   \smallskip

WEM possess strong covering properties  making them close to $V$:
\begin{theorem}[Woodin, {\cite{midrasha}}]
  If $M$ is $\mathrm{WEM}$ for $\kappa$ then $M$ has the $\kappa$-cover and $\kappa$-approximation properties.  In particular, every singular cardinal $\lambda>\kappa$ is singular in $M$ and $\lambda^{+M}=\lambda$. 
\end{theorem}
In \S\ref{sec:PreservingHODsupercompacts} we shall be concerned with the following notion due to Woodin:
\begin{definition}[Woodin, {\cite{WOEM1}}]\label{def:hodsupercompactness}
     $\kappa$ is $\HOD$-supercompact if for all strong limit cardinals $\lambda$, there exists an elementary embedding $j\colon V\to M$ with critical point $\kappa$ such that ${}^{\lambda}M\subseteq M$ and $\HOD^M\cap V_{\lambda}=\HOD\cap V_{\lambda}$.
\end{definition}
Every $\HOD$-supercompact cardinal is supercompact, but the converse is  false as shown by Sargsyan in \cite{Sar08}. For more results see \cite{MoreonHOD}.

\smallskip

Woodin's structural theory of HOD hinges on the following hypothesis:

\begin{definition}[\cite{midrasha}]
    The \textit{$\HOD$ hypothesis} states that there is a proper class of regular cardinals that are not \(\omega\)-strongly measurable in \(\HOD\).
\end{definition}
One of his central theorems is the following:
\begin{theorem}[Woodin~\cite{WOEM1}]\label{thm:HODsupercompactWoodin}
Assume the $HOD$ hypothesis and let $\kappa$ be $\HOD$-supercompact. Then $\HOD$ is a weak extender model for the supercompactness of $\kappa$. In particular, $\kappa$ is supercompact in $\HOD$.
\end{theorem}
The above result is actually sharp:
\begin{theorem}[Goldberg--Poveda, {\cite{GOP}}]
 The following configuration is consistent relative to large cardinals:   the $\HOD$ hypothesis holds, $\kappa$ is supercompact but  not supercompact in $\HOD$.   
\end{theorem}
The $\HOD$ hypothesis together with milder form of extendibility, such as supercompactness, already entail some agreement between $V$ and $\HOD$: 
\begin{theorem}[Woodin; Poveda {\cite{GOP}}]  Assume that $\kappa$ is a supercompact cardinal and the $\HOD$ hypothesis holds. Then, for unboundedly many singular cardinals $\theta < \kappa$, $\theta$ is singular in $\HOD$ and $\theta^{+\HOD}=\theta^+$.\footnote{This theorem was first observed by Woodin. Poveda observed it independently later while working on \cite{GOP}.}  
\end{theorem}
Once again, the above form of agreement is arguably  optimal:
\begin{theorem}[Goldberg--Poveda, {\cite{GOP}}]
   The following configuration is consistent relative to large cardinals:   the $\HOD$ hypothesis holds, $\kappa$ is supercompact and there is a club $C\s \kappa$ consisting of $\HOD$-regular cardinals $\theta$ such that $\theta^{+\HOD}<\theta^+$. 
\end{theorem}

In \S\ref{sec:PreservingHODsupercompacts} we shall prove Theorem~\ref{thm:finally} which produces a similar configuration replacing supercompactness by the stronger notion of $\HOD$-supercompactness. This configuration, however, cannot take place below the first extendible, thus evidencing a sharp disctintion between extendibility and HOD-supercom\-pactness when it comes to the structure of HOD as a weak extender model.

\begin{prop}\label{prop:extendible-weak-covering}
Assume the $\HOD$ hypothesis and let $\kappa$ be an extendible cardinal.
Then there is some $\kappa_*<\kappa$ such that, for every singular cardinal
$\theta\in(\kappa_*,\kappa)$, $\theta$ is singular in $\HOD$ and $\theta^{+\HOD}=\theta^+$.
\end{prop}

\begin{proof}
This follows easily from Goldberg's strongly compact version of the HOD dichotomy
\cite[Corollary~2.9]{GoldbergHODStrongCompact}. Concretely,   if $\kappa_*$
is the first strongly compact cardinal (which necessarily is less than $\kappa$) then the $\HOD$ hypothesis implies that  every  singular cardinal $\theta>\kappa_*$ remains singular in $\HOD$ and
satisfies $\theta^{+\HOD}=\theta^+$.
\end{proof}

\subsection{Forcing}\label{sec: forcing}
In this section we  present some forcing preliminaries needed in the main theorem of this paper. Any related concept  or notation that has not been defined  can be found in  Kunen's text \cite{Kun}.

\smallskip

\begin{definition}
    Let $\dP$ and $\dQ$ be posets and $\Gamma\colon\dP\to\dQ$ any map. We define $\Gamma^*\colon V^{\dP}\to V^{\dQ}$ by induction on the name-rank, letting
    $$\Gamma^*(\tau):=\{\langle\Gamma^*(\sigma),\Gamma(p)\rangle\;|\;\langle\sigma,p\rangle\in\tau\}$$
\end{definition}

Such maps are usually defined in the context of dense embeddings or automorphisms. Typically, homogeneity of posets, which is crucial when dealing with $\HOD$, is defined in terms of the existence of automorphisms. However, in the case of forcing iterations, outright automorphisms are difficult to obtain, since we are dealing with names, so we must weaken the notion somewhat. That this is possible is provided by the following lemma:

\begin{lemma}\label{Lemma: Self-Inverse Auto mod Equiv is DEmb}
    Let $\dP$ be a poset and $\Gamma\colon\dP\to\dP$ an order-preserving map such that for every $p\in\dP$, $p\leq\Gamma^2(p)\leq p$. Then $\Gamma$ is a dense embedding.
\end{lemma}

\begin{proof}
    It  remains to show that the image of $\Gamma$ is dense. However, clearly, for every $p\in\dP$, $\Gamma(\Gamma(p))\leq p$ by assumption.
\end{proof}

Next, we introduce a slight variation of the usual notion of automorphism of a poset -- that will be suitable for our purposes.

\begin{definition}
    Let $\dP$ be a poset and let $\Gamma\colon\dP\to\dP$ be a map. We say that $\Gamma$ is a \emph{self-inverse automorphism modulo equivalence} if $\Gamma$ preserves the order and for every $p\in\dP$, $p\leq\Gamma^2(p)\leq p$.
\end{definition}

\begin{lemma}\label{Lemma: Self-Inverse Auto mod Equiv gives Generic}
    Let $\dP$ be a poset and $\Gamma\colon\dP\to\dP$ a self-inverse automorphism modulo equivalence. Let $G\subseteq\dP$ be a generic filter. Then $\Gamma^{-1}[G]$ and the upward closure of $\Gamma[G]$ are equal. In particular, the upward closure of $\Gamma[G]$ is a $\dP$-generic filter.
\end{lemma}

\begin{proof}
    First suppose that $q\in\Gamma^{-1}[G]$, i.e. $\Gamma(q)\in G$. Then $\Gamma(\Gamma(q))\leq q$ by assumption and $\Gamma(\Gamma(q))\in\Gamma[G]$, showing that $q$ belongs to the upward closure of $\Gamma[G]$.

    Now suppose that $q$ belongs to the upward closure of $\Gamma[G]$, i.e. $\Gamma(p)\leq q$ for $p\in G$. Then $p\leq\Gamma(\Gamma(p))\leq\Gamma(q)$, showing that $\Gamma(q)\in G$ and thus $q\in\Gamma^{-1}[G]$.

    To see the ``in particular'' part, note that, since $\Gamma\colon\dP\to\dP$ is a dense embedding, $\Gamma^{-1}[G]$ is a $\dP$-generic filter.
\end{proof}

\begin{lemma}[Homogeneity]{\label{Lemma: Almost Homogeneity}}
    Suppose that $\dP$ is a poset such that for every $p,q\in\dP$ there is a dense set $D\subseteq\dP$, a condition $p'\leq p$ with $p'\in D$ and a self-inverse automorphism modulo equivalence $\Gamma\colon D\to D$ such that $\Gamma(p')$ and $q$ are compatible. Then whenever $x_0,\dots,x_{n-1}$ are sets, $\phi$ is a statement in the forcing language and $p\in\dP$, we have
    $$p\Vdash\phi(\check{x}_0,\dots,\check{x}_{n-1})\text{ if and only if }\one_{\dP}\Vdash\phi(\check{x}_0,\dots,\check{x}_{n-1})$$
\end{lemma}

\begin{proof}
    Clearly, the condition $\one_{\dP}\Vdash\phi(\check{x}_0,\dots,\check{x}_{n-1})$ is sufficient, let us show that it is also necessary. To this end, suppose that $p\Vdash\phi(\check{x}_0,\dots,\check{x}_{n-1})$ and assume toward a contradiction that $\one_{\dP}\not\Vdash\phi(\check{x}_0,\dots,\check{x}_{n-1})$, i.e. there is some $q\in\dP$ such that $q\Vdash\neg\phi(\check{x}_0,\dots,\check{x}_{n-1})$. By our assumption there is a dense set $D\subseteq\dP$, a condition $p'\leq p$ with $p'\in D$ and a self-inverse automorphism modulo equivalence $\Gamma\colon D\to D$ such that $\Gamma(p')$ and $q$ are compatible. First of all, we have
    $$p'\Vdash_{\dP}\phi(\check{x}_0,\dots,\check{x}_{n-1})$$
    Since $D$ is dense, this implies
    $$p'\Vdash_D\phi(\check{x}_0,\dots,\check{x}_{n-1})$$
    Since $\Gamma$ is a dense embedding from $D$ to $D$ by Lemma \ref{Lemma: Self-Inverse Auto mod Equiv is DEmb}, this implies that
    $$\Gamma(p')\Vdash_D\phi(\check{x}_0,\dots,\check{x}_{n-1})$$
    which, again by density of $D$, implies
    $$\Gamma(p')\Vdash_{\dP}\phi(\check{x}_0,\dots,\check{x}_{n-1})$$
    but this  contradicts our assumption of compatibility of $\Gamma(p')$ and $q$.
\end{proof}

\begin{definition}
    Let $\dP$ be a forcing poset. We say that $\dP$ is \emph{homogeneous modulo equivalence} if for all $p,q\in\dP$ there is a dense set $D\subseteq\dP$, a condition $p'\leq p$ with $p'\in D$ and a self-inverse automorphism modulo equivalence $\Gamma\colon D\to D$ such that $\Gamma(p')$  and $q$ are compatible.
\end{definition}

Using Lemma \ref{Lemma: Almost Homogeneity} and classical arguments, one sees easily:

\begin{lemma}\label{lemma:HODincluded}
    Suppose that $\dP$ is an ordinal-definable forcing that is homogeneous modulo equivalence. Then, $\one_{\dP}\Vdash\HOD^{V[\dot{G}]}\subseteq\HOD^V$.\qed
\end{lemma}

A poset that will be relevant in the main construction of this paper is the forcing shooting a club through a given stationary set $S$:

\begin{definition}\label{def: clubshooting}
    Let $\kappa$ be a regular cardinal and $S \s \kappa$ a stationary set. We denote by $\mathbb{CS}(S)$ the forcing whose conditions are closed bounded subsets $c\s \kappa$ with a maximal element $\max(c)$ such that $c\s  S$. The order of $\mathbb{CS}(S)$ is reversed end-extension.
\end{definition}

Assuming $\kappa^{<\kappa}=\kappa$ it follows that $|\mathbb{CS}(S)|=\kappa$ and so the poset has the $\kappa^+$-cc. However, $\mathbb{CS}(S)$ may collapse cardinals less than or equal to $\kappa$ unless the stationary set $S$ is taken with care.  The stationary sets $S$ for which $\mathbb{CS}(S)$ behaves as desired were isolated by Abraham--Shelah in \cite{AbrahamShelahForcingClub}:

\begin{definition}\label{def: fat stationary}
    A stationary set $S\s \kappa$ is called \emph{fat} if for every club $C\s \kappa$  the intersection $C\cap S$ contains a closed set of  order-type $\beta$ for all $\beta<\kappa$.   
\end{definition}

Their main result reads as follows:

\begin{fact}[Abraham--Shelah]{\label{Fact: Abraham-Shelah}}
    Let $\kappa$ be a strongly inaccessible cardinal and  $S\s \kappa$ a fat stationary set. Then, $\mathbb{CS}(S)$ is $\kappa$-distributive; namely, forcing with $\mathbb{CS}(S)$ does not add new sequences of ordinals  of length ${<}\kappa$.\qed. 
\end{fact}

The main family of  posets utilized in this manuscript are  of \emph{Prikry-type}:
\begin{definition}\label{GitikPrikrydef}
    A triple $\langle \mathbb{P},\leq,\leq^*\rangle$ is called a \emph{Prikry-type forcing} if
    \begin{enumerate}
        \item  $\leq^*\subseteq \leq$,
        \item $\leq^*$ has the \emph{Prikry property}; to wit, for every $p\in {P}$ and a statement $\varphi$ in the forcing language of $\langle \mathbb{P},\leq\rangle$ there is  $q\leq^* p$ such that $$\text{$q\forces_{\mathbb{P}}\varphi$ or $q\forces_{\mathbb{P}}\neg\varphi.$}$$
    \end{enumerate}
\end{definition}

Any poset is of Prikry-type just letting setting $\leq^* := \leq$. The Prikry-type forcings we shall be interested in  will satisfy $\leq^* \subsetneq \leq$ and the direct extension order $\leq^*$ will exhibit better closure properties than the regular order $\leq.$

\smallskip

In \cite{GitikChanging}, Gitik introduced Easton-supported iterations of Prikry-type forcings. In this paper, we consider a modification due to Ben-Neria and Unger \cite{BenUng}, where the iteration has non-stationary support. This modification allows for finer control over the types of ultrafilters available in the generic extension (see e.g., \cite{HP, GitikKaplan,GitikKaplan2,HP2}).
\begin{definition}[Non-stationary support iteration, \cite{BenUng}]\label{Gitikiteration}
Let $\varrho$ be an ordinal. We define an iteration $\langle\mathbb{P}_\alpha, \dot{\mathbb{Q}}_\beta\mid \beta<\alpha\leq \varrho\rangle$ recursively as follows:  For  $\alpha\leq \varrho$ members of $\mathbb{P}_\alpha$ are sequences $p=\langle \dot{p}_\beta\mid \beta\in \supp(p)\rangle$ such that:
\begin{enumerate}[label=(\alph*)]
    \item $\supp(p)\s \alpha$ has \emph{non-stationary support}: namely, $\supp(p)\cap \beta$ is \emph{non-stationary} in $\beta$ for all inaccessible cardinals $\beta\leq \alpha$.
    \item for every $\beta\in \supp(p)$, $p\restriction\beta:=\langle \dot{p}_\beta\mid \beta\in \supp(p)\cap\beta\rangle\in \mathbb{P}_\beta$ and      $$p\restriction\beta\forces_{\mathbb{P}_\beta}\dot{p}_\beta\in \dot{\mathbb{Q}}_\beta.$$
    \item For every $\beta\in \supp(p)$, the weakest condition of $\mathbb{P}_\beta$ forces that:
    \begin{itemize}[label=-]
        \item $\langle \dot{\mathbb{Q}}_\beta,\leq_\beta, \leq_\beta^*\rangle$  is a Prikry-type forcing notion.
        \item $|\dot{\mathbb{Q}}_\beta|\leq 2^\beta$.
        \item There exists a $\leq^*_\alpha$-dense subset $\dot{\mathbb{D}}_{\alpha}\subseteq \dot{\mathbb{Q}}_\alpha$ such that $\langle \dot{\mathbb{D}}_\alpha, \leq^*_\alpha \rangle$ is $\alpha$-closed.\footnote{Recall that $\beta$-closure stands for the following property: For every $\leq^*_\beta$-decreasing sequence $\langle r_\gamma\mid \gamma<\bar{\gamma}\rangle$ of members of $\mathbb{Q}_\beta$ (with $\bar{\gamma}<\beta$) there is $r\in\mathbb{Q}_\beta$ such that $r\leq^*_\beta r_\gamma.$ }
    \end{itemize}

\end{enumerate}
Let $p=\langle \dot{p}_\beta\mid \beta\in \supp(p)\rangle$ and $q=\langle \dot{q}_\beta\mid \beta\in \supp(q)\rangle$ be elements of $\mathbb{P}_\alpha$.

We write $p\leq_\alpha q$  if and only if the following hold:
\begin{enumerate}
    \item $\supp(p)\supseteq \supp(q)$,
    \item $p\restriction\beta\forces_{\mathbb{P}_\beta}\dot{p}_\beta\leq_\beta \dot{q}_\beta$, for every $\beta\in \supp(q)$;
    \item there is $b\subseteq \supp(q)$, finite, such that for all $\beta\in \supp(q)\setminus b$
    $$p\restriction\beta\forces_{\mathbb{P}_\beta} \dot{p}_\beta\leq^*_\beta \dot{q}_\beta.$$
\end{enumerate}
The Prikry ordering $\leq^*_\alpha$ is defined by the case in which $b$ above is empty.
\end{definition}
\begin{remark}
   About (a): While members of $p\in \mathbb{P}_\varrho$ are names for the various posets involved note that $\supp(p)$ is assumed to be a set in the ground model. About (c): 
  the requirement concerning the $\leq^*$-dense sets $\mathbb{D}_\alpha$ replaces the standard assumption that $\langle \mathbb{Q}_\alpha, \leq^*_{\alpha}\rangle$ is $\alpha$-closed.
\end{remark}

A fundamental theorem regarding the aforementioned iterations is the Prikry property. The result for Easton-supported iterations was proved by Gitik \cite[\S1]{GitikChanging}; the corresponding one for non-stationary-supported iterations is due to Ben-Neria and Unger \cite[\S2]{BenUng}.

\begin{lemma}[{\cite[Lemma~2.2]{BenUng}}]\label{lemma: prikry lemma}
For each $\alpha\leq \varrho$,
$\langle \mathbb{P}_\alpha,\leq_\alpha, \leq^*_\alpha\rangle$ has the Prikry property.\qed
\end{lemma}

The analysis of nonstationary support iterations heavily relies on the following Fusion Lemma.

\begin{lemma}[Fusion Lemma, {\cite[Lemma 2.3]{BenUng}}]\label{lemma: fusion lemma}
  Suppose that $\alpha\leq \kappa$ is an inaccessible cardinal and $d\colon \alpha \to V$ is a function such that, for each $\beta<\alpha$, $d(\beta)$ is a $\mathbb{P}_{\beta+1}$-name forced by the trivial condition of  ${\mathbb{P}_{\beta+1}}$ to be a dense open subset of $\langle \mathbb{P}_{\alpha}/\dot{G}_{\beta+1},\leq^*\rangle$. Then, given $p\in \mathbb{P}_{\alpha}$ and $\nu<\alpha$ there is $p^*\leq^* p$ and a club $C\s \alpha$ such that $p^*\restriction \nu=p\restriction \nu$ and for each $\beta\in C$, $p^*\restriction \beta+1\forces \dot{p}^*\setminus \beta+1\in d(\beta)$.\qed
\end{lemma}

\begin{lemma}[{\cite[Lemma 2.5 and Corollary 2.6]{gitik2024fresh}}]\label{Lem: EvaluatingNamesOfOrdinals}
    Assume {\rm GCH}. Let $\alpha$ be an inaccessible cardinal and let $\mathbb{P}_\alpha$ be an iteration of Prikry-type forcing satisfying the hypotheses of Definition \ref{Gitikiteration}. 
    \begin{enumerate}
        \item Assume that $\dot{\eta}$ is a $\mathbb{P}_\alpha$-name for an ordinal. Then for every $p\in \mathbb{P}_\alpha$ there exists $q\leq^* p$ and a set $X\in V$, $|X|<\alpha$ such that $q^*\Vdash \dot{\eta}\in \check{X}$.
        \item $\mathbb{P}_\alpha$ preserves the inaccessibility of $\alpha$.\qed
    \end{enumerate}
\end{lemma}

\begin{cor}[{\cite[Theorem 6.5]{gitik2024fresh}}]\label{corollary: cover property}
    Assume {\rm GCH}. Let $\alpha$ be an inaccessible cardinal. Let $\mathbb{P}_\alpha$ be a nonstationary support iteration of Prikry-type forcing satisfying the hypotheses of Definition \ref{Gitikiteration}. Then $\mathbb{P}_\alpha$ satisfies the $\alpha$ and $\alpha^+$-cover property.\footnote{We say that a poset $\mathbb{P}$ has the $\kappa$-cover property if $V\subseteq V[G]$ has the $\kappa$-cover property for every $V$-generic filter $G\subseteq\mathbb{P}$. That is, for every such $G$ and every set of ordinals $A\in V[G]$ with $|A|^{V[G]}<\kappa$, there is a set of ordinals $B\in V$ such that $A\subseteq B$ and $|B|^{V[G]}<\kappa$.} 
\end{cor}

\begin{proof}
    We begin with the $\alpha$-cover property of $\mathbb{P}_\alpha$. Assume that $A\in V[G]$ is a set of size strictly below $\alpha$. Let $\dot{\eta}$ be a $\mathbb{P}_\alpha$-name for an ordinal below $\alpha$ which bounds the cardinality of $\dot{A}$ in $V[G]$. By Lemma \ref{Lem: EvaluatingNamesOfOrdinals}, there exist $\eta^*< \alpha$ and a condition $p\in G$ such that $p\Vdash |\dot{A}|<\eta^*$. Work in $V[G\restriction (\eta^*+1)]$ and let 
    $$\mathbb{R}:= \{ r\in  \mathbb{P}_\alpha\setminus (\eta^*+1) \colon \forall \beta\in \supp(r), r\restriction \beta\Vdash r(\beta)\in \dot{\mathbb{D}}_{\beta} \}.$$
    Using Lemma \ref{Lem: EvaluatingNamesOfOrdinals} again, build a $\leq^*$-decreasing sequence of conditions $\langle r_i \colon i\leq \eta^* \rangle\subseteq \mathbb{R}$ below $p\setminus (\eta^*+1)$ and sets $\langle X_i \colon i<\eta \rangle $ such that $r_i$ forces that for every $i<\eta^*$, the $i$-th element of $\dot{A}$ (if exists) belongs to $X_i$, and $|X_i|<\alpha$. Note that we have sufficient closure in $\mathbb{R}$ at limit stages, including at the very final stage. By density, we can find a lower bound for such a sequence in $G\setminus (\eta^*+1)$, and this lower bound forces that $\dot{A}$ is covered by $X^*=\bigcup_{i<\eta^*} X_i$. In particular, there exists a condition $p^*\in G$ such that
    $$ p^*\restriction (\eta^*+1) \Vdash \exists X^*, |X^*|<\alpha \land   \left( p^*\setminus (\eta^*+1) \Vdash \dot{A}\subseteq X^* \right). $$
    Since $p^*\in G$ and $\mathbb{P}_{\eta^*+1}$ is small relative to $\alpha$, there exists $B\in V$ with $|B|<\alpha$ such that $p^*\Vdash \dot{A}\subseteq B$, as desired.

    Next, let us prove the $\alpha^+$-cover property of $\mathbb{P}_\alpha$. Let $A\in V[G]$ be a set of size $\alpha$. For simplicity, we identify $\dot{A}$ with a name for an enumeration of $\dot{A}$ of order type $\alpha$. In particular, for every $\beta<\alpha$, let $\dot{A}(\beta)$ denote the $\beta$-th element in the enumeration. For every $\beta<\alpha$, let 
    $$ e(\alpha) = \{ r\in \mathbb{P}_\alpha\setminus (\beta+1) \colon \exists X \left( |X|<\alpha \ \land r\Vdash \dot{A}\restriction \beta\subseteq  \check{X} \right) \}. $$
    Each set $e(\beta)$ is forced to be a $\leq^*$-dense open subset of $\mathbb{P}_\alpha\setminus (\beta+1)$, by Lemma \ref{Lem: EvaluatingNamesOfOrdinals} and since $\mathbb{P}\setminus (\beta+1)$ has a $\leq^*$-dense subset which is $(\beta+1)$-closed. By the Fusion Lemma (Lemma \ref{lemma: fusion lemma}), there exists $p\in G$ and a club $C\subseteq \alpha$ such that, for every $\beta\in C$, there exists a $\mathbb{P}_{\beta+1}$ name for a set of ordinals $\dot{X}_\beta$ of size strictly below $\alpha$, such that $p\Vdash \dot{A}\restriction \beta \subseteq \dot{X}_\beta$. Since $\mathbb{P}_{\beta+1}$ is a small forcing, we can find in $V$ a sequence $\langle X^*_\beta \colon \beta\in C \rangle$ such that each set $X^*_\beta$ has size strictly below $\alpha$, and $p\Vdash \dot{X}_\beta\subseteq \check{X}^*_\beta$. It thus follows that $B:= \bigcup_{\beta\in C} X^*_\beta$ is a set in $V$, $|B| = \alpha$, and $A\subseteq B$.
\end{proof}

\begin{lemma}\label{lemma: preserving stationaries}
    Assume {\rm GCH}. Let $\alpha$ be an inaccessible cardinal. Let $\mathbb{P}_\alpha$ be a nonstationary support iteration of Prikry-type forcing satisfying the hypotheses of Definition \ref{Gitikiteration}. Let $G\subseteq \mathbb{P}_\alpha$ be generic over $V$. Then  every club in $\alpha$ in $V[G]$ contains a club in $\alpha$ which belongs to $V$. In particular, forcing with $\mathbb{P}_\alpha$ preserves stationary subsets of $\alpha$.
\end{lemma}
\begin{proof}
Let $\dot{C}$ be a $\mathbb{P}_\alpha$-name for a club, and let $\dot{f}$ be a $\mathbb{P}_\alpha$-name for the increasing enumeration of $\dot{C}$. For every $\beta<\alpha$, let $d(\beta)$ be the $\mathbb{P}_{\beta+1}$-name for the set
$$ \{ q\in \mathbb{P}_\alpha \setminus (\beta+1 ) \mid \exists \gamma<\alpha ( q\Vdash \dot{f}(\beta)<\gamma ) \}. $$
Each such set $d(\beta)$ is forced to be a $\leq^*$-dense open subset of $\mathbb{P}_\alpha\setminus (\beta+1)$ by Lemma \ref{Lem: EvaluatingNamesOfOrdinals}. Thus, by the Fusion Lemma (Lemma \ref{lemma: fusion lemma}), there exists $p\in G$ and a club $D\subseteq \alpha$ such that, for every $\beta\in D$, 
$$p\restriction (\beta+1)\Vdash \exists \gamma(\beta)<\alpha  \left( p\setminus (\beta+1)\Vdash \dot{f}(\beta)< \gamma(\beta) \right).$$
Let $F\colon D\to \alpha$ in $V$ be a function such that, for every $\beta\in D$, $p\restriction (\beta+1)\Vdash \dot{\gamma}(\beta)< F(\beta)$ (such $F$ exists since $|\mathbb{P}_{\beta+1}|<\alpha$ for every $\beta<\alpha$). In $V$, let
$$C^*=\left\{ \beta\in\text{acc}(D)\mid F[D\cap\beta]\subseteq\beta\right\}.$$
Then $p\Vdash C^*\subseteq \dot{C}$.
\end{proof}

\begin{lemma}[{\cite[Corollary 2.6]{BenUng}}]\label{Lemma: Nonstatiterationspreservecardinals}
    Assume {\rm GCH}. Let $\alpha$ be an inaccessible cardinal. Let $\mathbb{P}_\alpha$ be a nonstationary support iteration of Prikry-type forcing satisfying the hypotheses of Definition \ref{Gitikiteration}. Assume also that for every $\beta<\alpha$, $\Vdash_{\mathbb{P}_\beta} \dot{\mathbb{Q}}_\beta \text{ preserves cardinals.}$
    Then $\mathbb{P}_\alpha$ preserves cardinals.\qed
\end{lemma}

\subsection{Coherent Sequences of Measures}

Our main iteration will be built using Mitchell \emph{coherent sequences of measures} \cite{MitIter}. In this subsection, we will collect the main definitions and prove that, from the consistency of a supercompact cardinal, we can obtain a coherent sequence of measures satisfying our desired properties.

\begin{definition}
    Let $\kappa$ be a cardinal and let $\mathcal{U}=\langle\mathcal{U}_{\alpha,\beta}\mid\alpha\leq\kappa,\beta<o^{\mathcal{U}}(\alpha)\rangle$ be a sequence. We say that $\mathcal{U}$ is a \emph{coherent sequence of measures} if the following holds:
    \begin{enumerate}
        \item For every $\alpha\leq\kappa$ and every $\beta<o^{\mathcal{U}}(\alpha)$, $\mathcal{U}_{\alpha,\beta}$ is a measure over $\alpha$;\footnote{I.e.,  an $\alpha$-complete normal ultrafilter on $\alpha$.}
        \item For every $\alpha\leq\kappa$ and every $\beta<o^{\mathcal{U}}(\alpha)$, letting $j_{\mathcal{U}_{\alpha,\beta}}$ be the ultrapower embedding associated to $\mathcal{U}_{\alpha,\beta}$, we have
        $$j_{\mathcal{U}_{\alpha,\beta}}(\mathcal{U})\restriction(\alpha+1,0)=\mathcal{U}\restriction(\alpha,\beta)$$
    \end{enumerate}
\end{definition}

In more detail, point (2) of the above definition states that the following holds in $M_{\mathcal{U}_{\alpha,\beta}}$:

\begin{enumerate}
    \item $o^{j_{\mathcal{U}_{\alpha,\beta}}(\mathcal{U})}(\alpha)=\beta$;
    \item For every $\gamma<\beta$, $j_{\mathcal{U}_{\alpha,\beta}}(\mathcal{U})_{\alpha,\gamma}=\mathcal{U}_{\alpha,\gamma}$.
\end{enumerate}

To simplify the notation, we will denote $j_{\mathcal{U}_{\alpha,\beta}}$ simply by $j_{\alpha,\beta}$ and $M_{\mathcal{U}_{\alpha,\beta}}$ by $M_{\alpha,\beta}$ if the coherent sequence $\mathcal{U}$ is clear from context (from Section \ref{sec: the iteration} onwards, we will fix a single coherent sequence).

\begin{definition}
    Let $\kappa$ be a cardinal and let $\mathcal{U}=\langle\mathcal{U}_{\alpha,\beta}\mid\alpha\leq\kappa,\beta<o^{\mathcal{U}}(\alpha)\rangle$ be a coherent sequence. We let $\dom_1(\mathcal{U})$ consist of all $\alpha$ such that $o^{\mathcal{U}}(\alpha)\neq 0$.
\end{definition}

We will now show that, from a supercompact cardinal, we can derive a nice coherent sequence of measures:

\begin{lemma}\label{Lemma: Construction of the coherent sequence}
    Let $\kappa$ be a supercompact cardinal and assume $V={\rm{HOD}}$. There exists a coherent sequence of measures $$\mathcal{U} = \langle \mathcal{U}_{\alpha,\beta} \colon \alpha\leq \kappa, \beta<o^{\mathcal{U}}(\alpha) \rangle$$  with the following properties:
    \begin{enumerate}
        \item For every $\alpha\leq \kappa$, $o^{\mathcal{U}}(\alpha)\leq \alpha^+$, and  $o^{\mathcal{U}}(\kappa) = \kappa^+$.
        \item Each measure $\mathcal{U}_{\alpha,\beta}$ with $o^{\mathcal{U}}(\alpha)\neq 0$ concentrates on the set of measurable cardinals $\alpha^*<\kappa$ such that $\ell ''\alpha^*\subseteq V_{\alpha^*}$.
        \item For every $(\alpha, \beta)\in \dom(\mathcal{U})$, $$j_{\mathcal{U}_{\alpha,\beta}}(\prec)\restriction H_{\alpha^+} = \prec \restriction H_{\alpha^+},$$ where $\prec$ denotes the canonical well order of {\rm HOD} (see p.\pageref{wellorderofHOD}).
        \item There are arbitrarily large cardinals $\lambda>\kappa$ such that every fine, normal measure $\mathcal{V}$ on $\mathcal{P}_{\kappa}(\lambda)$ satisfies $j_{\mathcal{V}}(\mathcal{U})\restriction \kappa+1 = \mathcal{U} $.
    \end{enumerate}
\end{lemma}

\begin{proof}
    Let us construct a coherent sequence of measures $\mathcal{U}$ by induction on pairs $(\alpha,\beta)$ as follows. Assume $\alpha\leq \kappa$ is measurable and $\mathcal{U}\restriction (\alpha,\beta)$ has been constructed for some ordinal $\beta\leq \alpha^+$. If $\beta=\alpha^+$, set $o^{\mathcal{U}}(\alpha) := \alpha^+$. Else, let $\mathcal{U}_{\alpha,\beta}$ be (if it exists) the $\prec$-least normal measure $W$ on $\alpha$ which concentrates on measurable cardinals $\alpha^*<\alpha$ with $\ell''\alpha^*\subseteq V_{\alpha^*}$, and satisfies $$j_{W}( \mathcal{U}\restriction \alpha )(\alpha) = \langle \mathcal{U}_{\alpha,\gamma} \mid \gamma<\beta \rangle$$
    and
    $$ j_{W}(\prec)\restriction H_{\alpha^+} = \prec \restriction H_{\alpha^+}.$$
     If such measure  $W$ does not exist, we stipulate $o^{\mathcal{U}}(\alpha) := \beta$.

     \smallskip
    
    Note that, by construction,  $o^{\mathcal{U}}(\alpha)\leq \alpha^+$ for every $\alpha<\kappa$.

    It remains to prove that $\mathcal{U}$ is as desired. It is not hard to verify that $\mathcal{U}$ is coherent. Let us claim that $o^{\mathcal{U}}(\kappa) = \kappa^+$. Assume otherwise. Fix any Beth-fixed point $\lambda>2^{2^\kappa}$ so that $\lambda$ is a closure point of the map $\alpha\mapsto \sup_{x\in V_\alpha}\text{rank}_{\prec}(x)$. Let $\mathcal{V}$ be a supercompactness measure on $\mathcal{P}_{\kappa}(\lambda)$. Let $j_{\mathcal{V}}\colon V\to M$ be the associated elementary embedding, with critical point $\kappa$, such that ${}^\lambda M \subseteq M$. In particular, since $(V_{\lambda})^M = V_{\lambda}$ and by the definition of $\prec$, we obtain that $j_{\mathcal{V}}(\prec)\restriction V_{\lambda} = \prec \restriction{V_{\lambda}}$.\footnote{Indeed, for $x_0,x_1\in V_\lambda$, the $j_{\mathcal{V}}(\prec)$ relation between $x_0,x_1$ depends only on $V_{\lambda}$:  Recall that the order relies on the lexicographic order between tuples $\langle \theta^{x_i}, n^{x_i}, \alpha^{x_i}_0,\ldots, \alpha^{x_i}_{m_{x_i}}\rangle $ (see p.\pageref{wellorderofHOD}), and each such tuple belongs to $V_\lambda$ by our choice of $\lambda$. The tuples associated to $x_i$ are computed the same over $V$ and $M$ since $(V_\lambda)^M = V_\lambda$.} In particular, the measures $$\langle \mathcal{U}_{\kappa, \beta} \mid \beta < o^{\mathcal{U}}(\kappa)\rangle$$ are chosen in the exact same way over $V$ and $M$, and are definable in the latter model. Next, let $W = \{ X\subseteq \kappa \mid \kappa\in j_{\mathcal{V}}(X) \}$ be the normal measure derived from $\mathcal{V}$. Note that $W$ concentrates on $\{\alpha^*<\kappa \mid \ell'' \alpha^*\subseteq V_{\alpha^*} \}$ since $j_{\mathcal{V}}(\ell)\restriction \kappa = \ell$. Let $k\colon M_W\to M_{\mathcal{V}}$ be the standard  embedding mapping $[f]_W$ to $j_{\mathcal{V}}(f)(\kappa)$. We have that $\text{crit}(k)= (\kappa^{++})^{M_{W}}$. In particular, 
    $$ j_{W}(\prec)\restriction H_{\kappa^+} =j_{\mathcal{V}}(\prec)\restriction H_{\kappa^+} = \prec\restriction H_{\kappa^+}.$$
    Furthermore, by definability of $\mathcal{U}$, we have that $\langle \mathcal{U}_{\kappa, \beta} \mid  \beta<o^{\mathcal{U}}(\kappa) \rangle$ is in $\text{range}(k)$. In particular, since $\text{crit}(k)> \kappa^+$, we have $\langle \mathcal{U}_{\kappa, \beta}\mid  \beta<o^{\mathcal{U}}(\kappa) \rangle$ is in  $M_{W}$, and $j_{W}(\mathcal{U})\restriction \kappa+1 = \mathcal{U}$. Thus we could not have had $o^{\mathcal{U}}(\kappa) <\kappa^+$,
    since otherwise $W$ was eligible for the role of $\mathcal{U}_{\kappa, o^{\mathcal{U}}(\kappa)}$.

\smallskip

The above proof also establishes that $j_{\mathcal{V}}(\mathcal{U})\restriction \kappa+1 = \mathcal{U}$, again by the facts that $(V_\lambda)^M = V_{\lambda}$ and the orders $j_{\mathcal{V}}(\prec)$ and $\prec$ agree on $V_{\lambda}$. 
\end{proof}
\begin{remark}\label{rem:flexibility}
    From the proof of Lemma \ref{Lemma: Construction of the coherent sequence}, it is easy to see that there is a lot of flexibility regarding the choice of what the measures should concentrate on. The only requirement is that this should hold for $\kappa$ in the target model of the chosen $\lambda$-supercompact embedding.\footnote{This hypothesis can be realized as follows: Let $\kappa$ be a $C^{(2)}$-extendible cardinal \cite{Bag}. By a result of Bagaria--Goldberg \cite{BagGol}, for every $\lambda>\kappa$ with $V_\lambda\prec_{\Sigma_2} V$ there is a $\lambda$-supercompact embedding $j\colon V\to M$ with $M\models ``V_\lambda\prec_{\Sigma_2} V$''. Since $\kappa$ is supercompact and this is a $\Pi_2$-definable property, $V_\lambda$ thinks that $\kappa$ is supercompact and then, since $\lambda$ is $\Sigma_2$-elementary inside $M$, the same is thought by $M$ as well.} Therefore, with the same large cardinal assumption, we could have that the measures concentrate on cardinals with an arbitrary (but bounded) degree of supercompactness. With a stronger large cardinal assumption, we may even have the measures concentrate on fully supercompact cardinals.
\end{remark}

The forcing we shall introduce in the next section requires canonical functions to be defined (intuitively, canonical functions will allow us to reflect the order of measurable cardinals $\alpha\in \dom_1(\mathcal{U})$ that go past $\alpha$). Bearing this issue in mind  we fix a sequence of functions
$$\vec{f}=\langle f^\delta_\tau\mid \delta\in \mathrm{Reg}\cap \kappa^+\,\wedge\, \tau\leq \delta^+\rangle$$
which is defined as follows: For each regular cardinal $\delta\leq\alpha$ and $\tau\leq\delta^+$, $f_{\tau}^{\delta}\colon \delta\to \delta$ will denote a \emph{canonical function for $\tau$}; namely,   for $\tau<\delta$, $f_{\tau}^{\delta}$ is the constant function with value $\tau$, $f_{\delta}^{\delta}$ is the identity, $f^{\delta}_{\delta^+}$ is the cardinal successor function $\xi\mapsto|\xi|^+$. To ensure that the forcing is well-defined, for $\tau\in (\delta,\delta^+)$, we  require that   $f^{\delta}_{\tau}$ is the $\prec$-minimal function in the $=_{\rm NS}$-class of the $\tau$-th canonical function. This  ensures the following coherency property:
\begin{equation}\label{eq: coherency}
   j_{\mathcal{U}_{\alpha,\sigma}}(\vec{f})_{\tau}^{\alpha}=f_{\tau}^{\alpha},\;\;\text{for every $\sigma\leq o^{\mathcal{U}}(\alpha)$ and $\tau\leq \alpha^+$.}\footnote{Recall that, by our choice, $j_{\mathcal{U}_{\alpha,\sigma}}(\prec)$ and $\prec$ coincide on $H_{\alpha^+}$} 
\end{equation}

\section{The main iteration}\label{sec: the iteration}

The following are the setup assumptions for the rest of this section:

\begin{setup}\label{setup}
   Assume GCH holds, $V = g{\rm HOD}$,  $\kappa$ is a supercompact cardinal, $\ell\colon \kappa\to V_\kappa$ is a Laver function for $\kappa$,  $$\mathcal{U}=\langle \mathcal{U}_{\alpha,\tau}\mid \alpha\leq \kappa\,\wedge\, \tau<o^{\mathcal{U}}(\alpha)\rangle$$ is  a coherent sequence of normal measures on cardinals  $\alpha\leq \kappa$ satisfying the hypotheses of Lemma \ref{Lemma: Construction of the coherent sequence}, and $\vec{f}$ is a sequence of functions as in \eqref{eq: coherency}.  To simplify notations we  write $j_{\alpha,\beta}:=j_{\mathcal{U}_{\alpha,\beta}}.$
\end{setup}

Our main forcing  $\mathbb{P} = \mathbb{P}_\kappa *\dot{\mathbb{Q}}_\kappa$ will be a $\dom_1(\mathcal{U})$-spaced nonstationary support iteration of Prikry type forcings of length $(\kappa+1)$. The construction of the iteration is by induction on $\alpha\leq \kappa$. Let us describe the forcings used at each stage  $\alpha\leq \kappa$. The picture is the following: 
\begin{definition}[The iteration]\label{definition of the iteration}
We shall denote by $\mathbb{P}_{\kappa+1}$ the $\dom_1(\mathcal{U})$-spaced nonstationary support iteration $\langle \mathbb{P}_\alpha, \dot{\mathbb{Q}}_\beta\mid\beta<\alpha\leq \kappa\rangle$ which for each $\alpha\leq \kappa$ is defined as follows:

    \smallskip

    $\br$ If $\alpha\notin \dom_1(\mathcal{U})$ then $\one\forces_{\mathbb{P}_\alpha}``\dot{\mathbb{Q}}_\alpha$ is the trivial forcing".

    \smallskip

     $\br$ If $\alpha\in \dom_1(\mathcal{U})$ but $o^{\mathcal{U}}(\alpha)<\alpha^+$ then $\one\forces_{\mathbb{P}_\alpha}``\dot{\mathbb{Q}}_\alpha=\mathbb{Q}^{o^{\mathcal{U}}(\alpha)}_\alpha$" (see \S\ref{sec: basic building block}).

     \smallskip

     $\br$ If $\alpha\in \dom_1(\mathcal{U})$ and $o^{\mathcal{U}}(\alpha)=\alpha^+$ then $\one\forces_{\mathbb{P}_\alpha}`` \dot{\mathbb{Q}}_\alpha=\langle \mathbb{R}_\alpha,\leq_{\dR_{\alpha}},\leq^*_{\dR_{\alpha}}\rangle$".

     \smallskip

     $\br$ If $\alpha=\kappa$ then $\one\forces_{\mathbb{P}_\alpha}`` \dot{\mathbb{Q}}_\alpha=\langle \mathbb{R}_{\alpha},\leq_{\dR_{\alpha}},\leq^*_{\dR_{\alpha}}\rangle$".
\end{definition}
If $o^{\mathcal{U}}(\alpha)<\alpha^+$  then we force with the adaptation of Gitik's poset $\mathbb{Q}^{o^{\mathcal{U}(\alpha)}}_\alpha$ from \cite[\S3]{GitikChanging} presented in \S\ref{sec: basic building block}. 
In case $o^{\mathcal{U}}(\alpha)=\alpha^+$,   $\langle {\mathbb{R}}_\alpha,\leq_{\dR_{\alpha}},\leq^*_{\dR_{\alpha}}\rangle$ is going to be a Prikry-type poset such that $\langle {\mathbb{R}}_\alpha,\leq_{\dR_{\alpha}}\rangle$ is forcing  equivalent to $\mathbb{C}_\alpha$\footnote{In particular, forcing with $\langle \mathbb{R}_\alpha,\leq_{\dR_{\alpha}}\rangle$ will add a club subset of $\alpha$ through the set of $V$-regular cardinals. The reason we use $\mathbb{C}_\alpha$, the $\alpha$-support product of $\alpha^+$ copies of club shooting $\mathbb{CS}_\alpha$, rather than a single copy of club shooting forcing, is related to the verification of homogeneity of $\mathbb{P}_\kappa$. See Remark~\ref{keyremark}.} and $\langle \mathbb{R}_\alpha,\leq^*_{\dR_{\alpha}}\rangle$ has a $\leq^*_{\mathbb{R}_\alpha}$-dense subset which is $\alpha$-closed (see Lemma~\ref{lemma: adding clubs}). We introduce $\dR_{\alpha}$ in Subsection~\ref{Subsection: Prikry Club-Shooting}. 

\subsection{The basic building block}\label{sec: basic building block}

We now define a variant of Gitik's forcing $\mathbb{P}(\alpha, o^{\mathcal{U}}(\alpha))$ that works for measurable cardinals $\alpha$ such that $o^{\mathcal{U}}(\alpha)\leq \alpha^+$. In \cite{GitikChanging}, Gitik assumes that $o^{\mathcal{U}}(\alpha)\leq \alpha$ for all $\alpha\leq \kappa$ but we shall show that his construction can be naturally modified to encompass measurables with larger Mitchell order. Our assumptions hereafter are as in Setup \ref{setup}.

Fix $\alpha\leq \kappa$. Assume that $\alpha\in \dom_1(\mathcal{U})$, $\mathbb{P}_\alpha$ has been defined, and $G_\alpha\subseteq \mathbb{P}_\alpha$ is generic over $V$. Assume that for every $\beta\in \dom_1(\mathcal{U})\cap \alpha$ with $o^{\mathcal{U}}(\beta)< \beta^+$, a set $b_{\beta}\subseteq \beta$ has been added in $V[G_\alpha]$, such that $b_\beta\subseteq \beta$ is a closed unbounded subset of $\beta$, and
$$ \text{otp}(b_\beta) = \begin{cases}
    \omega^{  o^{\mathcal{U}}(\beta) } & \text{If }  o^{\mathcal{U}}(\beta)<\beta\\
    \beta & \text{If }  o^{\mathcal{U}}(\beta)\geq \beta
\end{cases}.$$

 For the rest of the section we we work in $V[G_{\alpha}]$.

\begin{definition}\label{def: coherency}
    Let $\tau\leq o^{\mathcal{U}}(\alpha)$. A sequence $\langle \delta_0,\dots, \delta_{n-1}\rangle$  in $ [\alpha]^{<\omega}$ is called $\tau$-coherent if it is increasing and the following hold:
    \begin{enumerate}
        \item $o^{\mathcal{U}}(\delta_i)<f^\alpha _\tau(\delta_i)$ for all $i<n$;
      \item

   {For each $i<n$ let $n_i\leq i$ be the smallest index such that for each $j\in[n_i,i)$, $o^{\mathcal{U}}(\delta_j)<f^{\delta_i}_{o^{\mathcal{U}}(\delta_i)}(\delta_j)$.}\footnote{Note that $n_i$ could be $i$, in which case the set of those $j$'s is empty.} Then, $\textstyle \bigcup_{j\in [n_i,i)}b_{\delta_j}\cup\{\delta_j\}\sqsubseteq b_{\delta_i}.$
   \item {For each $0<i<n$, if $o^{\mathcal{U}}(\delta_{i-1})\geq f^{\delta_{i}}_{o^{\mathcal{U}}(\delta_i)}(\delta_{i-1})$ then $\min(b_{\delta_i})>\delta_{i-1}$.}
    \end{enumerate}
\end{definition}

\begin{definition}
    Let $t=\langle\delta_0,\dots,\delta_{n-1}\rangle$ be a $\tau$-coherent sequence and let $\sigma<\tau$. We define $t\restriction\sigma:=\langle\delta_i,\dots,\delta_{n-1}\rangle$, where $i\leq n$ is the smallest index such that for each $j\in[i,n)$, $o^{\mathcal{U}}(\delta_j)<f^\alpha_{\sigma}(\delta_j)$.
\end{definition}

\begin{definition}
    Let $t=\langle\delta_0,\dots,\delta_{n-1}\rangle$ be a $\tau$-coherent sequence. We say that $t$ is \emph{minimal} if for all $0<i<n$, $o^{\mathcal{U}}(\delta_{i-1})\geq f_{o^{\mathcal{U}(\delta_i)}}^{\delta_i}(\delta_{i-1})$.
\end{definition}

\begin{remark}
    The above definitions provide  a natural extension of the notion of coherency considered by Gitik  \cite{GitikChanging}. In Gitik's setup every $\alpha\in \dom_1(\mathcal{U})$ is assumed to satisfy $o^{\mathcal{U}}(\alpha)<\alpha$. This  amounts to saying that, for each $(\theta,\tau)\in \dom(\mathcal{U})$,
    $f^\theta_\tau$ is the constant function with value $\tau$, which yields  $o^{\mathcal{U}}(\delta_i)<\tau$ in (1)  and $o^{\mathcal{U}}(\delta_j)<o^{\mathcal{U}}(\delta_i)$ in (2) of  Definition~\ref{def: coherency}. Thus, we recover Gitik's notion of $\tau$-coherency.  Similarly, in the definition of $t\restriction\sigma$, $f^\alpha_\sigma(\delta_j)=\sigma$, which recovers Gitik's definition of $t\restriction\sigma$. Finally, our notion of minimal sequence is the natural adaptation of what Ben-Neria–Unger call \emph{maximal sequence}.\footnote{We deemed more appropriate the terminology ``minimal sequence'' in light of Lemma~\ref{lemma: minimal} below.}
\end{remark}

\begin{definition}
       For each $\tau$-coherent sequence $t$ and $\sigma<\tau$ we define
       $\mathcal{U}_{\alpha,\sigma}(t)$ as the collection of those $X\in \mathcal{P}^{V[G_{\alpha}]}(\alpha)$ for which there are $p\in G_{\alpha}$ and a $\mathbb{P}_\alpha$-name $\dot{T}$ such that $p\forces_{\mathbb{P}_\alpha}\langle t\restriction\sigma, \dot{T}\rangle\in \mathbb{Q}^\alpha_\sigma$ such that 
       $$p{}^\smallfrown \langle t\restriction\sigma, \dot{T}\rangle{}^\smallfrown j_{\alpha,\sigma}(p)\setminus (\alpha+1)\forces_{j_{\alpha,\sigma}(\mathbb{P}_\alpha)}\alpha\in j_{\alpha,\sigma}(\dot{X}).$$
\end{definition}

\begin{remark}
  {For each $\tau$-coherent sequence $t$ it makes sense to consider $\mathcal{U}_{\alpha,\sigma}(t\restriction\sigma)$ because $t\restriction\sigma $ is also $\tau$-coherent. Moreover, $\mathcal{U}_{\alpha,\sigma}(t)=\mathcal{U}_{\alpha,\sigma}(t\restriction\sigma)$.}
\end{remark}

\begin{lemma}\label{Lemma: Filter Properties}
Let $\sigma<\tau\leq o^{\mathcal{U}}(\alpha)$ and $t$ be a  $\tau$-coherent sequence. Then: 
\begin{enumerate}
    \item $\mathcal{U}_{\alpha, \sigma}(t)$ is an $\alpha$-complete ultrafilter in $V[G_{\alpha}]$ extending $\mathcal{U}_{\alpha,\sigma}$.
     \item $\{\delta<\alpha\mid f_\sigma(\delta)=o^{\mathcal{U}}(\delta)\}\in \mathcal{U}_{\alpha,\sigma}(t)$.%
    \item
{$\{\nu<\alpha\mid \text{$t^\smallfrown \langle \nu\rangle$ is $\tau$-coherent}\}\in \mathcal{U}_{\alpha,\sigma}(t)$.}
\end{enumerate}
\end{lemma}
\begin{proof}
    (1). This follows from \cite[Proposition 3.3]{BenUng}.

\smallskip

(2). This is because $\{\nu<\alpha\mid f_\sigma(\delta)=o^{\mathcal{U}}(\delta)\}\in \mathcal{U}_{\alpha,\sigma}.$

\smallskip
    {(3). Here is where the canonicity of the functions $f^\delta_{\tau}$'s and moreover the canonicity of the choice of the sequence $\vec{f}$ plays a key role. Let $\dot{X}$ be a $\mathbb{P}_\alpha$-name for the above displayed set and $p\in G_{\alpha}$ forcing that this name consist of all $\nu<\alpha$ such that $t{}^\smallfrown \langle \nu\rangle$ is $\tau$-coherent. We check that $p^{\frown}\langle t\restriction\sigma,\dot{T}\rangle^{\frown}(j_{\alpha,\sigma}(p)\setminus (\alpha+1))$ (for any possible $\mathbb{P}_\alpha$-name  $\dot{T}$) forces $``\check{\alpha}\in j_{\alpha,\sigma}(\dot{X})$", i.e. it forces that the sequence $t^{\frown}\langle\check{\alpha}\rangle$ is $j_{\alpha,\sigma}(\tau)$-coherent. Since $$o^{j_{\alpha,\sigma}(\mathcal{U})}(\alpha)=\sigma<\tau=j_{\alpha,\sigma}(f^\alpha_\tau)(\alpha)=j_{\alpha,\sigma}(f)^{j_{\alpha,\sigma}(\alpha)}_{j_{\alpha,\sigma}(\tau)}(\alpha),$$ the only thing to check is condition (2), which is moreover obvious for $i<|t|$. For $i=|t|$, i.e. $\delta_i=\alpha$, $n_i$ chosen in (2) is equal to the $i$ chosen in the definition of $t\restriction\sigma$, since $j_{\alpha,\sigma}(\vec{f})_{\sigma}^{\alpha}=f_{\sigma}^{\alpha}$. By definition, $p^{\frown}\langle t\restriction\sigma,\dot{T}\rangle$ forces that $b_{t\restriction\sigma}\sqsubseteq \dot{b}_{\alpha}$ and so $p^{\frown}\langle t\restriction\sigma,\dot{T}\rangle^{\frown} (j_{\alpha,\sigma}(p)\smallsetminus(\alpha+1))$ forces $\check{\alpha}\in j_{\alpha,\sigma}(\dot{X})$.}
\end{proof}

\begin{definition}
    Conditions in $\mathbb{Q}^\tau_\alpha$ consists of pairs $\langle t, T\rangle$ such that $t$ is a $\tau$-coherent sequence and $T$ is a tree on $[\alpha]^{<\omega}$, with root $\emptyset$, such that:
    \begin{itemize}
        \item Every $\nu\in\Suc_T(\emptyset)$ is such that $\max(t)<\nu$.
        \item
       {For each $\vec{\nu}\in T$, $\Suc_T(\vec{\nu})\subseteq\{\nu<\alpha\mid\text{$t^{\frown}\vec{\nu}^{\frown}\langle\nu\rangle$ is $\tau$-coherent}\}$ and for every $\sigma<\tau$, $\Suc_T(\vec{\nu})\in\mathcal{U}_{\alpha,\sigma}(t^\frown\vec{\nu}$).}
    \end{itemize}
    We say that $\langle t, T\rangle, \langle s, S\rangle$ are \emph{equivalent} if $b_t=b_s$ and $S=T$.
    
    Given  $\langle t, T\rangle, \langle s, S\rangle\in \mathbb{Q}^\tau_\alpha$ we write $\langle t, T\rangle\leq^* \langle s, S\rangle$ if $b_t=b_s$ and $T\s S$.
\end{definition}

$\mathbb{Q}_{\alpha}^{\tau}$ is nonempty: Let $t:=\emptyset$ and let $T$ be the tree on $[\alpha]^{<\omega}$ defined inductively such that, for every $\vec{\nu}\in T$, $\Suc_T(\vec{\nu})$ consists of all $\nu<\alpha$ such that $t^{\frown}\vec{\nu}^{\frown}\nu$ is $\tau$-coherent. By Lemma \ref{Lemma: Filter Properties}, this  is in $\mathcal{U}_{\alpha,\sigma}(t^{\frown}\vec{\nu})$ for  $\sigma<\tau$.

\begin{remark}[On our $\leq^*$-order]
    Note that our pure extension ordering differs from Ben-Neria--Unger:  In \cite{BenUng}, $\langle t,T\rangle\leq^*\langle s,S\rangle$ if $t=s$ and $T\subseteq S$, so our pure extension ordering is weaker. This is done so that, in the proof of homogeneity of Section~\ref{sec: homogeneity}, we may work with stems which are of minimal length.

    However, this change is a non-issue: Clearly, if $\leq_0\subseteq\leq_1\subseteq\leq$ and $\langle\dP,\leq,\leq_0\rangle$ has the Prikry property, so does $\langle\dP,\leq,\leq_1\rangle$. Therefore, the only potential problem from relaxing the pure extension ordering is losing the closure of $\langle\dQ_{\alpha}^{\tau},\leq^*\rangle$, since all other properties we need just depend on $\langle\dQ_{\alpha}^{\tau},\leq\rangle$.
\end{remark}

\begin{lemma}[Minimal sequences]\label{lemma: minimal}
    The following hold:
    \begin{enumerate}
        \item If $t,s$ are minimal $\tau$-coherent sequences, $b_t=b_s$ implies $t=s$.
        \item A $\tau$-coherent sequence $t$ is minimal if and only if $|t|\leq|s|$ whenever $s$ is another $\tau$-coherent sequence with $b_t=b_s$.
    \end{enumerate}
\end{lemma}

\begin{proof}
 (1).   We prove this by induction on the length of $t$. Clearly, if $|t|=0$ then $b_t=\emptyset=b_s$, ergo $s=\emptyset=t$, as sought.
    
    So suppose the statement holds for all sequences $t$ of length $n$. Let $t$ and $s$ be minimal $\tau$-coherent sequences,  with $|t|=(n+1)$ and $b_t=b_s$. Suppose $t=\langle\delta_0,\dots,\delta_n\rangle$ and $s=\langle\gamma_0,\dots,\gamma_m\rangle$. Clearly, we have have $$\gamma_m=\max(s)=\max(b_s)=\max(b_t)=\max(t)=\delta_{n}.$$ By minimality, $o^{\mathcal{U}}(\delta_{n-1})\geq f_{o^{\mathcal{U}}(\delta_n)}^{\delta_n}(\delta_{n-1})$. In particular, $\delta_{n-1}$ is not a member of the Magidor sequence $b_{\delta_n}$ (as $\min(b_{\delta_n})>\delta_{n-1}$) and so $b_{t\restriction n}$ is the initial segment of $b_t$ obtained by deleting $b_{\delta_n}$. By the same argument, $b_{s\restriction m}$ is the initial segment of $b_s$ obtained by deleting $b_{\gamma_m}=b_{\delta_n}$. Therefore, $b_{t\restriction n}=b_{s\restriction m}$ and so $t\restriction n=s\restriction m$, by the inductive hypothesis, showing that $t=s$.

    \smallskip

   (2). We first show the backward direction. So suppose that $t=\langle\delta_0,\dots,\delta_n\rangle$ is not minimal. Then there is $0<i<|t|$ such that $o^{\mathcal{U}}(\delta_{i-1})<f_{o^{\mathcal{U}(\delta_i)}}^{\delta_i}(\delta_{i-1})$. But then, letting $s:=\langle\delta_0,\dots,\delta_{i-2},\delta_i,\dots,\delta_n\rangle$ we have that $b_s=b_t$ and $|s|<|t|$, which yields a contradiction.

    Now we prove the forward direction of (2). So suppose $t$ is a minimal $\tau$-coherent sequence and $s$ is another $\tau$-coherent sequence such that $b_t=b_s$ and $|s|<|t|$. By repeating the procedure in the backward direction, we can reduce $s$ over and over again to make it minimal and therefore, without loss of generality, assume that $s$ is already minimal. But then, by (1), $s=t$, a contradiction.
\end{proof}

The following definition is enabled by the previous lemma:

\begin{definition}
    Let $t$ be a $\tau$-coherent sequence. The \emph{regularization of $t$},  $\reg(t)$, is the unique minimal sequence such that $b_{\reg(t)}=b_t$.
\end{definition}

\begin{remark}
    The proof of Lemma \ref{lemma: minimal} shows that $\reg(t)$ can be constructed by the following procedure: Define $t_0:=t$. If $t_n=\langle\delta_0^n,\dots,\delta_{m_n-1}^n\rangle$ has been defined, check if $t_n$ is minimal. If so, set $t_n:=\reg(t)$. Otherwise, there is a least index $0<i<m_n$ such that $o^{\mathcal{U}}(\delta_{i-1})< f_{o^{\mathcal{U}}(\delta_i)}^{\delta_i}(\delta_{i-1})$. Let
    $$t_{n+1}:=\langle\delta_0^n,\dots,\delta_{i-2}^n,\delta_i^n,\dots,\delta_{m_n-1}^n\rangle$$
   Because $|t|$ is finite, there must be some $n<|t|$ such that $t_n$ is minimal and thus equal to $\reg(t)$.
\end{remark}

\begin{definition}[One-point extensions]
    Given $\langle t, T\rangle\in \mathbb{Q}^\tau_\alpha$ and $\langle \nu\rangle\in T$ we denote by $\langle t, T\rangle\cat \langle \nu\rangle$ the pair $\langle t{}^\smallfrown \langle\nu\rangle, T_{\langle \nu\rangle}\rangle$ where $T_{\langle \nu\rangle}$ denotes the tree $T$ ``above $\nu$''; to wit, $\{\vec\eta\in [\alpha]^{<\omega}\mid \langle \nu\rangle{}^\smallfrown\vec\eta\in T\}$.  More generally, for $\vec\nu \in T$, $\langle t, T\rangle\cat \vec\nu$ is defined by recursion  setting $\langle t, T\rangle\cat \emptyset := \langle t, T\rangle$. 
\end{definition}
The following is easy to see: 

\begin{fact}
    Given  $\langle t,T\rangle\in\dQ_{\alpha}^{\tau}$ and $\vec{\nu}\in T$, $\langle t,T\rangle\cat\vec{\nu}\in\dQ_{\alpha}^{\tau}$.\qed
\end{fact}

\begin{definition}
     Given  $\langle t, T\rangle, \langle s, S\rangle\in \mathbb{Q}^\tau_\alpha$ we write $\langle t, T\rangle\leq \langle s, S\rangle$ if there is $\vec\nu \in S$ such that $\langle t, T\rangle\leq^* \langle s, S\rangle\cat \vec\nu$.
\end{definition}

\begin{remark}\label{remark: same measures}
    Note that if $\langle t, T\rangle, \langle s, S\rangle$ are equivalent  then $\langle t, T\rangle\leq^* \langle s, S\rangle$ and $\langle s, S\rangle\leq^* \langle t, T\rangle$. In particular, they force the same information. So, whenever $t$ and $s$ are equivalent $\tau$-coherent sequences and $\sigma<o^{\mathcal{U}}(\alpha)$, $$\mathcal{U}_{\alpha,\sigma}(t)=\mathcal{U}_{\alpha,\sigma}(s).$$
\end{remark}
\begin{lemma}[Basic Properties]\label{lemma: basic properties of Gitik's}\hfill
\begin{enumerate}
    \item  $\langle \mathbb{Q}^\tau_\alpha,\leq\rangle $  is $\alpha^+$-cc;
    \item $\langle \mathbb{Q}_\alpha^\tau,\leq^*\rangle$ is $\alpha$-directed-closed;
    \item $\langle \mathbb{Q}^\tau_\alpha,\leq\rangle$ adds a club set $b_\alpha \s \alpha$ with $\otp(b_\alpha)=\omega^{\tau}$ consisting of $V_0$-measurable cardinals.
\end{enumerate}  
\end{lemma}

\begin{proof}
    For (1) and (3), see e.g.  \cite[Lemma~5.31]{HP2}.

    To see (2), let $\{\langle t_i,T_i\rangle\;|\;i<j\}$ for $j<\alpha$ be a $\leq^*$-directed set. In particular, for every $i,i'<j$, $t_i$ and $t_{i'}$ are equivalent. This  implies that the trees $T_i$ split according to the same measures (see Remark~\ref{remark: same measures}), so, letting $T^*:=\bigcap_{i<j}T_i$ we obtain a legitimate tree. Thus,  $\langle t_0,T^*\rangle$ is a lower bound of the sequence $\{\langle t_i,T_i\rangle\;|\;i<j\}$.
\end{proof}

Let us now tackle the usual Prikry-type analysis of the poset.

\begin{definition}
    Fix $\tau\leq o^{\mathcal{U}}(\alpha)$. Given a condition $p=\langle t, T\rangle$ in $\mathbb{Q}^\tau_\alpha$, a tree $\Upsilon\in [\alpha]^{\leq n}$ (with $n<\omega$)  is called a \emph{$p$-fat tree (of height $n$)} if $\Upsilon$ consists of increasing sequences,  $\Upsilon\s T$, and for each  $\vec\nu\in [\Upsilon]^{\leq m}$ with $m<n$, there is an ordinal $\tau(\vec\nu)<\tau$ such that $\Suc_{\Upsilon}(\vec\nu)\in \mathcal{U}_{\alpha,\tau(\vec\nu)}(t{}^\smallfrown \vec\nu)$.

   We shall denote by $[\Upsilon]$ the set of maximal branches in $\Upsilon$.
\end{definition}

\begin{remark}
    Let $D\s \mathbb{Q}^\tau_\alpha$ be a dense open set. Notice that $D$ is closed under equivalence; to wit, if $p\in D$ and $q\in \mathbb{Q}^\tau_\alpha$ is equivalent to $p$ then $q\in D$. This holds simply because $q\leq p$ and $D$ is open.
\end{remark}

\begin{lemma}[Strong Prikry Property]\label{lemma: SPP}
Fix $\tau\leq o^{\mathcal{U}}(\alpha)$.
   Let $p\in \mathbb{Q}^\tau_\alpha$ and $D\s \mathbb{Q}^\tau_\alpha$ be a dense open set. Then, there is $p^*\leq^* p$ and a $p^*$-fat tree $\Upsilon$ such that $p^*\cat \vec\nu\in D$ for all $\vec\nu\in [\Upsilon]$. 
\end{lemma}
\begin{proof}
  Let $p=\langle t, T\rangle$ and $D\s \mathbb{Q}^{\tau}_\alpha$ be as in the lemma.  
  
  If there is $p^*\leq^* p$ in $D$ then we choose such a condition and let $\Upsilon := \emptyset$ and we are done. So, suppose that no $\leq^*$-extension of $p$ enters the dense set $D$. 
  
  For each $\tau<o^{\mathcal{U}}(\alpha)$ let us consider 
  $$X_\tau^0:=\{\delta\in \Suc_{T,\tau}(\emptyset)\mid \exists S\, (\langle t{}^\smallfrown \delta, S\rangle\leq^* p\cat \langle \delta\rangle\,\text{and}\, \langle t{}^\smallfrown \delta, S\rangle\in D)\}$$
  and let $X^1_\tau$ be the complement of $X^0_\tau$.

  For each $\tau<o^{\mathcal{U}}(\alpha)$ there is a unique  $i_{\tau}\in \{0,1\}$ such that $X^{i_\tau}_\tau\in \mathcal{U}_{\alpha,\tau}(t)$. 

  Let $T_1$ be  the tree whose first level is $X=\bigcup_{\tau<o^{\mathcal{U}}(\alpha)}X^{i_\tau}_\tau$ and for each $\delta\in X$ if there is some $\tau<o^{\mathcal{U}}(\alpha)$ such that $i_\tau=0$ then $(T_1)_{\langle \delta\rangle}$ is some tree $S(\delta)$ that witnesses it; otherwise, stipulate that $(T_1)_{\langle \delta\rangle}=T_{\langle \delta\rangle}.$

  Set $p_1:=\langle t, T_1\rangle$. Note that $p_1\leq^* p$, ergo there is no $q\leq^* p_1$ in $D$. However, if there is a one-step-extension $q\leq p_1$ in $D$ then there is a $p_1$-fat tree $\Upsilon$ of height $1$ such that $p_1\cat \langle \nu\rangle\in D$ for all $\langle \nu\rangle\in \Upsilon$. Indeed, if there is such a $q$ then $q\leq^* p_1\cat \langle \delta\rangle$ for some $\delta\in X^{i_\tau}_\tau$. Note that $i_\tau=0$ because the tree $T^q$ witnesses that $\delta\in X^0_\tau$.  Ergo, for every $\delta'\in X^0_\tau$, $$p_1\cat \langle \delta'\rangle = \langle t{}^\smallfrown \langle \delta'\rangle, (T_1)_{\langle \delta'\rangle}\rangle = \langle t{}^\smallfrown \langle \delta'\rangle, S(\delta')\rangle\in D.$$
  In this scenario, choosing $p^*=p_1$ and $\Upsilon = X^0_\tau$ we are done.

  \smallskip

  Suppose that no one-point extension of $p_1$ enters $D$. In that case we continue by induction, thus producing a $\leq^*$-decreasing sequence of conditions $\langle p_n\mid n<\omega\rangle$, $\leq^*$-below $p$, such that the following hold:
  \begin{equation*}
      \text{$\exists q\,(q\leq^n p_n$ and $q\in D\,\Rightarrow\, \exists \Upsilon(\Upsilon$ is a $p_n$-fat tree and $p_n\cat \vec\nu\in D$ for $\vec\nu\in [\Upsilon])$.}
  \end{equation*}
Let us spell out the details of the construction for the case $n=2$.

Suppose that $p_1$ is given. For each $\nu\in \mathrm{Lev}_1(T_1)$ and $\tau<o^{\mathcal{U}}(\alpha)$ let 
 $$X_\tau^0(\nu):=\{\delta\in \Suc_{T,\tau}(\langle\nu\rangle)\mid \exists S\, (\langle t{}^\smallfrown \langle\nu, \delta\rangle, S\rangle\leq^* p_1\cat \langle \nu, \delta\rangle\,\text{and}\, \langle t{}^\smallfrown \langle \nu, \delta\rangle, S\rangle\in D)\}$$
 and $X^1_\tau(\nu)$ be the complement of $X^0_\tau(\nu).$

 For each $\nu\in \mathrm{Lev}_1(T_1)$ and $\tau<o^{\mathcal{U}}(\alpha)$ let $i_{\nu, \tau}\in\{0,1\}$ be the unique index such that $X^{i_\tau}_\tau(\nu)\in \mathcal{U}_{\alpha,\tau}(t{}^\smallfrown\langle\nu\rangle).$ Set $X(\nu):=\bigcup_{\tau<o^{\mathcal{U}}(\alpha)}X_{\tau}^{i_\tau}(\nu)$. Clearly, this belongs to the intersection of the ultrafilters $\mathcal{U}_{\alpha,\tau}(t{}^\smallfrown\langle\nu\rangle).$

 \smallskip

 Now, for each $\tau,\sigma<o^{\mathcal{U}}(\alpha)$ and $j\in \{0,1\}$, look at 
 $$Y_{\tau}^{\sigma,j}:=\{\nu\in \mathrm{Succ}_{T_1,\sigma}(\emptyset)\mid X^{j}_\tau(\nu)=X^{i_{\nu,\tau}}_\tau(\nu)\}.$$
For each $\tau,\sigma$ there is a unique $j_{\tau,\sigma}$ such that $Y_\tau^{\sigma, j_{\sigma,\tau}}\in \mathcal{U}_{\alpha,\sigma}(t)$.

Let $Y^{\sigma}:=\bigcup_{\tau<o^{\mathcal{U}}(\alpha)}Y^{\sigma, j_{\sigma,\tau}}_\tau$. Clearly, this belongs to $\mathcal{U}_{\alpha,\sigma}(t)$.

\smallskip

Let $T_2$ be the tree whose first level is $\bigcup_{\sigma<o^{\mathcal{U}}(\alpha)} Y^\sigma$.

For each $\nu\in \mathrm{Lev}_1(T_2)$ let   $$\textstyle \Suc_{T_2}(\langle\nu\rangle):=\bigcup_{\{\sigma\mid \nu\in Y^\sigma\}}\bigcup_{\tau<o^{\mathcal{U}}(\alpha)} X^{j_{\sigma,\tau}}_\tau( \nu).$$
This belongs to the intersection of the filters $\mathcal{U}_{\alpha,\tau}(t{}^\smallfrown\langle \nu\rangle)$ for $\tau<o^{\mathcal{U}}(\alpha).$

Now, for each $\delta\in \Suc_{T_2}(\langle \nu\rangle)$, if there are $\sigma,\tau$ such that $j_{\sigma,\tau}=0$ and $\delta\in X^0_\tau(\langle \nu\rangle)$ then we stipulate that 
$(T_2)_{\langle \nu,\delta\rangle}:=S(\nu,\delta)$ for some tree witnessing it. Otherwise, $(T_2)_{\langle \nu,\delta\rangle}:=(T_1)_{\langle \nu,\delta\rangle}.$ Note that $T_2$ is a tree for a condition.

\smallskip

Let $p_2:=\langle t, T_2\rangle.$ Clearly, $p_2\leq^*p_1$.

Suppose that there was a two-step-extension $q\leq p_2$ inside $D$. Then, $q\leq^* p_2\cat \langle \nu, \delta\rangle$ where $\nu\in Y^\sigma$ and $\delta\in X^{j_{\sigma',\tau}}_\tau(\nu)$ for some $\sigma, \tau,\sigma'<o^{\mathcal{U}}(\alpha)$.

Observe that $T^q$ witnesses that $\delta\in X^0_\tau(\nu)$ ergo $j_{\sigma',\tau}=0$. Taking any $\nu'\in Y^{\sigma',0}_\tau$ and any $\delta'\in X^{0}_\tau(\nu')$ we have that 
$$p_2\cat \langle \nu',\delta'\rangle = \langle t{}^\smallfrown \langle\nu',\delta'\rangle, (T_2)_{\langle \nu',\delta'\rangle}\rangle= \langle t{}^\smallfrown \langle\nu',\delta'\rangle, S(\nu',\delta')\rangle\in D.$$
Let $\Upsilon$ be the $p_2$-fat tree whose first level is $Y^{\sigma',0}_\tau$ and for each $\nu'$ there in $\Suc_{\Upsilon}(\langle \nu\rangle)=X^0_\tau(\nu')$. It follows that $p_2$ and $\Upsilon$ witness the desired property.

\smallskip

Once the sequence $\langle p_n\mid n<\omega\rangle$ is defined we let $p^*$ a $\leq^*$-lower bound for the $p_n$'s. Since $D$ is dense there is, for some $n<\omega$, an $n$-step-extension of $p^*$ inside $D$. Then, there is an $n$-step extension of $p^n$ inside $D$. By the induction hypothesis, there is a $p_n$-fat tree $\Upsilon$ such that $p_n\cat \vec\nu\in D$ for all $\vec\nu\in [\Upsilon]$. Note that $\Upsilon \cap T^{p^*}$ is a $p^*$-fat tree and that $p^*\cat\vec\nu\in D$ for all $\vec\nu\in [\Upsilon]$. This completes the proof of the lemma.
\end{proof}

\begin{theorem}\label{thm: regularity}
${\mathbb{Q}^{\alpha^+}_\alpha}$ forces $``\alpha$ is regular".
\end{theorem}
\begin{proof}
Suppose otherwise. Let $p=\langle t, T\rangle\in \mathbb{Q}^{\alpha^+}_\alpha$ and $\dot{f}$ be a $\mathbb{Q}^{\alpha^+}_\alpha$-name for a function such that $p\forces_{\mathbb{Q}^{\alpha^+}_\alpha}``\dot{f}\colon \mu \to \alpha$ is cofinal" for some $\mu < \alpha$.

For each $\xi < \mu$ consider the dense open set defined by 
$$D_\xi := \{q\in \mathbb{Q}^{\alpha^+}_\alpha\mid q\perp p\,\vee\, q\forces_{\mathbb{Q}^{\alpha^+}_\alpha}\dot{f}(\xi) < \max(\dot{b}_q)\}.$$
By the Strong Prikry Property (Lemma~\ref{lemma: SPP}) there is $p_\xi \leq^* p$ and $\Upsilon_\xi$ a $p_\xi$-fat tree of height $n_\xi$. Since $\mu < \alpha$ and $\langle \mathbb{Q}^{\alpha^+}_\alpha, \leq^*\rangle$ is $\alpha$-closed, we can arrange  $\langle p_\xi\mid \xi <\mu \rangle$ be $\leq^*$-decreasing and $p^*$ a $\leq^*$-lower bound.

\smallskip

For each $\xi<\mu$ and $\vec\nu\in \Upsilon_\xi \setminus [\Upsilon_\xi]$ there is $\tau(\xi,\vec\nu)<o^{\mathcal{U}}(\alpha)=\alpha^+$ such that $$\Suc_{\Upsilon_\xi}(\vec\nu)\in \mathcal{U}_{\alpha,\tau(\xi,\vec\nu)}(t\cat \vec\nu).$$ Let $\tau^* :=\left(\sup_{\xi <\mu}\sup_{\vec\nu\in \Upsilon_\xi}\tau(\xi,\vec\nu)\right)+1$ noting that $\tau^* < \alpha^+$.

\smallskip

The following claim is standard, but we include it to convince the reader that our slight variation of coherent sequence is appropriate:
\begin{claim}
    For each $\sigma_0\leq\sigma_1$ and a $\sigma_0$-coherent  $s\in [\alpha]^{<\omega}$, the map $\pi^s_{\sigma_1, \sigma_0}\colon \alpha \to \alpha$ given by
    $$\pi^s_{\sigma_1,\sigma_0}(\nu):=\min \{x\in b_\nu\mid x>\max(s)\,\wedge\, o^{\mathcal{U}}(x)=f_{\sigma_0}(x)\}$$
    defines a Rudin-Keisler projection from $\mathcal{U}_{\alpha,\sigma_1}(s)$ to $\mathcal{U}_{\alpha,\sigma_0}(s)$.\footnote{Note that for $\sigma_0=\sigma_1$, the projection is simply the identity on a measure-one set, since for $o^{\mathcal{U}}(\nu)=f_{\sigma_0}(\nu)$, the only point $x$ in $b_{\nu}$ of order $f_{\sigma_0}(x)$ is $\nu$ itself.}
\end{claim}

\begin{proof}
    We have to show that whenever $A\in\mathcal{U}_{\alpha,\sigma_0}(s)$, $(\pi_{\sigma_1,\sigma_0}^s)^{-1}[A]\in\mathcal{U}_{\alpha,\sigma_1}$. We prove the following stronger statement by induction on $\sigma_1$: Whenever $A\in\mathcal{U}_{\alpha,\sigma_0}(s)$ the set of all $\nu<\alpha$ such that either $o^{\mathcal{U}}(\nu)<f_{\sigma_0}(\nu)$ or $\pi_{o^{\mathcal{U}}(\nu),\sigma}(\nu)\in A$ is an element of $\bigcap_{\sigma<\sigma_1}\mathcal{U}_{\alpha,\sigma}(s)$. To this end, we may assume without loss of generality that $A$ only contains points above $\max(s)$.

    The statement clearly holds for $\sigma_1=\sigma_0$ since there is nothing to do. Also, there is nothing to do in limit steps. So suppose that $\sigma_1=\sigma_1'+1$. We must show that $A\in\mathcal{U}_{\alpha,\sigma_1'}$. To this end, consider the condition $(s\restriction\sigma_1',T)$, where $T$ is a tree on $[\alpha]^{<\omega}$ with $\Suc_T(\vec{\nu})\subseteq A$ for each $\vec{\nu}\in T$. By the inductive hypothesis, $(s\restriction\sigma_1',T)$ is a condition in $\mathbb{Q}_{\alpha}^{\sigma_1'}$. Also, it is easy to see that whenever $H$ is $\mathbb{Q}_{\alpha}^{\nu}$-generic containing $(s\restriction\sigma_1',T)$, $j(\pi_{\sigma_1',\sigma_0})(\alpha)\in A$.
\end{proof}

For $\xi<\mu$ and $\nu\in\dom_1(\mathcal{U})$, by induction on $i<n_{\xi}$, let $$\text{$\pi_0^{\xi}(\nu):=\pi_{\tau^*,\tau(\xi,\emptyset)}^t(\nu)\;$ and $\;\pi_{i+1}^{\xi}(\nu):=\pi_{\tau^*,\tau(\xi,\langle\pi_0^{\xi}(\nu),\dots,\pi_i^{\xi}(\nu)\rangle)}^{t^{\frown}\langle\pi_0^{\xi}(\nu),\dots,\pi_i^{\xi}(\nu)\rangle}(\nu)$.}$$

Assume for the time being that for each $\xi < \mu$
\begin{equation}
  \tag{$\star_\xi$} X_\xi:=\{\nu < \alpha\mid \langle \pi_0^{\xi}(\nu),\dots, \pi_{n_{\xi}-1}^{\xi}(\nu)\rangle\in [\Upsilon_\xi]\}\in \mathcal{U}_{\alpha,\tau^*}(t) 
\end{equation}

The next claim yields the desired contradiction:
\begin{claim}
If $(\star_\xi)$ holds for all $\xi < \mu$ then there is $p^{**}\leq^* p^*$ such that for each $\nu\in \mathrm{Suc}_{T^{p^{**}},\tau^*}(\emptyset)$,  $p^{**}\cat \langle \nu\rangle \forces_{\mathbb{Q}^{\alpha^+}_\alpha}\mathrm{ran}(\dot{f}) <\nu$.
\end{claim}
\begin{proof}[Proof of claim]
Let $p^{**}\leq^* p^*$ be such that $\mathrm{Suc}_{T^{p^{**}},\tau^*}(\emptyset)\s \bigcap_{\xi < \mu} X_\xi$. Let  $\nu$ be in this set and $\xi < \mu$ arbitrary. Since $\nu\in X_\xi$ it follows that $$\langle \nu^\xi_0,\dots,\nu^\xi_{n_\xi-1}\rangle{}^\smallfrown\langle \nu\rangle\equiv \langle \nu\rangle;$$
namely, both sequences induce the same Prikry/Magidor sequence $b_{\nu}\cup\{\nu\}$.

As a result, $p^{**}\cat \langle \nu^\xi_0,\dots,\nu^\xi_{n_\xi-1}\rangle{}^\smallfrown\langle \nu\rangle$ forces the same as $p^{**}\cat \langle \nu\rangle$. However,
$$p^{**}\cat \langle \nu^\xi_0,\dots,\nu^\xi_{n_\xi-1}\rangle{}^\smallfrown\langle \nu\rangle \leq p_\xi \cat \langle \nu^\xi_0,\dots,\nu^\xi_{n_\xi-1}\rangle\forces_{\mathbb{Q}^{\alpha^+}_\alpha}\dot{f}(\xi) < \nu^\xi_{n_\xi-1}< \nu.$$
This gives the desired result.
\end{proof}
\begin{claim}
   For every $\xi<\mu$, $(\star_{\xi})$ holds in $V[G]$.
\end{claim}

\begin{proof}[Proof of Claim]%
    Fix $\xi<\mu$ and assume toward a contradiction that $X_{\xi}\notin\mathcal{U}_{\alpha,\tau^*}(t)$. Hence, $\alpha\smallsetminus X_{\xi}\in\mathcal{U}_{\alpha,\tau^*}(t)$ and thus there are $p\in G_{\alpha}$ and a $\dP_{\alpha}$-name $\dot{T}$ such that $p$ decides the height of $\dot{\Upsilon}_{\alpha}$ to be $n_{\xi}$, forces $\dot{X}_{\xi}$ to be given by the definition above, forces $\langle t,\dot{T}\rangle\in\dot{\dQ}_{\alpha}^{\tau^*}$ and such that
    $$p^{\frown}\langle t,\dot{T}\rangle^{\frown}(j_{\alpha,\tau^*}(p)\smallsetminus\alpha+1)\Vdash_{j_{\alpha,\tau^*}(\dP_{\alpha})}\check{\alpha}\notin j_{\alpha,\tau^*}(\dot{X}_{\xi})$$
    We will show that this is impossible.

    \smallskip
    
    Work in $V[G_{\alpha}]$. Denote $T:=\dot{T}^G$. We will first find a sequence $\langle\nu_0,\dots,\nu_{n_{\xi}-1}\rangle$ in $ T\cap[\Upsilon_{\xi}]$. We construct the sequence by induction, so suppose that $\langle\nu_0,\dots,\nu_{i-1}\rangle$ has been defined for $i<n_{\xi}$. We have
    $$\Suc_T(\langle\nu_0,\dots,\nu_{i-1}\rangle)\in\mathcal{U}_{\alpha,\tau(\xi,\langle\nu_0,\dots,\nu_{i-1}\rangle)}(t^{\frown}\langle\nu_0,\dots,\nu_{i-1}\rangle)$$
    by the definition of $\dQ_{\alpha}^{\alpha^+}$ and
    $$\Suc_{\Upsilon_{\xi}}(\langle\nu_0,\dots,\nu_{i-1}\rangle)\in\mathcal{U}_{\alpha,\tau(\xi,\langle\nu_0,\dots,\nu_{i-1}\rangle)}(t^{\frown}\langle\nu_0,\dots,\nu_{i-1}\rangle)$$
    by the definition of $\Upsilon_{\xi}$. Therefore, there must exist $\nu_i$ as required.

    \smallskip

    By extending $p$ inside of $G$ we may assume that $$\text{$p\forces_{\mathbb{P}_\alpha} \langle\check{\nu}_0,\dots,\check{\nu}_{n_{\xi}-1}\rangle\in\dot{T}\cap [\dot{\Upsilon}_{\xi}]$},$$  in particular, by elementarity $j_{\alpha,\tau^*}(p)$ forces $\langle\check{\nu}_0,\dots,\check{\nu}_{n-1}\rangle\in[j_{\alpha,\tau^*}(\dot{\Upsilon}_{\xi})]$ over $j_{\alpha,\tau^*}(\dP_{\alpha})$ (since the critical point of $j_{\alpha,\tau^*}$ is $\alpha$, greater than all of the $\nu_i$). Furthermore, it follows easily from the definition of the projection maps that $$p^{\frown}\langle t^{\frown}\langle\nu_0,\dots,\nu_{n_{\xi}-1}\rangle,\dot{T}_{\langle\nu_0,\dots,\nu_{n_{\xi}-1}\rangle}\rangle\forces_{\dP_{\alpha}*\dot{\dQ}_{\alpha}^{\tau^*}}\langle\dot{\pi}_0^{\xi}(\check{\alpha}),\dots,\dot{\pi}_{n_{\xi}-1}^{\xi}(\check{\alpha})\rangle=\langle\check{\nu}_0,\dots,\check{\nu}_{n_{\xi}-1}\rangle.$$ Indeed, we see inductively that the first element $x$ of $\dot{b}_{\alpha}$ above $\max(t)$ with $o^{\mathcal{U}}(x)=f_{\tau(\xi,\emptyset)}(x)$ is forced to be $\check{\nu}_0$, the first element of $\dot{b}_{\alpha}$ above $\max(t^{\frown}\langle\nu_0\rangle)$ with order $\tau(\xi,\langle\nu_0\rangle)$ is forced to be $\nu_1$ and so on. As a result, the condition $p^{\frown}\langle t^{\frown}\langle\nu_0,\dots,\nu_{n_{\xi}-1}\rangle,\dot{T}_{\langle\nu_0,\dots,\nu_{n_{\xi}-1}\rangle}\rangle^{\frown}(j_{\alpha,\tau^*}(p)\smallsetminus\alpha+1)$    forces  $\langle\dot{\pi}_0^{\xi}(\check{\alpha}),\dots,\dot{\pi}_{n_{\xi}-1}^{\xi}(\check{\alpha})\rangle=\langle\check{\nu}_0,\dots,\check{\nu}_{n_{\xi}-1}\rangle$ and $\langle\check{\nu}_0,\dots,\check{\nu}_{n_{\xi}-1}\rangle\in[j_{\alpha,\tau^*}(\dot{\Upsilon}_{\xi})]$.

    Therefore, that condition forces $\langle\dot{\pi}_0^{\xi}(\check{\alpha}),\dots,\dot{\pi}_{n_{\xi}-1}^{\xi}(\check{\alpha})\rangle\in[j_{\alpha,\tau^*}(\dot{\Upsilon}_{\xi})]$ and thus $\check{\alpha}\in\dot{X}_{\xi}$, a contradiction.
\end{proof}

This completes the proof of the theorem.
\end{proof}

\begin{prop}\label{prop: distributivity}
Suppose that $\mathbb{Q}$ is an $\alpha^+$-proper forcing that preserves the regularity of $\alpha$ and does not add bounded subsets to $\alpha$. Then, $\mathbb{Q}$ is $\alpha$-distributive.  In particular, $\mathbb{Q}^{\alpha^+}_\alpha$ is $\alpha$-distributive.
\end{prop}%
\begin{proof}
Let $G\s \mathbb{Q}$ be a generic filter. 
   We show that $({}^\mu\mathrm{Ord})^{V[G]}=({}^\mu\mathrm{Ord})^{V}$.

   Let $p\in G$ and $\dot{f}$ be a $\mathbb{Q}$-name such that $p\forces_{\mathbb{Q}}\dot{f}\colon \mu\to \mathrm{Ord}$. Let $\mathcal{M}\prec H_\Theta$  such that $|\mathcal{M}|=\alpha$, ${}^{<\alpha}\mathcal{M}\s \mathcal{M}$, $\mathbb{Q}, p,\dot{f}\in \mathcal{M}$. For each $\beta < \mu$ consider
   $$D_\beta=\{q\in \mathbb{Q}\mid \exists \rho\, (q\forces_{\mathbb{Q}}\dot{f}(\beta)=\rho)\}.$$
   This is a dense open subset of $\mathbb{Q}$ in $\mathcal{M}$. By properness, there is $q\leq p$, $q\in G$, such that $q\forces \forall \beta < \mu\, (D_\beta \cap \mathcal{M}\cap \dot{G} \neq \emptyset)$. By elementarity of $\mathcal{M}$,
$$q\forces \mathrm{range}(\dot{f})\s \mathcal{M}$$
hence $\mathrm{range}(\dot{f}_G)\s \mathcal{M}$. Working in $V$, let $\varphi\colon \alpha\leftrightarrow \mathcal{M}$ be a bijection. Then, $\varphi^{-1}\circ \dot{f}_G\colon \mu\to \alpha$. Since $\mathbb{Q}$ preserves the regularity of $\alpha$ there is  $\xi < \alpha$ such that $\mathrm{range}(\psi)\s \xi$. As such, $\psi$ can be coded as a bounded subset of $\alpha$ (for instance, using Gödel's pairing function) thus obtaining $\psi\in V$. Finally, $\dot{f}_G = \varphi\circ \psi$ is in the ground model, as claimed.
\end{proof}

\begin{remark}
    By the same argument, if $\dQ$ is a forcing which has the $\alpha^+$-covering property, preserves the regularity of $\alpha$ and does not add bounded subsets to $\alpha$, then $\dQ$ is $\alpha$-distributive.
\end{remark}

\subsection{The iteration when $o^{\mathcal{U}}(\alpha)=\alpha^+$}\label{Subsection: Prikry Club-Shooting}

We now describe the iterand $\dot{\dQ}_{\alpha}$ for $o^{\mathcal{U}}(\alpha)=\alpha^+$. Our goal is the following: At stage $\kappa$, we want to shoot a club through the set of $V$-regular cardinals below $\kappa$. In order to ensure that this preserves the supercompactness of $\kappa$, we need to have performed the same forcing already at many cardinals below $\kappa$. However, the forcing shooting a club trough the set of $V$-regular cardinals is typically going to destroy stationary subsets of $\kappa$ and therefore cannot be $\kappa$-closed and will merely be $\kappa$-distributive. As merely distributive forcing notions typically cannot be iterated well, we realize our club-shooting forcing as a Prikry-type forcing notion with a $\kappa$-closed pure extension ordering. However, while doing so, we need to ensure that this poset remains homogeneous (at least in the sense stipulated in Definition~\ref{Definition: Automorphism Sequence}). To do so, instead of using a single club-shooting forcing, we use a ${<}\,\alpha$-supported product (see Remark~\ref{keyremark} for details). To streamline the presentation we shall adopt the following convention for the scope of this section:

\begin{notation}
     We write $\forces$ for  $\forces_{\mathbb{Q}^{\alpha^+}_\alpha}$ and our ground model will be $V[\mathbb{P}_\alpha]$. When we want to refer to the initial ground model $V$, we shall write $V_0$· 
\end{notation}

\begin{definition}
    Let $\dC_{\alpha}$ be the following poset: Conditions are partial functions $c\colon\beta^+\to V$ of size ${<}\,\beta$ such that, for every $\xi\in\dom(c)$, $c(\xi)$ is a closed bounded subset of $\beta$ consisting of $V_0$-measurable cardinals.
\end{definition}

From Fact \ref{Fact: Abraham-Shelah}, we can see the following:

\begin{lemma}\label{Lemma: Abraham-Shelah New}
    Suppose that the collection of $V_0$-measurable cardinals below $\alpha$ is a fat stationary set. Then $\dC_{\alpha}$ is $\alpha$-distributive.
\end{lemma}

\begin{proof}
    We can regard $\dC_{\alpha}$ as the product $\dCS_{\alpha}\times\dC_{\alpha}$ and 
    therefore the iteration $\dCS_{\alpha}*\dC_{\alpha}^V$ as well. By Fact \ref{Fact: Abraham-Shelah}, $\dCS_{\alpha}$ is $\alpha$-distributive. Let $G$ be $\dCS_{\alpha}$-generic. In $V[G]$, let $C\subseteq\alpha$ be a club consisting of $V_0$-measurable cardinals. By the $\alpha$-distributivity of $\dCS_{\alpha}$, $\dC_{\alpha}^V=\dC_{\alpha}^{V[G]}$. Moreover, in $V[G]$, $\dC_{\alpha}$ contains an $\alpha$-closed dense subset consisting of those $c\in\dC_{\alpha}$ such that for every $\xi\in\dom(c)$, $\max(c(\xi))\in C$. In particular, $\dC_{\alpha}$ is $\alpha$-strategically closed in $V[G]$. In particular, $\dC_{\alpha}$ is $\alpha$-distributive in $V[G]$ and so, since $G$ was arbitrary, $\dCS_{\alpha}*\dC_{\alpha}^V$ is $\alpha$-distributive in $V$. This yields the lemma.
\end{proof}

Now we can define our preliminary forcing $\mathbb{R}^0_\alpha$:

\begin{definition}
    Let $\mathbb{R}^0_\alpha:=\dQ_{\alpha}^{\alpha^+}*\dot{\dC}_{\alpha}$. For $(p,\dot{c}),(p',\dot{c}')\in\mathbb{R}^0_\alpha$, we let $(p',\dot{c}')\leq_{\mathbb{R}^0_\alpha}(p,\dot{c})$ if $p'\leq_{\dQ_{\alpha}^{\alpha^+}}p$ and $p'\Vdash\dot{c}'\leq_{\dot{\dC}_{\alpha}}\dot{c}$ and we let $(p',\dot{c}')\leq_{\mathbb{R}^0_\alpha}^*(p,\dot{c})$ if $p'\leq_{\dQ_{\alpha}^{\alpha^+}}^*p$ and $p'\Vdash\dot{c}'\leq_{\dot{\dC}_{\alpha}}\dot{c}$.
\end{definition}

The next  describes the  properties of the preliminary triple $\langle \mathbb{R}^0_\alpha,\leq_{\mathbb{R}^0_\alpha}, \leq^*_{\mathbb{R}^0_\alpha}\rangle$:
\begin{lemma}[Properties of $\langle \mathbb{R}^0_\alpha,\leq_{\mathbb{R}^0_\alpha}, \leq^*_{\mathbb{R}^0_\alpha}\rangle$]{\label{lemma: Some properties}}\hfill
    \begin{enumerate}
        \item  $\langle \mathbb{R}^0_\alpha,\leq_{\mathbb{R}^0_\alpha}, \leq^*_{\mathbb{R}^0_\alpha}\rangle$ is of Prikry-type.  
    
        \item $\langle \mathbb{R}^0_\alpha,\leq^*_\alpha\rangle$ is $\alpha$-strategically-closed.

        \item  There is a projection between $\langle \mathbb{R}^0_\alpha,\leq_{\mathbb{R}^0_\alpha}\rangle$ and $\ro(\mathbb{C}_\alpha).$

    \end{enumerate}
\end{lemma}
    
\begin{proof}

(1). This is obvious.

\smallskip

(2). To establish this it suffices to show that $\langle \mathbb{R}^0_\alpha,\leq^*_{\mathbb{R}^0_\alpha}\rangle$ possesses a $\leq^*_{\mathbb{R}^0_\alpha}$-dense subset that is $\alpha$-directed-closed. First, we identify the $\leq^*_{\mathbb{R}^0_\alpha}$-dense set:
\begin{claim}
The set
    $$\{(p,\dot{c})\in \mathbb{R}^0_\alpha\mid p\forces \forall i\in \dom(\dot{c})\, (\max(\dot{c}(i))\in \dot{b}_\alpha)\}$$ is $\leq^*_{\mathbb{R}^0_\alpha}$-dense.
\end{claim}
\begin{proof}[Proof of claim]
 Let $(p,\dot{c})\in \mathbb{R}^0_\alpha$. Since $p$ forces
  $``\alpha$ is regular''  (Theorem~\ref{thm: regularity}), $``\dot{b}_\alpha$ is a club on $\alpha$''
      and $``|\dom(\dot{c})|{<\alpha}$", the maximality principle gives a $\mathbb{Q}^{\alpha^+}_\alpha$-name $\dot{\gamma}$ that $p$ forces to be above every $\dot{c}(i)$ and inside $\dot{b}_\alpha$. Let $\dot{d}$ be the standard name for the condition in $\dot{\mathbb{C}}_\alpha$ with same domain as $\dot{c}$ such that it is forced that $\dot{d}(i)=\dot{c}(i)\cup\{\dot{\gamma}\}$. It follows that $(p,\dot{d})$ is a condition within the set under consideration such that $(p,\dot{d})\leq^*_{\mathbb{R}^0_\alpha} (p,\dot{c})$.
\end{proof}
For the scope of this prove denote by $D$ the above-displayed set. We claim that $\langle D,\leq^*_{\mathbb{R}^0_\alpha}\rangle$ is $\alpha$-directed-closed: Let $\{(\langle t_i,T_i\rangle, \tau_i)\mid i<j\}$ be a directed subset of $D$ of size ${<}\alpha$. Then, the stems $t_i$ are equivalent. This implies that the trees $T_i$, $T_k$ split according to the same measures (see Remark~\ref{remark: same measures}) and therefore $T_i\cap T_{k}$ is a valid tree. As a result, 
$p_j:=\langle t_0, \bigcap_{i<j} T_i\rangle$ provides  a $\leq^*$-lower bound for all the $\langle t_i,T_i\rangle$'s. It follows that
$$p_j\Vdash\forall\xi\in\bigcup_{i<j}\supp(\tau_i)\left(\sup_{i<j}\max(\tau_i(\xi))\in\dot{b}_{\alpha}\right)$$
Since $\dot{b}_{\alpha}$ is forced to contain only $V$-measurable cardinals, $p_j$ forces that for all $\xi\in\bigcup_{i<j}\supp(\tau_i)$, the above supremum can be appended to $\bigcup_{i<j}\tau_i(\xi)$. Letting $\tau$ be a name such that
$$p_j\Vdash\forall\xi\in\bigcup_{i<j}\tau(\xi)=\bigcup_{i<j}\tau_i(j)\cup\sup_{i<j}\max(\tau_i(\xi))$$
we deduce that $(p_j,\tau)$ is a $\leq_{\dR_{\alpha}^0}^*$-lower bound in $D$.

\smallskip

(3). We show that every $\langle \mathbb{R}^0_\alpha,\leq_{\mathbb{R}^0_\alpha}\rangle$-generic $G$ induces a $\ro(\mathbb{C}_\alpha)$-generic.  Once this is established the projection is obtained as follows: Let $\dot{G}^*$ be the standard $\langle \mathbb{R}^0_\alpha,\leq_{\mathbb{R}^0_\alpha}\rangle$-name for the induced $\mathbb{C}_\alpha$-generic.  Standard forcing  facts (see e.g., \cite[\S2]{DPT}) show that the map 
$$a\mapsto \bigwedge\{c\in \ro(\mathbb{C}_\alpha)\mid a\forces_{\langle \mathbb{R}^0_\alpha,\leq_{\mathbb{R}^0_\alpha}\rangle} \check{c}\in \dot{G}^*\}$$
is a projection.

\smallskip

Let $G*H$ be a  $\langle \mathbb{R}^0_\alpha,\leq_{\mathbb{R}^0_\alpha}\rangle$-generic filter and work in $V[G*H]$. Define 
$$G^*:= \{c\in \mathbb{C}_\alpha\mid \exists a\in G\, (a=(p,\check{c}))\}.$$
Easily $G^*$ is a filter. To show genericity, we let$D\s \mathbb{C}_\alpha$ be a dense set and we verify that 
$\{(p,\check{c})\in \mathbb{R}^0_\alpha\mid c\in D\}$ is dense in $\langle \mathbb{R}^0_\alpha,\leq_{\mathbb{R}^0_\alpha}\rangle$. Here the $\alpha$-distributivity of $\mathbb{Q}^{\alpha^+}_\alpha$ (see Proposition~\ref{prop: distributivity}) plays a key role -- indeed, it ensures that $\dC_{\alpha}^{V[G]}=\dC_{\alpha}^V$. As a result, the trivial condition of $\dQ_{\alpha}^{\alpha^+}$ forces $``\check{D}$ is dense in $\dot{\mathbb{C}}_\alpha$". Given $(p,\dot{c})\in \mathbb{R}^0_\alpha$ we can use the previous observation to find an extension of it $(q,\check{d})$ such that $q\forces \check{d}\in D$, as sought.
\end{proof}

\begin{definition}
    Let $\pi\colon\mathbb{R}^0_\alpha\to\ro(\dC_{\alpha})$ denote the map
    $$(p,\dot{c})\mapsto\bigwedge\{b\in\ro(\dC_{\alpha})\;|\;(p,\dot{c})\Vdash_{\langle\mathbb{R}^0_\alpha,\leq_{\mathbb{R}^0_\alpha}\rangle}\check{b}\in\dot{G}^*\}$$
\end{definition}

The following lemma will be helpful later, since it shows in essence that the definition of $\pi$ does not depend too much on the first coordinate.

\begin{lemma}\label{Lemma: Easier description of Projection}
    For every $(p,\dot{c})\in\mathbb{R}^0_\alpha$,
    $$\pi((p,\dot{c}))=\bigwedge\{b\in\ro(\dC_{\alpha})\;|\;p\Vdash_{\dQ_{\alpha}^{\alpha^+}}\dot{c}\leq_{\ro(\dot{\dC}_{\alpha})}\check{b}\}$$
\end{lemma}

\begin{proof}
    Recall that $\dot{G}^*$ was a name for the upward closure of $\pi_2[\Gamma_{\mathbb{R}^0_\alpha}]$. It suffices to show that for every $(p,\dot{c})\in\mathbb{R}^0_\alpha$, $(p,\dot{c})\Vdash\check{b}\in\dot{G}^*$ if and only if $p\Vdash_{\dQ_{\alpha}^{\alpha^+}}\dot{c}\leq_{\ro(\dC_{\alpha})}\check{b}$. The backward direction is immediate. So suppose that $(p,\dot{c})\Vdash_{\mathbb{R}^0_\alpha}\check{b}\in\dot{G}^*$ and suppose toward a contradiction that $p\not\Vdash_{\dQ_{\alpha}^{\alpha^+}}\dot{c}\leq\check{b}$. In particular, there is $p'\leq_{\dQ_{\alpha}^{\alpha^+}}p$ such that $p'\Vdash_{\dQ_{\alpha}^{\alpha^+}}\dot{c}\not\leq_{\ro(\dot{\dC}_{\alpha})}\check{b}$. Letting $\dot{c}'$ be a name for $\dot{c}-\check{b}$, which is forced by $p'$ to be nonzero, we have
    $$p'\Vdash_{\dQ_{\alpha}^{\alpha^+}}\dot{c}'\leq_{\ro(\dot{\dC}_{\alpha})}\dot{c}\wedge\dot{c}'\perp_{\ro(\dot{\dC}_{\alpha})}\check{b}$$
    But then $(p',\dot{c}')\leq_{\mathbb{R}^0_\alpha}(p,\dot{c})$ forces $\check{b}\notin\dot{G}^*$, a contradiction.
\end{proof}

Now we are ready to define our actual forcing $\dR_{\alpha}$.

\begin{definition}\label{def: Gitiksforcing}
    Define the forcing $\dR_{\alpha}$ as follows: $\dR_{\alpha}:=\mathbb{R}^0_\alpha$. For $r,r'\in\dR_{\alpha}$, we let $r'\leq_{\dR_{\alpha}} r$ if $\pi(r')\leq_{\ro(\dC_{\alpha})}\pi(r)$ and we let $r'\leq_{\dR_{\alpha}}^*r'$ if $r'\leq_{\mathbb{R}^0_\alpha}^*r$.
\end{definition}

    Clearly $\leq_{\mathbb{R}^0_\alpha}\s \leq_{\dR_\alpha}$ because $\pi$ is a projection. In particular, $\leq^*_\alpha\s \leq_\alpha$, so it  makes sense to wonder if $\langle \mathbb{R}_\alpha,\leq_{\dR_\alpha},\leq^*_{\dR_\alpha}\rangle$ is of Prikry-type.

\begin{lemma}[Prikry property]{\label{Lemma: Prikry Property for R}}
    Given $a\in \mathbb{R}_{\alpha}$ and $\varphi$ a sentence in the language of forcing of $\langle \mathbb{R}_\alpha,\leq_{\dR_{\alpha}}\rangle$, there is $b\leq_{\dR_{\alpha}}^* a$ such that $b\parallel_{\langle \mathbb{R}_\alpha,\leq_{\dR_\alpha}\rangle}\varphi.$
\end{lemma}
\begin{proof}
    Let $a\in \mathbb{R}_\alpha$ and $\varphi$ be as before. For each $i<2$, let  $$D_i=\{b\in \mathbb{R}_\alpha\mid \pi(b)\forces^i_{\ro(\mathbb{CS}_\alpha)}\varphi\}.$$

    Clearly, $D_0\cup D_1$ is dense. The Prikry property of $\langle \mathbb{R}^0_\alpha,\leq_{\mathbb{R}^0_\alpha}, \leq^*_{\dR_{\alpha}^0}\rangle$  yields a condition $b\leq^*_{\mathbb{R}^0_\alpha} a$ with $$b\forces^i_{\langle \mathbb{R}^0_\alpha,\leq_{\mathbb{R}^0_\alpha}\rangle} \dot{G}\cap \check{D}_0\neq \emptyset.$$ Since $D_0\cup D_1$ is dense, if $b$ forces the negation of the above sentence  then $b$ must force  $``\dot{G}\cap \check{D}_1\neq \emptyset$". Suppose, without loss of generality, that $b\forces_{\mathbb{R}^0_\alpha} \dot{G}\cap \check{D}_0\neq \emptyset$ (the other case is treated analogously).

    We claim that $$\pi(b)\forces_{\ro(\mathbb{C}_\alpha)}\varphi,$$
    which by the definition  of $\leq_{\dR_{\alpha}}$ would yield  $b\forces_{\langle\mathbb{R}_\alpha,\leq_{\dR_{\alpha}}\rangle}\varphi$.

    Otherwise, let $r\leq_{\ro(\mathbb{CS}_\alpha)} \pi(b)$ forcing $\neg \varphi$. Let $b'\leq_{\mathbb{R}^0_\alpha} b$ be such that $\pi(b')\leq_{\ro(\mathbb{CS}_\alpha)} r$. In particular, $\pi(b')$ forces $\neg \varphi$. Since $b'\leq_{\mathbb{R}^0_\alpha} b$ and the latter forces $``\dot{G}\cap \check{D}_0\neq \emptyset"$ in $\mathbb{R}^0_\alpha$,  there must be $b''\in D_0$, compatible with $b$, such that $\pi(b'')$ forces $\varphi$. But this is a contradiction, because $\pi(b)$ and $\pi(b'')$ are compatible, yet they force contradictory information.
\end{proof}

\begin{lemma}[Properties of $\langle \mathbb{R}_\alpha,\leq_{\dR_{\alpha}},\leq^*_{\dR_{\alpha}}\rangle$]\label{lemma: forcing equivalent}\label{lemma: adding clubs}\hfill
    \begin{enumerate}
        \item $\langle \mathbb{R}_\alpha,\leq_{\dR_{\alpha}},\leq^*_{\dR_{\alpha}} \rangle$ is of Prikry-type  and $\langle \mathbb{R}_\alpha,\leq^*_{\dR_{\alpha}}\rangle$ is $\alpha$-strategically-closed.
        \item  $\langle \mathbb{R}_\alpha,\leq_{\dR_{\alpha}},\leq^*_{\dR_{\alpha}} \rangle$ is forcing equivalent to $ \mathbb{C}_\alpha$. 
        \item  $\langle \mathbb{R}_\alpha,\leq_{\dR_{\alpha}},\leq^*_{\dR_{\alpha}} \rangle$ is $\alpha^+$-cc and $\alpha$-distributive. 
    \end{enumerate}
    In particular, forcing with $\langle \mathbb{R}^0_\alpha,\leq_\alpha\rangle$ preserves cardinals and adds a club on $\alpha$ consisting of $V_0$-regulars.
\end{lemma}

\begin{proof}
(1). The Prikry property is precisely Lemma~\ref{Lemma: Prikry Property for R}. The closure of of $\langle\dR_{\alpha},\leq_{\dR_{\alpha}}^*\rangle$ follows from Lemma~\ref{lemma: Some properties}, using the fact that that poset is literally equal to $\langle\mathbb{R}^0_\alpha,\leq_{\mathbb{R}^0_\alpha}^*\rangle$.

\smallskip

(2). Since $  \mathbb{C}_\alpha$ is forcing equivalent to $\ro(\mathbb{C}_\alpha)$ it suffices to check that $\langle \mathbb{R}_\alpha,\leq_{\dR_{\alpha}}\rangle$ is forcing equivalent to the latter. To justify this we show that the map $d\colon \langle\mathbb{R}_\alpha,\leq_{\dR_{\alpha}}\rangle\to \ro(\mathbb{C}_\alpha)$ given by $a\mapsto \pi(a)$  is a dense embedding. Clearly, $d$ is order-preserving and its range is dense in $\ro(\mathbb{C}_\alpha)$ because $\pi$ was a projection. The map also preserves incompatibilities: Suppose that $\pi(a)$ is $\leq_{\ro(\mathbb{C}_\alpha)}$-compatible with $\pi(b)$ and $c\in \ro(\mathbb{C}_\alpha)$  witnessing this. Since $\pi$ was a projection, there is a condition $u\leq_{\mathbb{R}^0_\alpha} a$ such that $\pi(u)\leq_{\ro(\mathbb{C}_\alpha)} c$. Then, by definition, $u\leq_\alpha a, b$.

\smallskip

(3). Since $\langle \mathbb{R}_\alpha,\leq_{\dR_{\alpha}}\rangle$ is forcing equivalent to $\mathbb{C}_\alpha$, if the latter is $\alpha^+$-cc and $\alpha$-distributive then the same would be true for the former. Certainly, $\mathbb{C}_\alpha$ is $\alpha^+$-cc. The $\alpha$-distributivity of $\mathbb{C}_\alpha$ follows from Lemma \ref{Lemma: Abraham-Shelah New} with the proviso that $\mathrm{Reg}^{V_0}\cap \alpha$ is a fat stationary set (see Definition~\ref{def: fat stationary}): Let  $D\s \alpha$  be a club and $\beta<\alpha$. Since $o^{\mathcal{U}}(\alpha)=\alpha^+$, $S=\{\delta<\alpha\mid o^{\mathcal{U}}(\delta)\geq \beta\}$ is stationary in $V_0$. Let $D_0\s D$ be a club in $V_0$ (see Lemma~\ref{lemma: preserving stationaries}). Then, $\acc(D_0)\cap S$ is non-empty and every $\delta$ therein carries a Prikry/Magidor sequence (relative to $V_0$) $b_\delta\s \delta$ with order-type $\geq \beta$. Therefore, $b_\delta\s^* D\cap  (\mathrm{Reg}^{V_0}\cap \alpha)$ is  closed with order-type ${\geq}\beta$. 
\end{proof}

It is easy to show that $\langle \mathbb{R}_\alpha,\leq_{\dR_{\alpha}}\rangle$ contains a $\leq_{\dR_{\alpha}}$-dense weakly homogeneous subset, i.e., $\{(p,\dot{c})\in \mathbb{R}_\alpha\mid \exists c\in \mathbb{C}_\alpha\, (p\forces \check{c}=\dot{c})\}$. However the existence of this set is not enough to carry out the proof of homogeneity of the iteration $\mathbb{P}_{\kappa+1}$. Indeed, the issue is that in extensions in nonstationary support iterations  must use the $\leq^*$-order in all but finitely many coordinates in the support of the condition that is being extended. If at various stages $\alpha\leq \kappa$ with $o^{\mathcal{U}}(\alpha)=\alpha^+$ we allow ourselves to take $\leq$-extensions after $\omega$-many steps there is a worry that the constructed sequence is not longer  a condition in $\mathbb{P}_{\kappa+1}$. This issue is resolved here by taking the ${<}\alpha$-support product of $\alpha^+$-many copies of $\mathbb{CS}_\alpha$. Intuitively, this choice confer us ``enough room'' to shuffle around conditions in the  $\mathbb{C}_\alpha$-part of a condition. 

\begin{remark}[\textbf{Why we use many copies of club shooting}]\label{keyremark}
Let us illustrate why the more natural approach, using a single copy of $\mathbb{CS}_\alpha$ in $\mathbb{R}^0_\alpha$, does not seem to work. Work in $V[G_\alpha]$, where $G_\alpha\subseteq\mathbb{P}_\alpha$ is $V$-generic, and let $\dot\eta_i$ be the standard $\mathbb{Q}^{\alpha^+}_\alpha$-name for the $i$th member of $\dot b_\alpha$. Consider the conditions $a=(\one,\{\dot\eta_0\})$ and $b=(\one,\{\dot\eta_1\})$ of $\mathbb{R}^0_\alpha$. They have no common direct extension: such a condition $(q,\dot c)$ would force $\dot c$ to end-extend both $\{\dot\eta_0\}$ and $\{\dot\eta_1\}$, so $\dot\eta_0=\min(\dot c)=\dot\eta_1$. Thus the Ben-Neria--Unger argument needs direct extensions $a^*\le^*_\alpha a$ and $b^*\le^*_\alpha b$, and an automorphism $\sigma$ of $\mathbb{R}^0_\alpha$ for both $\le_\alpha$ and $\le^*_\alpha$, such that $\sigma(a^*)$ and $b^*$ have a common direct extension.

The first option is $\sigma(q,\dot c)=(q,\Gamma(\dot c))$ for an automorphism $\Gamma\in V[G_\alpha]$ of $\mathbb{CS}_\alpha$. Write $a^*=(q,\dot c_0)$, and let $(r,\dot e)$ be a common direct extension of $\sigma(a^*)=(q,\Gamma(\dot c_0))$ and $b^*$. Then $r$ is a direct extension of $\one$. As $\dot c_0$ end-extends $\{\dot\eta_0\}$, and $\Gamma$ preserves end-extension and fixes $\emptyset$, the set $\Gamma(\{\dot\eta_0\})$ is a nonempty initial segment of $\Gamma(\dot c_0)$, and hence of $\dot e$. Since $\dot e$ also end-extends $\{\dot\eta_1\}$, we get $r\Vdash \min(\Gamma(\{\dot\eta_0\}))=\dot\eta_1$. That is, $r$ forces $\dot\eta_1=F(\dot\eta_0)$, where $F(\nu):=\min(\Gamma(\{\nu\}))$ is a function in $V[G_\alpha]$. No direct extension of $\one$ forces this, because once the first point of $b_\alpha$ is fixed, the second still ranges over a measure-one set.\footnote{Indeed, write $r=\langle\emptyset,T\rangle$, pick $\nu\in\mathrm{Suc}_T(\emptyset)$ with $o^{\mathcal U}(\nu)=0$, and pick distinct $\nu',\nu''\in\mathrm{Suc}_T(\langle\nu\rangle)$ with $o^{\mathcal U}(\nu')=o^{\mathcal U}(\nu'')=0$. Then $r\curvearrowright\langle\nu,\nu'\rangle$ and $r\curvearrowright\langle\nu,\nu''\rangle$ both force $\dot\eta_0=\nu$, but they force $\dot\eta_1=\nu'$ and $\dot\eta_1=\nu''$ respectively.} 

The other option is to let the automorphism of $\mathbb{CS}_\alpha$ depend on the auxiliary generic of $\mathbb{Q}^{\alpha^+}_\alpha$. For instance, let $\dot\Gamma$ name the automorphism of $\mathbb{CS}_\alpha$ that swaps the cones above $\{\dot\eta_0\}$ and $\{\dot\eta_1\}$, and put $\sigma(q,\dot c)=(q,\dot\Gamma(\dot c))$. Then $\sigma$ preserves $\mathbb{R}_\alpha$ and $\le^*_{\dR_{\alpha}}$, and $\sigma(a)=b$. However, $\sigma$ does not preserve $\le_{\dR_{\alpha}}$, which only sees the projection $\pi$ to $\mathrm{RO}^+(\mathbb{CS}_\alpha)$. To see this, fix $\mu<\nu$ with $o^{\mathcal U}(\mu)=o^{\mathcal U}(\nu)=0$, and let $q_0,q_1\in\mathbb{Q}^{\alpha^+}_\alpha$ have stems $\langle\nu\rangle$ and $\langle\mu,\nu\rangle$. The conditions $x=(q_0,\{\check\nu\})$ and $y=(q_1,\{\check\nu\})$ belong to $\mathbb{R}_\alpha$ and $\pi(x)=\pi(y)=\{\nu\}$, so they are $\le_{\dR_\alpha}$-equivalent. But $q_0$ forces $\{\nu\}=\{\dot\eta_0\}$, so $\sigma(x)=(q_0,\{\dot\eta_1\})$, and every possible value of its second coordinate has minimum greater than $\nu$. On the other hand, $q_1$ forces $\{\nu\}=\{\dot\eta_1\}$, so $\sigma(y)=(q_1,\{\check\mu\})$ and $\pi(\sigma(y))=\{\mu\}$. Hence $\sigma(x)$ and $\sigma(y)$ are $\le_{\dR_\alpha}$-incompatible. In general, an automorphism $\sigma$ for $\le_\alpha$ must satisfy that $\pi(\sigma(x))$ depends only on $\pi(x)$, and a $\dot\Gamma$ depending on the auxiliary generic breaks this.

Using $\alpha^+$ coordinates avoids both problems: a ground-model permutation moves $\{\dot\eta_0\}$ to a different coordinate, so no automorphism has to send $\dot\eta_0$ to $\dot\eta_1$. This idea lies at the core of Claim~\ref{Claim: Verifying aleph and beth for order betaplus}.
\end{remark}

\subsection{Relation with the constructions of Gitik and Ben-Neria--Unger}\label{sec:comparison with Gitik}
It may be useful to explain more precisely how the forcing developed here differs from the constructions on which it is based. Our forcing  combines two ingredients: Gitik's method~\cite{GitikOnClosed} for realizing club-shooting through the ground-model regular cardinals as a Prikry-type forcing, and the nonstationary support iteration of Ben-Neria--Unger~\cite{BenUng}, whose homogeneity properties make it possible to obtain a precise computation of $\HOD$. Neither construction, however, can be used verbatim for our purposes.

There are two difficulties on the Gitik side. In the limiting case $o(\alpha)=\alpha$, the forcing of~\cite{GitikOnClosed} is assembled from $\alpha$ many cells of the form
$$\mathbb Q^\beta_\alpha*\mathbb{CS}(\Reg^V\cap\alpha),$$
and the direct-extension order in the $\beta$th cell has only the corresponding $\beta$-closure. Thus one does not obtain uniformly $\alpha$-closed direct extensions at stage $\alpha$. Gitik is able to compensate for this loss of closure by exploiting the Easton support of his iteration; this feature enters essentially, for instance, in the inductive construction of~\cite[Claim~1.15.1]{GitikOnClosed}, where he shows that cardinals are preserved. The same argument is not available for a nonstationary-support iteration, since the support of a condition may be unbounded in regular cardinals and the fusion arguments require simultaneous control of direct extensions along such supports. We therefore redesign the local forcings so that they admit substantially stronger closure in the direct-extension order, and extend the construction to coherent sequences with Mitchell order as high as $o(\alpha)=\alpha^+$(see \S\ref{sec: basic building block}).\footnote{The need to allow high Mitchell orders is not a whim: Mitchell \cite{Mitchellokappa} shows that producing a club of former regulars while retaining measurability of $\kappa$ requires consistency strength $o(\kappa)=\kappa$.}

A different difficulty arises on the Ben-Neria--Unger side. Their nonstationary support machinery provides exactly the global fusion and homogeneity framework that we need, but the local forcings considered in~\cite{BenUng} live in a regime of Mitchell orders strictly below the relevant cardinal. This is not sufficient for our application, where the top large cardinal must survive and where club shooting at the top has to be supported by many corresponding high-order stages below it. We therefore have to incorporate Gitik-type club shooting at stages of large Mitchell order while preserving the Ben-Neria--Unger homogeneity analysis. This is particularly delicate because the most natural local forcing--using a single copy of $\mathbb{CS}(\Reg^V\cap\alpha)$--does not provide enough flexibility for the required automorphism argument (see Remark~\ref{keyremark}). 

Thus the forcing construction developed here is not simply the juxtaposition of the two earlier constructions. We simultaneously strengthen the local Prikry-type forcings and strengthen the homogeneity mechanism so that these more complicated local forcings can still be iterated with nonstationary support. It is this combination that allows us to shoot a club of ground-model measurable cardinals while preserving the top large cardinal--in our application, $\HOD$-supercompactness--and at the same time retain enough homogeneity to identify the $\HOD$ of the final model.

\section{Homogeneity}\label{sec: homogeneity}

In this section we establish that the iteration  $\mathbb{P}_{\kappa+1}$ is almost homogeneous (Theorem~\ref{Theorem: Iteration is Homogeneous}). Our automorphism will be built up by an inductive construction. Commonly, this is done by using names for automorphisms of $\dot{\dQ}_{\beta}$ and requiring that the previously built automorphism of $\dP_{\beta}$ preserves $\dot{\dQ}_{\beta}$ (see \cite{DF, BenUng, BenHay}). However, in our construction, it is vital that $\dot{\dQ}_{\beta}$ is moved by the previously constructed automorphism of $\dP_{\beta}$. Therefore, we will instead consider \emph{isomorphisms} between $\dot{\dQ}_{\beta}$ and the preimage of $\dot{\dQ}_{\beta}$ under the previously constructed automorphism of $\dP_{\beta}$. Furthermore, since we are using a Magidor-style iteration, where the definition of even the regular order on the iteration appeals to the pure extension ordering on the iterands, it is instrumental that these isomorphisms are forced to preserve both $\leq_{\dot{\dQ}_{\alpha}}$ \emph{and} $\leq_{\dot{\dQ}_{\alpha}}^*$. All of these ideas are captured in the following definition, implicit in Ben-Neria–Unger's work \cite{BenUng}:

\begin{definition}\label{Definition: Automorphism Sequence}
    Let $\langle\dP_{\alpha},\dot{\dQ}_{\beta}\;|\;\beta<\alpha\leq\varrho\rangle$ be a nonstationary support iteration of Prikry-type forcings. Let $(\dot{\pi}_{\beta},\dot{D}_{\beta})_{\beta<\rho}$ be a sequence. We define the objects $\mathcal{D}^{\alpha}(\vec{D})$ and $\Gamma_{\vec{\pi}}^{\alpha}$ by induction on $\alpha\leq\rho$, letting $$\text{$\mathcal{D}(\vec{D}):=\mathcal{D}^{\varrho}(\vec{D})$ and $\Gamma_{\vec{\pi}}:=\Gamma_{\vec{\pi}}^{\varrho}$.}\footnote{We will later show inductively that $\mathcal{D}^{\alpha}(\vec{D})$ is an open dense subset of $\dP_{\alpha}$ and $\Gamma_{\vec{\pi}}^{\alpha}$ is an automorphism of $\mathcal{D}^{\alpha}(\vec{D})$ so that the following definition makes sense.}$$

    Suppose $\mathcal{D}^{\beta}(\vec{D})$ and $\Gamma_{\vec{\pi}}^{\beta}$ have been defined for all $\beta<\alpha$.

    \smallskip
    \underline{\textbf{Case $\alpha=0$}:} Let $\mathcal{D}^0(\vec{D}):=\dP_0=\{\emptyset\}$ and $\Gamma_{\vec{\pi}}^0\colon \{\emptyset\}\to \{\emptyset\}$.

    \smallskip
    \underline{\textbf{Case $\alpha=\beta+1$:}} Suppose the following:
    
    \begin{itemize}
        \item $\dot{D}_{\beta}$ is forced to be a $\leq_{\dot{\dQ}_{\beta}}^*$-dense subset of $\dot{\dQ}_{\beta}$
        \item $\dot{\pi}_{\beta}$ is forced to be an isomorphism between $(\dot{D}_{\beta},\leq_{\dot{\dQ}_{\beta}},\leq_{\dot{\dQ}_{\beta}}^*)$ and $((\Gamma^\beta_{\vec{\pi}})^{-1})^*((\dot{D}_{\beta},\leq_{\dot{\dQ}_{\beta}},\leq_{\dot{\dQ}_{\beta}}^*))$ such that $\Vdash\dot{\pi}_{\beta}\circ(\Gamma_{\vec{\pi}}^{\beta})^*(\dot{\pi}_{\beta})=\id_{(\Gamma_{\vec{\pi}}^{\beta})^*(\dot{D}_{\beta})}$.
    \end{itemize}

    If this is not the case, we stop the construction.

    Suppose the two statements above hold. We will suppose that it is forced that $\one_{\dot{\dQ}_{\beta}}\in\dot{D}_{\beta}$. If this is not the case, we simply add it and let $\dot{\pi}_{\beta}$ map $\one_{\dot{\dQ}_{\beta}}$ to itself. Let $\mathcal{D}^{\alpha}(\vec{D})$ consist of all $p\in\dP_{\alpha}=\dP_{\beta}*\dot{\dQ}_{\beta}$ such that $p\restriction\beta\in\mathcal{D}^{\beta}(\vec{D})$ and $p\Vdash p(\beta)\in\dot{D}_{\beta}$. Define $\Gamma_{\vec{\pi}}^{\alpha}$ as follows: For every $p\in\mathcal{D}^{\alpha}(\vec{D})$,
    $$\Gamma_{\vec{\pi}}^{\alpha}(p)=(\Gamma_{\vec{\pi}}^{\beta}(p\restriction\beta),(\Gamma_{\vec{\pi}}^{\beta})^*(\dot{\pi}_{\beta}(p(\beta)))).\footnote{In \cite[Lemma 4]{DF}, the authors define the lifted automorphism the other way around, i.e. as $\dot{\pi}_\beta((\Gamma^\beta_{\vec\pi})^*(\dot{q}))$. However, their case is different, since $(\Gamma^\beta_{\vec\pi})^*$ does not move $\dot{\dQ}_{\beta}$, while in our case it is important that $(\Gamma^\beta_{\vec\pi})^*$ \emph{does} move $\dot{\dQ})_{\beta}$.}$$

    \underline{\textbf{Case $\alpha$ is a limit:}} Let $\mathcal{D}^{\alpha}(\vec{D})$ consist of all conditions $p\in\dP_{\alpha}$ such that $p\restriction\beta\in\mathcal{D}^{\beta}(\vec{D})$ for all $\beta<\alpha$. Define $\Gamma_{\vec{\pi}}^{\alpha}$ by letting, for $p\in\mathcal{D}^{\alpha}(\vec{D})$, $\Gamma_{\vec{\pi}}^{\alpha}(p)$ be such that for every $\beta<\alpha$, $$\Gamma_{\vec{\pi}}^{\alpha}(p)(\beta)=\Gamma_{\vec{\pi}}^{\beta+1}(p\restriction\beta+1)(\beta).$$

    \smallskip

    Finally, we say that $(\dot{\pi}_{\beta},\dot{D}_{\beta})_{\beta<\rho}$ is an \emph{automorphism sequence} if the construction can proceed until $\rho$.
\end{definition}

\begin{remark}\label{Remark: Easy Facts about Auto Sequences}
    It is easy to see that the following holds:
    \begin{enumerate}
        \item $\mathcal{D}(\vec{D})$ consists of all $p\in\dP_{\rho}$ such that for every $\beta<\rho$, $$p\restriction\beta\Vdash p(\beta)\in\dot{D}_{\beta}.$$
        \item The automorphisms $\Gamma_{\vec{\pi}}^{\beta}$ for $\beta\leq\rho$ are coherent in the following sense: Whenever $\beta<\alpha\leq\rho$ and $p\in\mathcal{D}^{\alpha}(\vec{D})$, $$\Gamma_{\vec{\pi}}^{\beta}(p\restriction\beta)=\Gamma_{\vec{\pi}}^{\alpha}(p)\restriction\beta.$$
    \end{enumerate}
\end{remark}
    
    It turns out that every automorphism sequence $(\vec{\pi},\vec{D})$ gives rise to a $\leq^*$-dense subset $\mathcal{D}(\vec{D})$ of $\dP_{\alpha}$ and a self-inverse automorphism of $\mathcal{D}(\vec{D})$ modulo equivalence.
    
    \begin{lemma}\label{lemma: GammaandD}
    Let $(\vec\pi, \vec D)=(\dot{\pi}_{\beta},\dot{D}_{\beta})_{\beta <\rho}$ be an automorphism sequence. Then the following hold:
        \begin{enumerate}
            \item $\mathcal{D}(\vec{D})$ is a $\leq^*$-dense subset of $\dP_{\rho}$;
            \item $\Gamma_{\vec{\pi}}$ is a self-inverse automorphism of $\mathcal{D}(\vec{D})$ modulo equivalence with respect to both $\leq$ and $\leq^*$.
        \end{enumerate}
    \end{lemma}

    \begin{proof}
        We prove the statements one-by-one.

        For (1), let $p\in\dP_{\rho}$. Using the maximum principle, let $q\in\dP_{\rho}$ be a condition with $\supp(q)=\supp(p)$ such that, for every $\beta\in\supp(p)$,
        $$\Vdash_{\dP_{\beta}}q(\beta)\in\dot{D}_{\beta}\wedge q(\beta)\leq_{\dot{\dQ}_{\beta}}^*p(\beta).$$
        Clearly, $q\leq_{\dP_{\rho}}^*p$ and $q\in\mathcal{D}(\vec{D})$.

\smallskip

        To show (2), we prove by induction on $\alpha\leq\rho$ that $\Gamma_{\vec{\pi}}^{\alpha}$ is a self-inverse automorphism of $\mathcal{D}^{\alpha}(\vec{D})$ modulo equivalence. Suppose the statement holds for all $\beta<\alpha$.
        
        \underline{\textbf{Case $\alpha=0$:}} This is easy.

        \underline{\textbf{Case $\alpha=\beta+1$:}} For notational simplicity, we write $\Gamma:=\Gamma_{\vec{\pi}}^{\beta'}$.

        \begin{enumerate}
            \item \textbf{$\Gamma_{\vec{\pi}}^{\alpha}$ preserves the orders:} Suppose that $p',p\in\mathcal{D}(\vec{D})$ and $p'\leq_{\dP_{\alpha}}p$, i.e. $p'\restriction\beta\leq_{\dP_{\beta}}p\restriction\beta$ and $p'\restriction\beta\Vdash_{\dP_{\beta'}}p'(\beta)\leq_{\dot{D}_{\beta}}p(\beta)$. It follows that $\Gamma(p'\restriction\beta)\leq_{\dP_{\beta}}\Gamma(p\restriction\beta)$ by the inductive hypothesis. By the hypothesis on $\dot{\pi}_{\beta}$, we have $p'\restriction\beta\Vdash\dot{\pi}_{\beta}(p'(\beta)\leq_{(\Gamma^{-1})^*(\dot{D}_{\beta})}\dot{\pi}_{\beta}(p(\beta))$. But then we have $\Gamma(p'\restriction\beta)\Vdash\Gamma^*(\dot{\pi}_{\beta}(p'(\beta)))\leq_{\dot{D}_{\beta}}\Gamma^*(\dot{\pi}_{\beta}(p(\beta)))$ as desired.

            The same argument applies when replacing $\leq_{\dot{\dQ}_{\beta}}$ with $\leq_{\dot{\dQ}_{\beta}}^*$.

            \item \textbf{$\Gamma_{\vec{\pi}}^{\alpha}$ maps from $\mathcal{D}^{\alpha}(\vec{D})$ into $\mathcal{D}^{\alpha}(\vec{D})$:} Suppose that $p\in\mathcal{D}^{\alpha}(\vec{D})$. Thus, $p\restriction\beta\in\mathcal{D}^{\beta}(\vec{D})$, implying $\Gamma(p\restriction\beta)\in\mathcal{D}^{\beta}(\vec{D})$ by the inductive assumption. Furthermore, $p\restriction\beta\Vdash p(\beta)\in\dot{D}_{\beta}$. By the assumption, $p\Vdash\dot{\pi}_{\beta}(p(\beta))\in(\Gamma^{-1})^*(\dot{D}_{\beta})$. But then $\Gamma(p\restriction\beta)\Vdash\Gamma^*(\dot{\pi}_{\beta}(\dot{q}))\in\dot{D}_{\beta}$. This easily implies that $\Gamma_{\vec{\pi}}^{\alpha}(p)\in\dot{D}^{\alpha}(\vec{D})$.

            \item \textbf{$\Gamma_{\vec{\pi}}^{\alpha}$ is self-inverse modulo equivalence:} Let $p\in\mathcal{D}^{\alpha}(\vec{D})$. By the inductive assumption,
            $$p\restriction\beta\leq_{\dP_{\beta}}\Gamma^2(p\restriction\beta)\leq_{\dP_{\beta}}p\restriction\beta$$
            It follows that $\Gamma_{\vec{\pi}}^{\alpha}(p)\restriction\beta=\Gamma(p\restriction\beta)$ and $\Gamma_{\vec{\pi}}^{\alpha}(p)(\beta)=\Gamma^*(\dot{\pi}_{\beta}(\Gamma^*(\dot{\pi}_{\beta}(p(\beta)))))$. Note that $\Vdash\Gamma^*(\dot{\pi}_{\beta}(p(\beta)))=\Gamma^*(\dot{\pi}_{\beta})(\Gamma^*(p(\beta)))$. Hence, by the assumptions placed on $\dot{\pi}_{\beta}$, we have that
            $$\Vdash\dot{\pi}_{\beta}(\Gamma^*(\dot{\pi}_{\beta}(p(\beta))))=\dot{\pi}_{\beta}(\Gamma^*(\dot{\pi}_{\beta})(\Gamma^*(p(\beta))))=\Gamma^*(p(\beta))$$
            This clearly implies that $\Vdash\Gamma^*(\dot{\pi}_{\beta'}(\Gamma^*(\dot{\pi}_{\beta'}(p(\beta)))))=\Gamma^*(\Gamma^*(p(\beta)))=p(\beta)$ which implies the desired statement.
            
            \item \textbf{$\Gamma_{\vec{\pi}}^{\alpha}$ is an automorphism modulo equivalence:} This follows from the fact that $\Gamma_{\vec{\pi}}^{\alpha}$ is order-preserving and self-inverse modulo equivalence.
        \end{enumerate}

        \underline{\textbf{Case $\alpha$ is a limit:}} Suppose the statement holds for all $\beta<\alpha$.

        \begin{enumerate}
            \item \textbf{$\Gamma_{\vec{\pi}}^{\alpha}$ preserves the orders:} Suppose that $p,p'\in\mathcal{D}^{\alpha}(\vec{D})$ and $p'\leq_{\dP_{\alpha}}p$. Thus, for every $\beta<\alpha$, $p'\restriction\beta\Vdash p'(\beta)\leq_{\dot{\dQ}_{\beta}}p(\beta)$ and there is $b\subseteq\supp(p')$ finite such that for all $\beta\in\supp(p')\smallsetminus b$, $p'\restriction\beta\Vdash p'(\beta)\leq_{\dot{\dQ}_{\beta}}p(\beta)$. It follows just as in the previous case that this implies that $\Gamma_{\vec{\pi}}^{\alpha}(p')\leq_{\dP_{\alpha}}\Gamma_{\vec{\pi}}^{\alpha}(p)$. The same applies to $\leq_{\dP_{\alpha}}^*$.

            \item \textbf{$\Gamma_{\vec{\pi}}^{\alpha}$ maps from $\mathcal{D}^{\alpha}(\vec{D})$ into $\mathcal{D}^{\alpha}(\vec{D})$:} We have that $p\in\mathcal{D}^{\alpha}(\vec{D})$ if and only if $p\in\dP_{\alpha}$ and $p\restriction\beta\in\mathcal{D}^{\beta}(\vec{D})$. By Remark \ref{Remark: Easy Facts about Auto Sequences}(2), we are only tasked with showing that $\Gamma_{\vec{\pi}}^{\alpha}(p)\in\dP_{\alpha}$. But it is clear from the construction that $\supp(\Gamma_{\vec{\pi}}^{\alpha}(p))=\supp(p)$, which implies the desired statement.

            \item \textbf{$\Gamma_{\vec{\pi}}^{\alpha}$ is self-inverse modulo equivalence:} This follows easily from the inductive assumption, using Remark \ref{Remark: Easy Facts about Auto Sequences}(2).

            \item \textbf{$\Gamma_{\vec{\pi}}^{\alpha}$ is an automorphism modulo equivalence:} This again follows easily from the fact that $\Gamma_{\vec{\pi}}^{\alpha}$ is order-preserving and self-inverse modulo equivalence.
        \end{enumerate}

        This finishes the proof.
    \end{proof}
    Since $\mathcal{D}(\vec{D})$ is $\leq^*$-dense in $\dP_{\rho}$, Lemma \ref{Lemma: Self-Inverse Auto mod Equiv gives Generic} easily implies:

    \begin{cor}\label{Corollary: Automorphism Sequence gives Generic}
        Let $(\dot{\pi}_{\beta},\dot{D}_{\beta})_{\beta<\alpha}$ be an automorphism sequence. Let $G$ be a $\dP_{\alpha}$-generic filter. Then the upward closure of $\Gamma_{\vec{\pi}}[G\cap\mathcal{D}(\vec{D})]$ defines a $\dP_{\alpha}$-generic filter. \qed
    \end{cor}
    Let us explain the construction behind the homogeneity. It is easy to see (and classical) that the usual Prikry-forcing (even tree-Prikry-forcing) is cone homogeneous: Given conditions $\langle t_0,T_0\rangle$ and $\langle t_1,T_1\rangle$, we  take $T^*:=(T_0\cap T_1)\setminus \max\{\max(t_0),\max(t_1)\}$ and define an isomorphims $\gamma$ between $\mathbb{Q}/\langle t_0, T^*\rangle$ and $\mathbb{Q}/\langle t_1, T^*\rangle$  taking $\gamma(\langle s, S\rangle):=\langle t_1{}^\smallfrown (s\setminus \max(t_1)+1), S\rangle)$.
    
    However, in the case of the poset $\dQ_{\alpha}^{\rho}$, it is not that simple: Given conditions $\langle t_0,T_0\rangle$ and $\langle t_1,T_1\rangle$, it is the case that  $\Suc_{T_0}(t_0)$ is large for $\mathcal{U}_{\alpha,\rho}(t_0)$ and $\Suc_{T_1}(t_1)$ is large for $\mathcal{U}_{\alpha,\rho}(t_1)$. If $t_0\neq t_1$, these two ultrafilters necessarily concentrate on disjoint sets: One concentrates on those ordinals whose Prikry/Magidor sequence begins with $t_0\restriction\rho$, the other one concentrates on those ordinals whose Prikry/Magidor sequence begins with $t_1\restriction\rho$. In fact, it can be shown that for $\rho$ a limit ordinal of uncountable cofinality, $\dQ_{\alpha}^{\rho}$ is not homogeneous.\footnote{An easy way to show this is  by using work of Hayut--Karagila \cite{HayKag}. Combining their Corollary~12 and Theorem 13, it follows that any weakly homogeneous  poset changing the cofinality of $\kappa$ to a regular cardinal $\mu>\omega$, without adding bounded subsets of $\kappa$, must change the cofinality of some $\lambda <\kappa$. Thus, if $\mathbb{Q}^\rho_\alpha$ were to be weakly homogeneous then it would have to change cofinalities below $\kappa$, which is not the case.} However, here is where the prior iteration $\mathbb{P}_\alpha$ comes in: According to the definition of an automorphism sequence, given a $\dP_{\alpha}$-generic filter $G_{\alpha}$ and letting $H_{\alpha}:=\Gamma_{\vec{\pi}}^{\alpha}[G_{\alpha}]$, we are not aiming for an automorphism of $(\dQ_{\alpha}^{\rho})^{G_\alpha}$, but instead for an \emph{isomorphism} of dense subsets of $(\dQ_{\alpha}^{\rho})^{G_\alpha}$ and $(\dot{\dQ}_{\alpha}^{\rho})^{H_{\alpha}}$.\footnote{It is easy to see that $(\Gamma_{\vec{\pi}}^{\alpha,-1})^*(\dot{\dQ}_{\alpha}^{\rho})^{G_{\alpha}}=(\dot{\dQ}_{\alpha}^{\rho})^{H_{\alpha}}$} Thus, by choosing our prior automorphism sequence correctly, for fixed conditions $\langle t_0,T_0\rangle$ and $\langle t_1,T_1\rangle$, we can make sure that, if $\langle t_0,T_0\rangle\in(\dQ_{\alpha}^{\rho})^{G_\alpha}$ and $\langle t_1,T_1\rangle\in(\dot{\dQ}_{\alpha}^{\rho})^{H_{\alpha}}$, $T_0$ splits according to the sequence of measures $$\langle \dot{\mathcal{U}}_{\alpha,\sigma}^{H_{\alpha}}((t_1{}^\smallfrown \vec\nu)\restriction \sigma)\mid \sigma<\rho\rangle$$ and vice versa, $T_1$ splits according to 
    $$\langle \dot{\mathcal{U}}_{\alpha,\sigma}^{G_{\alpha}}((t_0{}^\smallfrown \vec\nu)\restriction \sigma)\mid \sigma<\rho\rangle.$$
    
   This is achieved by swapping $t_0$ and $t_1$ \emph{below $\alpha$} so that whenever the Prikry/Magidor sequence of $\gamma$ according to $G_{\beta}$ begins with $t_0$, the Prikry sequence of $\gamma$ according to $H_{\beta}$ begins with $t_1$. Thanks to the fractal-like structure of the poset, this amounts to the following task: We want to perform the same automorphism at $\alpha$ in $V$ and at $\alpha$ in any $M_{\alpha,\rho}$. The easiest way to accomplish this is to perform an automorphism representing the one at $\alpha$ at a club of points below $\alpha$. However, in any condition in $\dP_{\kappa+1}$, we do not have to worry about just one condition, but unboundedly many. Fortunately, since the support of any condition is nowhere-stationary, we can partition the ordinals below $\kappa$ into disjoint clubs for any $\alpha$ in the support of a given condition in such a way that this aligns nicely with the relevant ultrapowers. To ensure that this is possible, we prove a lemma giving us \emph{coherent representations} for small sets which also align nicely with the relevant ultrapowers, allowing us to define an automorphism sequence such that, for every $(\alpha,\sigma)(\mathcal{U})$, the same automorphism is performed at $\alpha$ in $V[G_{\alpha}]$ and in $M_{\alpha,\sigma}[G_{\alpha}]$.

   \smallskip

   The first step is showing that nowhere-stationary sets admit disjoint ladder systems satisfying a strong zero-one law. Note that nowhere-stationarity here is crucial: If $A$ where stationary in any regular cardinal, we could apply Fodor's lemma to $\alpha\mapsto\min(C_{\alpha})$ which would show that the $C_{\alpha}$ cannot be pairwise disjoint.

\begin{lemma}\label{Lemma: C-Sequence for NWS}
    Suppose that $\delta$ is a regular cardinal, $\gamma<\delta$ is an ordinal and $A\subseteq\delta$ is a nowhere-stationary set consisting of regular cardinals such that $A\cap(\gamma+1)=\emptyset$. Then, there exists a sequence $(C_{\alpha})_{\alpha\in A}$ such that:
    \begin{enumerate}
        \item For every $\alpha\in A$, $C_{\alpha}\subseteq\alpha$ is club and $C_{\alpha}\cap\gamma=\emptyset$;
        \item For every distinct $\alpha,\beta\in A$, $(C_{\alpha}\cup\{\alpha\})\cap (C_{\beta}\cup\{\beta\})=\emptyset$;
        \item For every regular cardinal $\mu\leq\delta$, either $(\bigcup_{\alpha\in A}C_{\alpha})\cap \mu$ is nonstationary in $\mu$ or there is $\alpha\in A$ such that $C_{\alpha}\cap \mu$ is club in $\mu$. Moreover, for any $\alpha\in\delta$, if $\mu\in C_{\alpha}$ then $C_{\alpha}\cap\mu$ is club in $\mu$.
    \end{enumerate}
\end{lemma}

\begin{proof}
    We prove the lemma by induction on $\delta$.
    
\smallskip
    \underline{\textbf{Case  $\delta=\theta^+$:}} Let $\gamma<\delta$ be an ordinal and $A\subseteq\delta$ a nowhere-stationary set of regular cardinals. If $\theta\notin A$ then the induction hypothesis hands us a sequence $(C_\alpha)_{\alpha\in A\cap \theta}=(C_\alpha)_{\alpha\in A}$ witnessing Clauses~(1)--(3) above.

    So, let us assume that $\theta\in A$. Let $C_{\theta}\subseteq\theta$ be a club set disjoint from $A\cup\gamma$. We may assume that for every $\gamma\in C_\theta$ with $\cf(\gamma)>\omega$, $C_\theta\cap\gamma$ is unbounded in $\gamma$. In particular, for every $\alpha\in A$, $C_\theta\cap\alpha$ is bounded in $\alpha$. Fix a sequence $(D_\alpha)_{\alpha\in A\cap\theta}$ witnessing the inductive hypothesis for $A\cap\theta$ and $\gamma$. For $\alpha\in A\cap\theta$, let $C_{\alpha}:=D_\alpha\smallsetminus\sup(C_{\alpha}\cap\theta)+1$.

    \begin{claim}
     $(C_\alpha)_{\alpha\in A}$ satisfies the conclusion of the lemma.   
    \end{claim}
    \begin{proof}[Proof of claim]
     
     Both (1) and
    and (2) hold by the inductive hypothesis and the definition of $C_{\alpha}$. To see (3), let $\mu\leq\delta$ be a regular cardinal. Suppose that $(\bigcup_{\alpha\in A}C_{\alpha})\cap \mu$ is stationary in $\mu$. If $C_{\theta}\cap\mu$ is stationary in $\mu$, we have that $C_{\theta}\cap\mu$ is club in $\mu$ -- by assumption. So suppose that $\left(\bigcup_{\alpha\in A\cap\theta}C_{\alpha}\right)\cap\mu$ is stationary in $\mu$. In particular, $\left(\bigcup_{\alpha\in A\cap\theta}D_\alpha\right)\cap\mu$ is stationary in $\mu$. Additionally, we must have $\mu\leq\theta$. Therefore, there exists $\alpha\in A\cap\theta$ such that $D_\alpha\cap\mu$ is club in $\mu$. Note that $D_\alpha$ witnesses that $(\bigcup_{\beta\in (A\cap\theta)\smallsetminus\{\alpha\}}D_\alpha)\cap\mu$ and therefore also $$\left(\bigcup_{\beta\in (A\cap\theta)\smallsetminus\{\alpha\}}C_{\alpha}\right)\cap\mu$$ is nonstationary in $\mu$. Thus, $C_{\alpha}\cap\mu$ must be stationary in $\mu$, since we were assuming that $(\left(\bigcup_{\beta\in (A\cap\theta)\smallsetminus\{\alpha\}}C_{\beta}\right)\cup C_{\alpha})\cap\mu$ is stationary in $\mu$. But this implies that $\mu\in C_{\alpha}$ and so $C_{\alpha}\cap\mu$ is actually club in $\mu$.
    \end{proof}
    
    \underline{\textbf{Case $\delta$ is a limit:}} As before, find $C\subseteq\delta$ club such that $C\cap A=\emptyset$. Enumerate $C=\langle \delta_i\;|\;i<\delta\rangle$, where $\gamma=\delta_0$. Thus, $\delta_i\notin A$ for all $i<\delta$. For each $i<\delta$, apply the inductive hypothesis to $A\cap(\delta_i,\delta_{i+1})$ and $\delta_i$ to find a sequence $(C_{\alpha})_{\alpha\in A\cap(\delta_i,\delta_{i+1})}$ witnessing the inductive hypothesis for $A\cap(\delta_i,\delta_{i+1})$ and $\delta_i$, respectively. Let us show that $(C_{\alpha})_{\alpha\in A}$ fulfills (1)--(3).

\smallskip
    
    Clause (1) is clear.  For Clause~(2), let $\alpha<\beta$ be elements of $A$. If there is $i<\delta$ such that $\alpha,\beta\in (\delta_i,\delta_{i+1})$ then  $(C_{\alpha}\cup\{\alpha\})\cap (C_{\beta}\cup\{\beta\})=\emptyset$ because of the properties of the sequence $(C_{\alpha})_{\alpha\in A\cap(\delta_i,\delta_{i+1})}$. Otherwise, there is $i<\delta$ such that $\alpha<\delta_i<\beta$. In this case, we have $C_{\alpha}\subseteq\alpha<\delta_i$ and $C_{\beta}\cap\delta_i=\emptyset$, so $$(C_{\alpha}\cup\{\alpha\})\cap (C_{\beta}\cup\{\beta\})=\emptyset.$$
    
    For (3), let $\mu\leq\delta$ be a regular cardinal. For $\mu=\delta$, clearly, $C$ is club in $\delta$ and disjoint from $\bigcup_{\alpha\in A}C_{\alpha}$, ergo the latter is nonstationary. So, suppose that $\mu<\delta$. If $C\cap\mu$ is unbounded in $\mu$, it must be club in $\mu$ and, again, it is disjoint from $(\bigcup_{\alpha\in A}C_{\alpha})\cap \mu$. Otherwise, there is $i<\delta$ such that $\delta_i<\mu<\delta_{i+1}$. If $(\bigcup_{\alpha\in A}C_{\alpha})\cap \mu$ is stationary in $\mu$, $\bigcup_{\alpha\in A\cap(\delta_i,\delta_{i+1})}C_{\alpha}$ must be stationary in $\mu$ as well, for $\bigcup_{\alpha\in A\cap\delta_i}C_{\alpha}\subseteq\delta_i<\mu$ and $(\bigcup_{\alpha\in A\cap(\delta_{i+1},\delta)}C_{\alpha})\cap\delta_{i+1}=\emptyset$. It follows from the inductive hypothesis that there is $\alpha\in A\cap(\delta_i,\delta_{i+1})$ such that $C_{\alpha}\cap \mu$ is club in $\mu$, as sought.
\end{proof}

The second ingredient is the following lemma, stating that, for $\alpha\in\dom_1(\mathcal{U})$, sets in $V_{\alpha+1}$ can be represented using functions which cohere with respect to the ultrapowers according to $\mathcal{U}_{\beta,\sigma}$ for $\beta\leq\alpha$.

\begin{lemma}\label{Lemma: Representing small sets}
    Suppose $\alpha\in\dom_1(\mathcal{U})$ and $x\in V_{\alpha+1}$. Then, there is a function $h\colon \alpha\to V$ such that the following holds:
    \begin{enumerate}
        \item For every $\sigma<o^{\mathcal{U}}(\alpha)$, $j_{\alpha,\sigma}(h)(\alpha)=x$;
        \item For every $(\beta,\sigma)\in\dom(\mathcal{U})$, $j_{\beta,\sigma}(h)(\beta)=h(\beta)$.
    \end{enumerate}
\end{lemma}

\begin{proof}
    Since $x\in V_{\alpha+1}\s M_{\alpha,0}$, there is a function $f\colon\alpha\to V_{\alpha}$ such that $$j_{\alpha,0}(f)(\alpha)=x.$$ Let $g\colon\alpha\to V$ denote the following function:
    $$g(\beta):=\begin{cases}
        \emptyset, & o^{\mathcal{U}}(\beta)=0; \\
        j_{\beta,0}(f\restriction\beta)(\beta), & o^{\mathcal{U}}(\beta)>0.
    \end{cases}$$

    Let us show that $g$ satisfies the first requirement:
    
    \begin{claim}
        For every $0<\sigma<o^{\mathcal{U}}(\alpha)$, $$j_{\alpha,\sigma}(g)(\alpha)=x.$$
    \end{claim}
    \begin{proof}
        By coherency, $j_{\alpha,\sigma}(\mathcal{U})_{\alpha,0}=\mathcal{U}_{\alpha,0}$. Ergo,
        $$(j_{{j_{\alpha,\sigma}(\mathcal{U})_{\alpha,0}}})^{M_{\alpha,\sigma}}=(j_{\mathcal{U}_{\alpha,0}})^{M_{\alpha,\sigma}}=j_{\alpha,0}\restriction M_{\alpha,\sigma}.$$

        Therefore, we have
        $$j_{\alpha,\sigma}(g)(\alpha)=(j_{{j_{\alpha,\sigma}(\mathcal{U})_{\alpha,0}}})^{M_{\alpha,\sigma}}(j_{\alpha,\sigma}(f)\restriction\alpha)(\alpha)=j_{\alpha,0}(f)(\alpha)=x.$$
       The second of these equalities follows from two observations: First, we have  $j_{\alpha,\sigma}(f)\restriction\alpha = f$ as $\mathrm{range}(f)\s V_\alpha$ and $\crit(j_{\alpha,\sigma})=\alpha$. Second, the function $f\colon \alpha\to V_\alpha$ belongs to $M_{\alpha,\sigma}$ (in that this ultrapower is closed under  $\alpha$-sequences of its members) and so $(j_{{j_{\alpha,\sigma}(\mathcal{U})_{\alpha,0}}})^{M_{\alpha,\sigma}}(f)=j_{\alpha,0}(f)$.
    \end{proof}
    
    \begin{claim}\label{claim: representation}
        For every $\beta<\alpha$ with $o^{\mathcal{U}}(\beta)>0$ and every $0<\sigma<o^{\mathcal{U}}(\beta)$, $$j_{\beta,\sigma}(g\restriction\beta)(\beta)=j_{\beta,0}(f\restriction\beta)(\beta).$$
    \end{claim}

    \begin{proof}
      By the same argument as in the previous claim, 
        $$j_{\beta,\sigma}(g\restriction\beta)(\beta)=(j_{j_{\beta,\sigma}(\mathcal{U})_{\beta,0}})^{M_{\beta,\sigma}}(j_{\beta,\sigma}(f\restriction\beta)\restriction\beta)=j_{\beta,0}(f\restriction\beta)(\beta).$$
    \end{proof}

    So we can let $h\colon \alpha\to V$ be a function defined as
    $$h(\beta):=\begin{cases}
        f(\beta), & o^{\mathcal{U}}(\beta)=0;\\
        g(\beta), & o^{\mathcal{U}}(\beta)>0.
    \end{cases}$$

    Now we can show that $h$ is as required by proving the following two claims:

    \begin{claim}
        For every $\sigma<o^{\mathcal{U}}(\alpha)$, $j_{\alpha,\sigma}(h)(\alpha)=x$.
    \end{claim}

    \begin{proof}
        If $\sigma=0$ then $j_{\alpha,0}(h)(\alpha)=j_{\alpha,0}(f)(\alpha)=x$ (because, by coherency, $o^{j_{\alpha,0}(\mathcal{U})}(\alpha)=0$). Similarly, for  $\sigma>0$, $j_{\alpha,\sigma}(h)(\alpha)=j_{\alpha,\sigma}(g)(\alpha)=x$.
    \end{proof}

    \begin{claim}
        For every $\beta<\alpha$ with $o^{\mathcal{U}}(\beta)>0$ and every $\sigma<o^{\mathcal{U}}(\beta)$, $$j_{\beta,\sigma}(h\restriction\beta)(\beta)=h(\beta).$$
    \end{claim}

    \begin{proof}
        For $o^{\mathcal{U}}(\beta)>0$ and $\sigma<o^{\mathcal{U}}(\beta)$, the claim follows from the definition of $h$ and because  $o^{j_{\beta,0}(\mathcal{U})}(\beta)=0$. Indeed:
        $$h(\beta)=g(\beta)=j_{\beta,0}(f\restriction\beta)(\beta)=j_{\beta,0}(h\restriction\beta)(\beta).$$
        Suppose $\sigma>0$, then combining the Claim~\ref{claim: representation} with the definition of $h$:
        $$j_{\beta,\sigma}(h\restriction\beta)(\beta)=j_{\beta,\sigma}(g\restriction\beta)(\beta)=j_{\beta,0}(f\restriction\beta)(\beta)=g(\beta)=h(\beta).$$
        This completes the proof of the claim.
    \end{proof}
    We are done with the proof of the lemma.
\end{proof}

To finish the preparations, we combine the previous two lemmas to show that we can represent sequences indexed by nowhere-stationary subsets of $\dom_1(\mathcal{U})$ in a way that coheres with the ultrapowers by \emph{every} measure $\mathcal{U}_{\beta,\sigma}$ where $(\beta,\sigma)\in\dom(\mathcal{U})$.

\begin{lemma}\label{Lemma: Universal Representations}
    Let $A\subseteq\dom_1(\mathcal{U})$ be nowhere-stationary. Let $(t_{\alpha})_{\alpha\in A}$ be a sequence such that for every $\alpha\in A$, $t_{\alpha}\in V_{\alpha+1}$. Then there is a function $$h\colon\kappa+1\to V_{\kappa+1}$$ such that the following holds:

    \begin{enumerate}
        \item For every $\alpha\in A$, $h(\alpha)=t_{\alpha}$;
        \item For every $\beta\in\dom_1(\mathcal{U})$ and $\sigma\leq o^{\mathcal{U}}(\beta)$, $j_{\beta,\sigma}(h)(\beta)=h(\beta)$.
    \end{enumerate}
\end{lemma}

\begin{proof}
    Fix a sequence $(C_{\alpha})_{\alpha\in A}$ exemplifying Lemma \ref{Lemma: C-Sequence for NWS}, i.e. for every $\alpha\in A$, $C_{\alpha}\subseteq\alpha$ is club, for every distinct $\alpha,\beta\in A$, $(C_{\alpha}\cup\{\alpha\})\cap(C_{\beta}\cup\{\beta\})=\emptyset$ and for every regular cardinal $\mu\leq\kappa$, either $\left(\bigcup_{\alpha\in A}C_{\alpha}\right)\cap\mu$ is nonstationary in $\mu$ or there is $\alpha\in A$ such that $C_{\alpha}\cap\mu$ is club in $\mu$. Moreover, for any $\alpha\in\delta$, if $\mu\in C_{\alpha}$, then $C_{\alpha}\cap\mu$ is club in $\mu$.
    
    For every $\beta\in\kappa+1$, let $g(\beta)$ be the unique $\alpha$ with $\beta\in C_{\alpha}\cup\{\alpha\}$ if it exists, $0$ otherwise. Also fix, for every $\alpha\in A$, a function $f_{\alpha}\colon\kappa+1\to V_{\kappa+1}$ such that $f_{\alpha}(\alpha)=t_{\alpha}$ and for every $\beta\in\dom_1(\mathcal{U})$ and $\sigma<o^{\mathcal{U}}(\beta)$,
    $$j_{\beta,\sigma}(f_{\alpha})(\beta)=f_{\alpha}(\beta)$$
    
    This is possible by Lemma \ref{Lemma: Representing small sets}: Let $f_{\alpha}'\colon\alpha\to V_{\alpha+1}$ be a function exemplifying Lemma \ref{Lemma: Representing small sets} for $x=t_{\alpha}$. Define $f_{\alpha}\colon\kappa+1\to V_{\alpha+1}$ by letting $f_{\alpha}(\alpha):=t_{\alpha}$ and $f_{\alpha}\restriction(\alpha,\kappa+1)\equiv\emptyset$. For $\beta\in\dom_1(\mathcal{U})\cap\alpha$ and $\sigma<o^{\mathcal{U}}(\beta)$, we have
    $$j_{\beta,\sigma}(f_{\alpha})(\beta)=(j_{\alpha,\sigma}(f_{\alpha})\restriction j_{\beta,\sigma}(\alpha))=j_{\beta,\sigma}(f_{\alpha}\restriction\alpha)(\beta)=j_{\beta,\sigma}(f_{\alpha}')(\beta)=f_{\alpha}'(\beta)=f_{\alpha}(\beta).$$

    For $\beta=\alpha$ and $\sigma<o^{\mathcal{U}}(\beta)$, we have
    $$j_{\beta,\sigma}(f_{\alpha})(\beta)=t_{\alpha}=f_{\alpha}(\alpha)$$

    For $\beta\in\dom_1(\mathcal{U})\cap(\alpha,\kappa+1)$, we have
    $$j_{\beta,\sigma}(f_{\alpha})(\beta)=\emptyset=f_{\alpha}(\beta)$$
    where the first equality follows from the fact that $j_{\beta,\sigma}(\alpha)<j_{\beta,\sigma}(\beta)$ by elementarity.
    
    Now we define
    $$h(\beta):=\begin{cases}
        f_{g(\beta)}(\beta), & g(\beta)\neq 0;\\
        \emptyset, & g(\beta) = 0.
    \end{cases}$$
    Let us show that this works. Clearly, for $\alpha\in A$, $$h(\alpha)=f_{g(\alpha)}(\alpha)=f_{\alpha}(\alpha)=t_{\alpha}.$$ Now suppose that $(\beta,\sigma)\in\dom(\mathcal{U})$.
    
    First suppose that $g(\beta)=0$. Thus, $\left(\bigcup_{\alpha\in A}C_{\alpha}\right)\cap\beta$ is nonstationary in $\beta$ and so, since $A\cap\beta$ is nonstationary in $\beta$ as well,
    $$\beta\notin j_{\beta,\sigma}\left(\bigcup_{\alpha\in A}C_{\alpha}\cup\{\alpha\}\right)$$
    Therefore, $j_{\beta,\sigma}(g)(\beta)=0$, implying $j_{\beta,\sigma}(h)(\beta)=0=h(\beta)$.

    Suppose now that $g(\beta)=\alpha\neq 0$. Therefore, $C_{\alpha}\cap\beta$ is club in $\beta$ and thus $\beta\in j_{\beta,\sigma}(C_{\alpha})$. Since $j_{\beta,\sigma}(C_{\alpha})=j_{\beta,\sigma}((C_{\alpha})_{\alpha\in A})(j_{\beta,\sigma}(\alpha))$, this shows $j_{\beta,\sigma}(g)(\beta)=j_{\beta,\sigma}(\alpha)$. Now we have
    $$j_{\beta,\sigma}(h)(\beta)=j_{\beta,\sigma}((f_{\alpha})_{\alpha\in A})_{(j_{\beta,\sigma}(\alpha))}(\beta)=j_{\beta,\sigma}(f_{\alpha})(\beta){=f_{\alpha}(\beta)}=h(\beta)$$
    by assumption.
\end{proof}

Now that we are done with the preparations, we can move on to the prove of the main result of this section. 

\begin{theorem}\label{Theorem: Iteration is Homogeneous}
    $\dP_{\kappa+1}$ is almost homogeneous modulo equivalence.
\end{theorem}

\begin{proof}
    Let $q_0,q_1\in\dP_{\kappa+1}$. To simplify the notation, we may assume without loss of generality that $\supp(q_0)=\supp(q_1)$. Denote this common set by $A$. 
    
    Our goal is to find the following:
    \begin{enumerate}
        \item An automorphism sequence $(\vec{\pi},\vec{D})=(\dot{\pi}_{\beta},\dot{D}_{\beta})_{\beta <\kappa+1}$;
        \item Pure extensions $w_0\leq^*q_0$, $w_1\leq^*q_1$ in $\mathbb{P}_{\kappa+1}$ such that $w_0,w_1\in\mathcal{D}(\vec{D})$ and $\Gamma_{\vec{\pi}}(w_0)=w_1$.
    \end{enumerate}

    Let $\alpha\in A$. If $o(\alpha)<\alpha^+$, let $\dot{t}_{\alpha}^i$ be such that $q_i(\alpha)=\langle\dot{t}_{\alpha}^i,\dot{T}_{\alpha}^i\rangle$. If $o^{\mathcal{U}}(\alpha)\geq\alpha^+$, let $\dot{t}_{\alpha}^i$ be such that $q_i(\alpha)=\langle\langle\dot{t}_{\alpha}^i,\dot{T}_{\alpha}^i\rangle,\dot{r}_{\alpha}^i\rangle$. By Lemma \ref{Lemma: Universal Representations}, fix functions $h_0,h_1\colon\kappa+1\to V_{\kappa+1}$ such that, for every $\alpha\in A$, $h_i(\alpha)=\dot{t}_{\alpha}^i$ and for every $(\beta,\sigma)\in\dom(\mathcal{U})$,
    $$j_{\beta,\sigma}(h_i)(\beta)=h_i(\beta)$$

   \smallskip

    We construct the sequence $(\dot{\pi}_{\beta},\dot{D}_{\beta})_{\beta <\kappa+1}$ by induction on $\beta$. For now, simply suppose that $(\dot{\pi}_{\gamma},\dot{D}_{\gamma})_{\gamma < \beta}$ has been constructed and is an automorphism sequence. If $\beta\notin\dom_1(\mathcal{U})$, there is nothing to do, since the $\beta$th iterand $\dot{\mathbb{Q}}_\beta$ is forced to be trivial. So, moving forward we assume that  $\beta\in \dom_1(\mathcal{U})$. To describe the inductive hypotheses, work in $V[G_{\beta}]$, where $G_{\beta}$ is $\dP_{\beta}$-generic. Define $H_{\beta}$ to be the upward closure of $\Gamma_{\vec{\pi}}[H_{\beta}]$, which is $\dP_{\beta}$-generic as well, by Corollary \ref{Corollary: Automorphism Sequence gives Generic}.

    A central notion is the following:

    \begin{definition}
        Let $t$ be a finite sequence of ordinals. Let $K_{\beta}\subseteq\dP_{\beta}$ be a generic filter. Let $\tau\leq o^{\mathcal{U}}(\beta)$. We say that $t$ is \emph{$\tau$-coherent with respect to $K_{\beta}$} if it is a $\tau$-coherent sequence in $V[K_{\beta}]$ according to Definition \ref{def: coherency}. We say that $t$ is \emph{coherent with respect to $K_{\beta}$} if it is $o^{\mathcal{U}}(\beta)$-coherent with respect to $K_{\beta}$.
    \end{definition}

    \smallskip
    
  For every $\tau\leq o^{\mathcal{U}}(\beta)$, we define an auxiliary function $\sigma_{\beta}^{\tau}$  as follows:

    \begin{definition}\label{def: sigmatau}
    
        $\sigma_{\beta}^{\tau}$ is a function on the  set of all minimal $\tau$-coherent sequences of ordinals below $\beta$. The values of $\sigma^\tau_\beta$ are defined as follows:

        Let us denote
        $$\text{$s_0^G:=(h_0(\beta))^{G_{\beta}}$, $s_0^H:=(h_0(\beta))^{H_{\beta}}$ and  $s_1^G:=(h_1(\beta))^{G_{\beta}}$, $s_1^H:=(h_1(\beta))^{H_{\beta}}$}$$

        \begin{enumerate}
            \item If there exists $r$ such that $(s_0^G\restriction\tau)^{\frown}r$ is $\tau$-coherent with respect to $G_{\beta}$, $(s_1^H\restriction\tau)^{\frown}r$ is $\tau$-coherent with respect to $H_{\beta}$ and $t=\reg((s_0^G\restriction\tau)^{\frown}r)$, then $$\sigma_{\beta}(t)=\reg((s_1^H\restriction\tau)^{\frown}r).$$
            \item If there exists $r$ such that $(s_1^G\restriction\tau)^{\frown}r$ is $\tau$-coherent with respect to $G_{\beta}$, $(s_0^H\restriction\tau)^{\frown}r$ is $\tau$-coherent with respect to $H_{\beta}$ and $t=\reg((s_1^G\restriction\tau)^{\frown}r)$, then $$\sigma_{\beta}(t)=\reg((s_0^H\restriction\tau)^{\frown}r).$$
            \item Otherwise, $\sigma_{\beta}(t)=t$.
        \end{enumerate}
    \end{definition}

    \begin{remark}
        The idea behind the map $\sigma_{\beta}^{\tau}$ is as follows: At stage $\alpha$, $\sigma_{\alpha}^{\tau}$ is going to swap $h_0(\alpha)\restriction\tau$ and $h_1(\alpha)\restriction\tau$, i.e. $\dot{t}_{\alpha}^0\restriction\tau$ and $\dot{t}_{\alpha}^1\restriction\tau$. So that this is possible, we must do the same procedure at many places below $\alpha$ as well.

        Ordinarily, we would want to simply map $s_0^{\frown}r$ to $s_1^{\frown}r$ and vice versa, ignoring sequences which do not extend $s_0$ or $s_1$. However, due to the intricate definition of the forcing, there can be sequences which do not extend $s_0$ or $s_1$ per se, but are equivalent to sequences extending $s_0$ or $s_1$. So that we preserve the ordering, we must also swap sequences like this and therefore deal with equivalence of sequences.
    \end{remark}

    \begin{remark}
        Note that the definition of the map $\reg$ is independent of the choice of generic (assuming the sequence is coherent with respect to the correct one), since it simply gets rid of all ordinals which are superceded by an ordinal of larger order, which is decided in the ground model.
    \end{remark}

    We formulate our \textbf{induction hypothesis}\label{induction hypothesis} in $V$. We assume the following:
    \begin{enumerate}
        \item[$(\aleph)$] $(\vec\pi,\vec{D})=(\dot{\pi}_\gamma,\dot{D}_\gamma)_{\gamma <\beta}$ is an automorphism sequence.
        \item[$(\beth)$] For every $\gamma < \beta$ with $o(\gamma)<\gamma^+$, for every $\dP_{\gamma}$-generic filter $G_{\gamma}$, setting $\pi_{\gamma}:=(\dot{\pi}_{\gamma})^{G_{\gamma}}$ and $D_\gamma:=(\dot{D}_\gamma)^{G_\gamma}$, it is the case that 
        $$\pi_\gamma\colon D_\gamma\to (\dot{D}_{\gamma})^{H_{\gamma}}\footnote{Recall that $(\dot{D}_{\gamma})^{H_{\gamma}}=(\Gamma_{\vec{\pi}}^{\gamma,-1})^*(\dot{D}_{\gamma})^{G_{\gamma}}$}$$
        is an isomorphism that is defined as $$\pi_{\gamma}(\langle t, T\rangle):=\langle\sigma_{\gamma}^{o^{\mathcal{U}}(\gamma)}(t),T\rangle.$$

        For every $\gamma\in \dom_1(\mathcal{U})\cap\beta$ with $o^{\mathcal{U}}(\gamma)=\gamma^+$, for every $\dP_{\gamma}$-generic filter $G_{\gamma}$, setting $\pi_{\gamma}:=(\dot{\pi}_{\gamma})_{G_{\gamma}}$ and $D_\gamma:=(\dot{D}_\gamma)_{G_\gamma}$, the following holds: If $\gamma\notin A$, then
        $$D_{\gamma}=\dQ_{\gamma}\text{ and }\pi_{\gamma}=\id_{\dQ_{\gamma}}$$
        If $\gamma\in A$, then there is a $\leq_{\gamma}^*$-dense set $E_{\gamma}\subseteq\dQ_{\gamma}^{\gamma^+}$, an isomorphism $\tau_{\gamma}$ from $E_{\gamma}$ to $(\dot{E}_{\gamma})^{H_{\gamma}}$ and a self-inverse bijection $\iota_{\gamma}$ between $\gamma^+$ and $\gamma^+$ such that
        $$\Vdash_{\dQ_{\gamma}^{\gamma^+}}\iota_{\gamma}[\supp(\dot{r}_{\gamma}^i)]\cap\supp(\dot{r}_{\gamma}^i)=\emptyset\footnote{Here we regard $\dot{r}_{\gamma}^i$ as a $\dQ_{\gamma}^{\gamma^+}$-name inside $V[G_{\gamma}]$}$$
        and
        $$\tau_{\gamma}(\langle t,T\rangle)=\langle\sigma_{\gamma}^{\gamma^+}(t),T\rangle$$
        as well as
        $$\pi_{\gamma}(\langle t,T\rangle,\dot{r})=(\tau_{\gamma}(\langle t,T\rangle),\tau_{\gamma}^*(\dot{\sigma}_{\iota_{\gamma}}(\dot{r})))$$
        Here, $\dot{\sigma}_{\iota_{\gamma}}$ is a $\dQ_{\gamma}^{\gamma^+}$-name for the automorphism of $\dot{\dR}_{\gamma}$ that maps $\dot{r}$ to the unique condition with $\supp(\dot{\sigma}_{\iota_{\gamma}}(\dot{r}))=\check{\iota}_{\gamma}[\supp(\dot{r})]$ such that $\dot{\sigma}_{\iota_{\gamma}}(\dot{r})(\check{\iota}_{\gamma}(\xi)=\dot{r}(\xi)$.
    \end{enumerate}

    The main goal for the rest of the proof is showing that, on a $\leq_{\beta}^*$-dense subset of $\mathbb{Q}_\beta:= (\dot{\dQ}_{\beta})^{G_{\beta}}$, the sequence  $\langle \sigma_\beta^\tau\mid \tau\leq o^{\mathcal{U}}(\beta)\rangle$ is well-defined and that it leads to a  longer  automorphism sequence
    $(\dot{\pi}_\gamma, \dot{D}_\gamma)_{\gamma < \beta+1}$
    satisfying Clauses~$(\aleph)$ and $(\beth)$.

    To simplify the notation, we fix the following:
   
   \begin{notation}
      Let $t$ be a finite sequence of ordinals below $\beta$. If $t$ is $\tau$-coherent with respect to $G_{\beta}$, let $\mathcal{U}_{\beta,\tau}(t)$ be the lift of $\mathcal{U}_{\beta,\tau}$ according to $G_{\beta}$. If $t$ is $\tau$-coherent with respect to $H_{\beta}$, let $\mathcal{V}_{\beta,\tau}(t)$ be the lift of $\mathcal{U}_{\beta,\tau}$ according to $H_{\beta}$.  
   \end{notation}

    The crucial claim is the following:
   
   \begin{claim}\label{claim: equality of measures}
       Let $t$ be a finite sequence of ordinals below $\beta$ and $\tau\leq o^{\mathcal{U}}(\beta)$. Suppose that that $t$ $\tau$-coherent with respect to $G_{\beta}$ and $\sigma_{\beta}^{\tau}(t)$ is $\tau$-coherent with respect to $H_{\beta}$. Then $\mathcal{U}_{\beta,\tau}(t)=\mathcal{V}_{\beta,\tau}(\sigma_{\beta}^{\tau}(t))$.
   \end{claim}

    \begin{remark}
        Technically, $\sigma_{\beta}^{\tau}(t)$ is only defined when $t$ is minimal, i.e. order-decreasing. However, we may and will ignore this distinction by letting $\sigma_{\beta}^{\tau}(t):=\sigma_{\beta}^{\tau}(\reg(t))$ if $t$ is $\tau$-coherent with respect to $G_{\beta}$.
    \end{remark}
    
   \begin{proof}[Proof of Claim \ref{claim: equality of measures}]
       Suppose that $X\in\mathcal{U}_{\beta,\tau}(t)$. Let $\dot{X}$ be a $\dP_{\beta}$-name such that $\dot{X}_{G_{\beta}}=X$. By the definition of $\mathcal{U}_{\beta,\tau}(t)$, there exist $p\in G_{\beta}$ and a $\dP_{\beta}$-name $\dot{T}$ such that
       $$p^{\frown}\langle t,\dot{T}\rangle^{\frown}(j_{\beta,\tau}(p)\smallsetminus(\beta+1))\Vdash\check{\beta}\in j_{\beta,\tau}(\dot{X})$$
       Since $\mathcal{D}(\vec{D})$ is $\leq^*$-dense, we may find $q\leq_{\dP_{\beta}}p$ such that $q\in G_{\beta}\cap\mathcal{D}(\vec{D})$. Furthermore, we can find such a $q$ such that it decides the values of $(h_0(\beta))$, $(h_1(\beta))$ and their translates by $\Gamma_{\vec{\pi}}$. Therefore, $j_{\beta,\tau}(q)\in j_{\beta,\tau}(\mathcal{D}(\vec{D}))=\mathcal{D}(j_{\beta,\tau}(\vec{D}))$. Thus, for every $\gamma\in \supp(j_{\beta,\tau}(q))$, $j_{\beta,\tau}(q)\restriction\gamma\Vdash j_{\beta,\tau}(q)(\gamma)\in j_{\beta,\tau}(\vec{D})(\gamma)$. By the maximum principle, let $\langle t,\dot{T}'\rangle$ be a $\dP_{\beta}$-name for a pure extension of $\langle t,\dot{T}\rangle$ that lies in $j_{\beta,\tau}(\vec{D})(\beta)$. It follows that
       $$q^{\frown}\langle t,\dot{T}'\rangle^{\frown}(j_{\beta,\tau}(q)\smallsetminus\beta+1)\in j_{\beta,\tau}(\mathcal{D}(\vec{D}))$$
       Therefore, we may assume that already $p^{\frown}\langle t,\dot{T}'\rangle^{\frown} (j_{\beta,\tau}(p)\smallsetminus(\beta+1))\in j_{\beta,\tau}(\mathcal{D}(\vec{D}))$.
       By the inductive hypothesis, $\Gamma_{j_{\beta,\tau}(\vec{\pi})}$ is an almost automorphism modulo equivalence of the $\leq^*$-dense set $\mathcal{D}(j_{\beta,\tau}(\vec{D}))$ (see Lemma~\ref{lemma: GammaandD}). Therefore, we have
       $$\Gamma_{j_{\beta,\tau}(\vec{\pi})}\left(p^{\frown}\langle t,\dot{T}\rangle^{\frown}(j_{\beta,\tau}(p)\smallsetminus(\beta+1))\right)\Vdash\check{\beta}\in\Gamma_{j_{\beta,\tau}(\vec{\pi})}^*(j_{\beta,\sigma}(\dot{X})).$$
       The following is clear, by elementarity and the fact that $j_{\beta,\tau}$ has critical point $\beta$:

       \begin{enumerate}
           \item $\Gamma_{j_{\beta,\tau}(\vec{\pi})}=j_{\beta,\tau}\left(\Gamma_{\vec{\pi}}\right)$;
           \item $\Gamma_{j_{\beta,\tau}(\vec{\pi})}\restriction\dP_{\beta}=\Gamma_{\vec{\pi}}$.
       \end{enumerate}

        Now consider $j_{\beta,\tau}(\vec{\pi})(\beta)$. Work in $M_{\sigma,\tau}$ for now. By Clause~$(\beth)$ in the inductive hypothesis on  page~\pageref{induction hypothesis} and the elementarity of $j_{\beta,\tau}$, it is forced that
        $$j_{\beta,\tau}(\vec{\pi})(\beta)(\langle \check{t},\dot{T}\rangle)=\langle j_{\beta,\tau}(\vec{\sigma})_{\beta}^{o^{j_{\beta,\tau}(\mathcal{U})}(\beta)}(\check{t}),\dot{T}\rangle.$$
        Note that, in $M_{\beta,\tau}$, $o^{j_{\beta,\tau}(\mathcal{U})}(\beta)=\tau$. Here, $j_{\beta,\tau}(\vec{\sigma})_{\beta}^{\tau}$ is a $\dP_{\beta}$-name. This name is defined from the $\dP_{\beta}$-names $j_{\beta,\tau}(h_0)(\beta)$, $j_{\beta,\tau}(h_1)(\beta)$ and their $\Gamma_{\vec{\pi}}$-translates. By assumption, we have $j_{\beta,\tau}(h_0)(\beta)=h_0(\beta)$ and $j_{\beta,\tau}(h_1)(\beta)=h_1(\beta)$ and, by assumption, all of the relevant names are decided by $p$. Ergo, $j_{\beta,\tau}(\vec{\sigma})_{\beta}^{\tau}$ as computed in $M_{\beta,\tau}[G_{\beta}]$ coincides with $\sigma_{\beta}^{\tau}$ as computed in $V[G_{\beta}]$. Let $s:=\sigma_{\beta}^{\tau}(t)$, computed in $V[G_{\beta}]$. It follows that
        $$p\Vdash j_{\beta,\tau}(\vec{\pi})(\beta)(\langle t,\dot{T}\rangle)=\langle\check{s},\dot{T}\rangle$$

        Therefore, working again in $V[G_{\beta}]$, we have:
        \begin{align*}
            & \Gamma_{j_{\beta,\tau}(\vec{\pi})}\left(p^{\frown}\langle t,\dot{T}\rangle^{\frown}(j_{\beta,\tau}(p)\smallsetminus\beta+1)\right) \\
            =\; & \Gamma_{\vec{\pi}}(p)^{\frown}\Gamma_{\vec{\pi}}^*(j_{\beta,\sigma}(\vec{\pi})(\beta)(\langle t,\dot{T}\rangle))^{\frown}(j_{\beta,\tau}\left(\Gamma_{\vec{\pi}}(p)\right)\smallsetminus(\beta+1)) \\
            =\; & \Gamma_{\vec{\pi}}(p)^{\frown}\Gamma_{\vec{\pi}}^*(\langle\check{s},\dot{T}\rangle))^{\frown}(j_{\beta,\tau}\left(\Gamma_{\vec{\pi}}(p)\right)\smallsetminus(\beta+1)) \\
            =\; & \Gamma_{\vec{\pi}}(p)^{\frown}\langle\check{s},\Gamma_{\vec{\pi}}^*(\dot{T})\rangle)^{\frown}(j_{\beta,\tau}\left(\Gamma_{\vec{\pi}}(p)\right)\smallsetminus(\beta+1)).
        \end{align*}
        So $\Gamma_{\vec{\pi}}(p)\in H_{\beta}$ and $\Gamma_{\vec{\pi}}^*(\dot{T})$ witness that $\Gamma_{\vec{\pi}}^*(\dot{X})^{H_{\beta}}=X\in\mathcal{V}_{\beta,\sigma}(\sigma_{\beta}^{\tau}(t))$.
    \end{proof}

    We are now ready to define our $\leq_{\dQ_{\beta}}^*$- dense subset $D_\beta$.
    
     Let $\langle t,T\rangle\in \mathbb{Q}_\beta$ be arbitrary, and we define when $\langle t,T\rangle\in D_\beta$. For $\langle t,T\rangle\in D_{\beta}$, we require first of all that $t$ is a minimal $o^{\mathcal{U}}(\beta)$-coherent sequence.
     
     Recall that
     $$\text{$s_0^G:=(h_0(\beta))^{G_{\beta}}$, $s_0^H:=(h_0(\beta))^{H_{\beta}}$ and  $s_1^G:=(h_1(\beta))^{G_{\beta}}$, $s_1^H:=(h_1(\beta))^{H_{\beta}}$}$$ (see Definition~\ref{def: sigmatau}). Let
     $$\Xi:=\max(\max(s_0^G),\max(s_0^H),\max(s_1^G),\max(s_1^H))$$

    \begin{enumerate}
        \item \textbf{\textbf{Case A}:} Suppose that, for some finite sequence $r\in\beta^{<\omega}$, $(s_0^G)^{\frown}r$ is $o^{\mathcal{U}}(\beta)$-coherent with respect to $G_{\beta}$, $(s_1^H)^{\frown}r$ is $o^{\mathcal{U}}(\beta)$-coherent with respect to $H_{\beta}$ and $t=\reg((s_0^G)^{\frown}r)$. We let $\langle t,T\rangle\in D_\beta$ if and only if for every $s\in T$, $\reg((s_1^H)^{\frown}r)^{\frown}s$ is coherent with respect to $H_{\beta}$ and $\min(s)>\Xi$.
        \item \textbf{\textbf{Case B}:} Suppose that, for some finite sequence $r\in\beta^{<\omega}$, $(s_1^G)^{\frown}r$ is $o^{\mathcal{U}}(\beta)$-coherent with respect to $G_{\beta}$, $(s_0^H)^{\frown}r$ is $o^{\mathcal{U}}(\beta)$-coherent with respect to $H_{\beta}$ and $t=\reg((s_1^G)^{\frown}r)$. We let $\langle t,T\rangle\in D_\beta$ if and only if for every $s\in T$, $\reg((s_0^H)^{\frown}r)^{\frown}s$ is coherent with respect to $G_{\beta}$ and $\min(s)>\Xi$.
        \item \textbf{\textbf{Case C}:} Suppose otherwise. We let $\langle t,T\rangle\in D_\beta$ if and only if for every $s\in T$, $\min(s)>\Xi$ and $t^{\frown}s$ is coherent with respect to both $G_{\beta}$ and $H_{\beta}$.
    \end{enumerate}

    Our goal for now is to show:

    \begin{claim}\label{claim: Dbeta dense}
        $D_\beta$ is $\leq_{\dQ_\beta}^*$-open dense in $\dQ_\beta$.
    \end{claim}

    \begin{proof}
        Let $\langle t,T\rangle\in\dQ_{\beta}$. By replacing $t$ with $\reg(t)$, which results in a pure extension, we may suppose that $t$ is minimal. We define the set $X_t$ as follows:

        \begin{enumerate}
            \item Suppose that \textbf{\textbf{Case A}} holds, witnessed by $r$. We let $\xi\in X_t$ if and only if $\reg((s_0^G)^{\frown}r)^{\frown}\langle\xi\rangle$ is coherent with respect to $G_{\beta}$ and $\reg((s_1^H)^{\frown}r)^{\frown}\langle\xi\rangle$ is coherent with respect to $H_{\beta}$ and $\xi>\Xi$.

            \item Suppose that \textbf{\textbf{Case B}} holds, witnessed by $r$. We let $\xi\in X_t$ if and only if $\reg((s_1^G)^{\frown}r)^{\frown}\langle\xi\rangle$ is coherent with respect to $G_{\beta}$ and $\reg((s_0^H)^{\frown}r)^{\frown}\langle\xi\rangle$ is coherent with respect to $H_{\beta}$ and $\xi>\Xi$.

            \item Suppose that \textbf{\textbf{Case C}} holds. We let $\xi\in X_t$ if and only if $t^{\frown}\langle\xi\rangle$ is coherent with respect to both $G_{\beta}$ and $H_{\beta}$ and $\xi>\Xi$.
        \end{enumerate}

        \begin{subclaim}
            Let $\tau\leq o^{\mathcal{U}}(\beta)$. The set $X_t$ is a member of $\mathcal{U}_{\tau,\sigma}(t)$.
        \end{subclaim}

        \begin{proof}
            Since $\Xi<\beta$, the requirement $\xi>\Xi$ can easily be satisfied, so we ignore it.
        
            Suppose first that \textbf{Case A} holds, witnessed by $r$. By Lemma \ref{Lemma: Filter Properties} (3), the collection $Y_t$ of all ordinals $\xi$ such that $t^{\frown}\langle\xi\rangle=\reg((s_0^G)^{\frown}r)^{\frown}\langle\xi\rangle$ is coherent with respect to $G_{\beta}$ is a member of $\mathcal{U}_{\beta,\tau}(t)$. By the same argument, the collection $Z_t$ of all ordinals $\xi$ such that $\sigma_{\beta}^{\tau}(t)^{\frown}\langle\xi\rangle=\reg((s_1^H)^{\frown}r)^{\frown}\langle\xi\rangle$ is coherent with respect to $H_{\beta}$ is a member of $\mathcal{V}_{\beta,\tau}(\sigma_{\beta}^{\tau}(t))$, which is equal to $\mathcal{U}_{\beta,\tau}(t)$ (by Claim~\ref{claim: equality of measures}). Thus, $X_t=Y_t\cap Z_t\in\mathcal{U}_{\beta,\tau}(t)$.

        \smallskip

            \textbf{Case B} is entirely symmetric to \textbf{Case A}. 

            \smallskip

            Suppose that \textbf{\textbf{Case C}} holds. Again, the collection $Y_t$ of all ordinals $\xi$ such that $t^{\frown}\langle\xi\rangle$ is coherent with respect to $G_{\beta}$ is in $\mathcal{U}_{\beta,\tau}(t)$ and the collection $Z_t$ of all ordinals $\xi$ such that $t^{\frown}\langle\xi\rangle$ is coherent with respect to $H_{\beta}$ is in
            $$\mathcal{V}_{\beta,\tau}(t)=\mathcal{U}_{\beta,\tau}(\sigma_{\beta}^{\tau}(t))=\mathcal{U}_{\beta,\tau}(t).$$
            Here, the last equality follows since $\sigma^\tau_\beta(t)$ is equivalent to $t$. Therefore, $X_t:=Y_y\cap Z_t\in\mathcal{U}_{\beta,\tau}(t)$.
        \end{proof}

        Now we are ready to define the pure extension $\langle t,T^*\rangle$ lying in $D_\beta$.

        \begin{enumerate}
            \item Suppose that \textbf{Case A} holds, witnessed by $r$. We let $s\in T^*$ if and only if $s\in T$, $\reg((s_0^G)^{\frown}r)^{\frown}s$ is coherent with respect to $G_{\beta}$ and $\reg((s_1^H)^{\frown}r)^{\frown}s$ is coherent with respect to $H_{\beta}$ and $\min(s)>\Xi$.

             \item Suppose that \textbf{Case B} holds, witnessed by $r$. We let $s\in T^*$ if and only if $s\in T$, $\reg((s_1^G)^{\frown}r)^{\frown}s$ is coherent with respect to $G_{\beta}$ and $\reg((s_0^H)^{\frown}r)^{\frown}s$ is coherent with respect to $H_{\beta}$ and $\min(s)>\Xi$.

             \item Suppose that \textbf{Case C} holds. We let $s\in T^*$ if and only if $s\in T$, $t^{\frown}s$ is coherent with respect to both $G_{\beta}$ and $H_{\beta}$ and $\min(s)>\Xi$.
        \end{enumerate}

        It is clear that, if $\langle t,T^*\rangle\in\dQ_{\beta}$, it is an element of $D_\beta$ and a pure extension of $\langle t,T\rangle$. Therefore, we are done after showing:

        \begin{subclaim}\label{Subclaim: Pure extension well-defined}
            $\langle t,T^*\rangle\in\dQ_{\beta}$.
        \end{subclaim}

        \begin{proof}
            The only substantial claim is showing that for all $s\in T^*$, $\Suc_{T^*}(s)\in\mathcal{U}_{\beta,\tau}(t^{\frown}s)$ for every $\tau<o^{\mathcal{U}}(\beta)$. Crucially, the case we are in for $t$ is the same case we are in for $t^{\frown}s$: If we are in \textbf{Case A} or B, witnessed by $r$, $r^{\frown}s$ witnesses that we are in \textbf{Case A} or \textbf{Case B} respectively. If we are in \textbf{Case C}, we are still in \textbf{Case C} for $r^{\frown}s$, since $\min(r)>\Xi$.

            So suppose we are in \textbf{Case A}. Then $\Suc_{T^*}(s)$ consists of all $\xi\in\Suc_{T^*}(s)$ such that $\reg((s_0^G)^{\frown}r)^{\frown}r^{\frown}\langle\xi\rangle$ is coherent with respect to $G_{\beta}$ and the sequence $\reg((s_1^H)^{\frown}r)^{\frown}s^{\frown}\langle\xi\rangle$ is coherent with respect to $H_{\beta}$. Since $\reg(x_0)^{\frown}x_1^{\frown}\langle\zeta\rangle$ is coherent with respect to a filter if and only if this is true for $\reg(x_0^{\frown}x_1)^{\frown}\langle\zeta\rangle$, we must have $\Suc_{T^*}(s)=\Suc_T(s)\cap X_{t^{\frown}s}$, which is as required.

            The same argument applies equally in \textbf{Case B} and \textbf{Case C}.
            \renewcommand{\qedsymbol}{$\square_{\text{Subclaim~\ref{Subclaim: Pure extension well-defined}}}$}
        \end{proof}

        This finishes the proof of the claim.
        \renewcommand{\qedsymbol}{$\square_{\text{Claim~\ref{claim: Dbeta dense}}}$}
    \end{proof}

    Now define the putative isomorphism $\pi_{\beta}$ as in Clause~$(\beth)$; namely,
 $$\pi_{\beta}(\langle t,T\rangle)=\langle\sigma_{\beta}^{o^{\mathcal{U}}(\beta)}(t),T\rangle\;\;\text{for $\langle t,T\rangle\in D_\beta$.}$$

    \begin{claim}
        $\pi_{\beta}$ is a map from $D_{\beta}$ into $(\Gamma_{\vec{\beta}}^{-1})^*(D_{\beta})^{G_{\beta}}$.
    \end{claim}

    \begin{proof}
        First, notice that
        $(\Gamma_{\vec{\beta}}^{-1})^*(D_{\beta})^{G_{\beta}}=\dot{D}_{\beta}^{H_{\beta}}$, i.e. $(\Gamma_{\vec{\beta}}^{-1})^*(D_{\beta})^{G_{\beta}}$ is defined the same way as $D_{\beta}$, except replacing $G_{\beta}$ by $H_{\beta}$ and thus replacing $H_{\beta}$ by $G_{\beta}$ (because $\Gamma$ is self-inverse modulo equivalence).

        So suppose that $\langle t,T\rangle\in D_{\beta}$. We check the three cases.

        \begin{enumerate}
            \item Suppose that \textbf{Case A} holds for $t$, witnessed by $r$, i.e. $t=\reg((s_0^G)^{\frown}r)$ and $\sigma_{\beta}^{o(\beta)}(t)=\reg((s_1^H)^{\frown}r)$. Let $s:=\sigma_{\beta}^{o(\beta)}(t)$. It follows that, in $V[H_{\beta}]$, \textbf{Case B} holds for $s$, witnessed by $r$ as well, since the roles of $G_{\beta}$ and $H_{\beta}$ are switched in $V[H_{\beta}]$. Then, because \textbf{Case A} holds for $t$ in $V[G_{\beta}]$, whenever $u\in T$, $\reg((s_0^G)^{\frown}r)^{\frown}s$ is coherent with respect to $G_{\beta}$ and $\reg((s_1^H)^{\frown}r)^{\frown}s$ is coherent with respect to $H_{\beta}$. Therefore, $\langle s,T\rangle\in \dot{D}_{\beta}^{H_{\beta}}$ by the definition.

            \item Suppose that \textbf{Case B} holds for $t$, witnessed by $r$, i.e. $t=\reg((s_1^G)^{\frown}r)$ and $\sigma_{\beta}^{o(\beta)}(t)=\reg((s_0^H)^{\frown}r)$. Let $s:=\sigma_{\beta}^{o(\beta)}(t)$. It follows that, in $V[H_{\beta}]$, \textbf{Case A} holds for $s$, witnessed by $r$ as well, since the roles of $G_{\beta}$ and $H_{\beta}$ are switched in $V[H_{\beta}]$. Then, because \textbf{Case B} holds for $t$ in $V[G_{\beta}]$, whenever $u\in T$, $\reg((s_1^G)^{\frown}r)^{\frown}s$ is coherent with respect to $G_{\beta}$ and $\reg((s_0^H)^{\frown}r)^{\frown}s$ is coherent with respect to $H_{\beta}$. Therefore, $\langle s,T\rangle\in \dot{D}_{\beta}^{H_{\beta}}$ by the definition.

            \item Suppose that \textbf{Case C} holds for $t$. It follows that, in $V[H_{\beta}]$, \textbf{Case C} holds for $t$ as well, since a witness e.g. for \textbf{Case A} in $V[H_{\beta}]$ would yield a witness for \textbf{Case B} in $V[G_{\beta}]$ and vice versa. Thus, whenever $s\in T$, $t^{\frown}s$ is coherent with respect to both $G_{\beta}$ and $H_{\beta}$ and $\min(s)>\Xi$, showing that $\langle t,T\rangle\in\dot{D}_{\beta}^{G_{\beta}}$.\qedhere \end{enumerate}
    \end{proof}

    \begin{claim}\label{Claim: Pi Beta Self Inverse}
        $\pi_{\beta}\circ\Gamma^*(\dot{\pi}_{\beta})^{G_{\beta}}=\id_{\Gamma^*(\dot{D}_{\beta})^{G_{\beta}}}$.
    \end{claim}

    \begin{proof}
        Since $G_{\beta}$ was arbitrary, using the previous observation, it is sufficient (and equivalent) to show that $\dot{\pi}_{\beta}^{H_{\beta}}\circ\pi_{\beta}=\id_{D_{\beta}}$. So let $\langle t,T\rangle\in D_{\beta}$.

        \begin{enumerate}
            \item Suppose that \textbf{Case A} holds for $t$, witnessed by $r$, i.e. $t=\reg((s_0^G)^{\frown}r)$ and $\sigma_{\beta}^{o(\beta)}(t)=\reg((s_1^H)^{\frown}r)$. Let $s:=\sigma_{\beta}^{o(\beta)}(t)$. It follows that, in $V[H_{\beta}]$, \textbf{Case B} holds for $s$, witnessed by $r$ as well, since the roles of $G_{\beta}$ and $H_{\beta}$ are switched in $V[H_{\beta}]$. Since $s=\reg((s_1^H)^{\frown}r)$, $(\dot{\sigma}_{\beta}^{o(\beta)})^{H_{\beta}}(s)=\reg((s_0^G)^{\frown}r)=t$. Thus, $\dot{\pi}_{\beta}^{H_{\beta}}(\pi_{\beta}(\langle t,T\rangle))=\langle t,T\rangle$ by the definition.

            \item Suppose that \textbf{Case B} holds for $t$, witnessed by $r$, i.e. $t=\reg((s_1^G)^{\frown}r)$ and $\sigma_{\beta}^{o(\beta)}(t)=\reg((s_0^H)^{\frown}r)$. Let $s:=\sigma_{\beta}^{o(\beta)}(t)$. It follows that, in $V[H_{\beta}]$, \textbf{Case A} holds for $s$, witnessed by $r$ as well, since the roles of $G_{\beta}$ and $H_{\beta}$ are switched in $V[H_{\beta}]$. Since $s=\reg((s_0^H)^{\frown}r)$, $(\dot{\sigma}_{\beta}^{o(\beta)})^{H_{\beta}}(s)=\reg((s_1^G)^{\frown}r)=t$. Thus, $\dot{\pi}_{\beta}^{H_{\beta}}(\pi_{\beta}(\langle t,T\rangle))=\langle t,T\rangle$ by the definition.

            \item Suppose that \textbf{Case C} holds for $t$. It follows that, in $V[H_{\beta}]$, \textbf{Case C} holds for $t$ as well. Therefore,
            $$\dot{\pi}_{\beta}^{H_{\beta}}(\pi_{\beta}(\langle t,T\rangle))=\dot{\pi}_{\beta}^{H_{\beta}}(\langle t,T\rangle)=\langle t,T\rangle$$
            as desired.\qedhere   \end{enumerate}
    \end{proof}

    \smallskip
    
     It is clear from the definition that $\pi_{\beta}$ preserves the pure extension ordering. Let us show that $\pi_{\beta}$ preserves the normal ordering:

    \begin{claim}\label{claim: pi respects one point extensions}
        Let $\langle s,S\rangle,\langle t,T\rangle\in D_{\beta}$ and suppose that $\langle s,S\rangle\leq_{D_{\beta}}\langle t,T\rangle$. Then $\pi_{\beta}(\langle s,S\rangle)\leq_{\dot{D}_{\beta}^{H_{\beta}}}\pi_{\beta}(\langle t,T\rangle)$
    \end{claim}

    \begin{proof}
        By the definition of $\leq_{\dQ_{\beta}}$ there exists $\eta\in T$ such that $\langle s,S\rangle\leq_{\dQ_{\beta}^*}\langle t^{\frown}\eta,T_{\eta}\rangle$. By the definition of $D_{\beta}$, $s$ is minimal and so we must have $s=\reg(t^{\frown}\eta)$.

        So suppose that \textbf{Case A} holds for $t$, witnessed by $r$. By the definition of $D_{\beta}$ it follows that \textbf{Case A} also holds for $s=\reg(t^{\frown}\eta)$, witnessed by $r^{\frown}\eta$. Ergo we have
        $$\sigma_{\beta}^{o(\beta)}(s)=\reg((s_1^H)^{\frown}r^{\frown}\eta)=\reg(\reg((s_1^H)^{\frown}r)^{\frown}\eta)=\reg(\sigma_{\beta}^{o^{\mathcal{U}}(\beta)}(t)^{\frown}\eta)$$
        which shows that
        $$\pi_{\beta}(\langle s,S\rangle)=\langle\sigma_{\beta}^{o(\beta)}(s),S\rangle\leq_{\dot{D}_{\beta}^{H_{\beta}}}\langle\sigma_{\beta}^{o(\beta)}(t),T\rangle=\pi_{\beta}(\langle t,T\rangle)$$
        as desired.

        The same argument works for \textbf{Case B} and \textbf{Case C} holding for $t$ as well. In the latter case we use that, because $\min(s)>\xi$, if \textbf{Case C} holds for $t$, \textbf{Case C} also holds for $t^{\frown}\eta$.
    \end{proof}
    
    Now suppose that $o^{\mathcal{U}}(\beta)<\beta^+$. We are in a position to conclude the proof of the induction:
    
    \begin{claim}
        For $o^{\mathcal{U}}(\beta)<\beta^+$, $(\dot{\pi}_\gamma, \dot{D}_\gamma)_{\gamma < \beta+1}$ witnesses Clauses~$(\aleph)$--$(\beth)$.
    \end{claim}
    
    \begin{proof}[Proof of claim]
        Clause~$(\beth)$ is clear, by construction. Let us check Clause~$(\aleph)$; namely,  $(\dot{\pi}_\gamma, \dot{D}_\gamma)_{\gamma < \beta+1}$ is an automorphism sequence: First, observe that the previous construction works irrespective of the chosen generic filter $G_\beta$, ergo the trivial condition of $\mathbb{P}_\beta$ forces all the claims established so far. In particular, invoking Claim~\ref{claim: Dbeta dense}  $\dot{D}_\beta$ is forced to be an open $\leq^*_{{\mathbb Q}_\beta}$-dense subset of ${\mathbb Q}_\beta$. Second, Claim~\ref{claim: pi respects one point extensions} and the comments preceding it  showed that $\dot{\pi}_\beta$ is forced to be an isomorphism from $\dot{D}_\beta$ onto $\Gamma^*_{\vec\pi \restriction\beta}(\dot{D}_\beta)$ such that  $\dot\pi_\beta\circ \Gamma_{\vec{\pi}\restriction\beta}^*(\dot\pi_\beta) = \id_{\dot{D}_\beta}$.   
    \end{proof}

    Now suppose that $o^{\mathcal{U}}(\beta)=\beta^+$. As before, we have that $\pi_{\beta}$ is an isomorphism of $D_{\beta}$ and $(\Gamma_{\vec{\pi}}^{-1})^*(\dot{D}_{\beta})^{G_{\beta}}$. Let us call the dense set $E_{\beta}$ and the automorphism $\tau_{\beta}$ instead so that we can use $D_{\beta}$ and $\pi_{\beta}$ again. Suppose first that $\beta\notin A$. We let $\pi_{\beta}$ be the automorphism of $\dR_{\beta}=\dQ_{\beta}^{\beta^+}*\dot{\dC}_{\beta}$ defined by $\pi_{\beta}(\langle t,T\rangle,\dot{r})=(\langle t,T\rangle,\dot{r})$. Now suppose that $\beta\in A$. For now, work in $V$. Recall that $q_i(\beta)=(\langle\dot{t}_{\beta}^i,\dot{T}_{\beta}^i\rangle,\dot{r}_{\beta}^i)$.

    Recall that $\dP_{\beta}$ has the $\beta^+$-cover property by Corollary \ref{corollary: cover property} and forces that $\dot{\dQ}_{\beta}^{\beta^+}$ is $\beta^+$-cc by Lemma~\ref{lemma: basic properties of Gitik's}. In $V$, fix a well-ordering $\prec_0$ on $[\beta^+]^{\beta}$ and fix a well-ordering $\prec_1$ on ${}^{\beta^+}\beta^+$. In $V[G_{\beta}]$, let $X_{\beta}\in V$ be the $\prec_0$-least set such that for some $p\in G_{\beta}$,
    $$(p,\one_{\dot{\dQ}_{\beta}^{\beta^+}})\Vdash_{\dP_{\beta}*\dot{\dQ}_{\beta}^{\beta^+}}\supp(\dot{r}^0_{\beta})\cup\supp(\dot{r}^1_{\beta})\subseteq\check{X}_{\beta}$$
    Also let $Y_{\beta}$ be the $\prec_0$-least set such that for some $q\in H_{\beta}$,
    $$(q,\one_{\dot{\dQ}_{\beta}^{\beta^+}})\Vdash_{\dP_{\beta}*\dot{\dQ}_{\beta}^{\beta^+}}\supp(\dot{r}^0_{\beta})\cup\supp(\dot{r}^1_{\beta})\subseteq\check{Y}_{\beta}$$
    
    In $V[G_{\beta}]$, let $\iota_{\beta}$ be the $\prec_1$-least self-inverse bijection from $\beta^+$ to $\beta^+$ such that $\iota_{\beta}[X_{\beta}\cup Y_{\beta}]$ is disjoint from $X_{\beta}\cup Y_{\beta}$. Clearly, $\iota_{\beta}^2=\id_{X_{\beta}}$. In $V[G_{\beta}][K_{\beta}]$, where $K_{\beta}$ is $\dQ_{\beta}^{\beta^+}$-generic over $V[G_{\beta}]$, let $\sigma_{\iota_{\beta}}$ be the automorphism of $\dot{\dR}_{\beta}^{K_{\beta}}$ such that, for every $r\in\dot{\dR}_{\beta}^{K_{\beta}}$, $\supp(\sigma_{\iota_{\beta}}(r))=\iota_{\beta}[\supp(r)]$ and $\sigma_{\iota_{\beta}}(r)(\iota_{\beta}(\gamma))=r(\gamma)$. Let $\dot{\sigma}_{\iota_{\beta}}$ be a name for $\sigma_{\iota_{\beta}}$. Let $D_{\beta}$ consist of all $(q,r)\in\dQ_{\beta}^{\beta^+}*\dot{\dR}_{\beta}$ such that $q\in E_{\beta}$. Since $E_{\beta}$ is $\leq_{\dQ_{\beta}^{\beta^+}}^*$-dense by Claim \ref{claim: Dbeta dense}, it follows that $D_{\beta}$ is $\leq_{\beta}^*$-dense. For $(q,\dot{r})\in D_{\beta}$, let $\pi_{\beta}(q,\dot{r}):=(\tau_{\beta}(q),\tau_{\beta}^*(\dot{\sigma}_{\iota_{\beta}}(\dot{r}))$.

    \begin{claim}\label{Claim: Verifying aleph and beth for order betaplus}
        For $o^{\mathcal{U}}(\beta)\geq\beta^+$, $(\dot{\pi}_{\gamma},\dot{D}_{\gamma})_{\gamma<\beta+1}$ witnesses Clauses~$(\aleph)$--$(\beth)$.
    \end{claim}

    \begin{proof}
        As in the previous proof, since $G_{\beta}$ was arbitrary, it suffices to show that the desired statements hold in $V[G_{\beta}]$.
    
        \begin{subclaim}\label{Subclaim: Pi Beta Isomorphism}
            $\pi_{\beta}$ is an isomorphism of $D_{\beta}$ and $(\Gamma_{\vec{\pi}}^{-1})^*(\dot{D}_{\beta})^{G_{\beta}}$.
        \end{subclaim}
    
        \begin{proof}
            By the definition, $(\Gamma_{\vec{\pi}}^{-1})^*(\dot{D}_{\beta})^{G_{\beta}}=\dot{D}_{\beta}^{H_{\beta}}$ consists of all $(q,\dot{c})$ such that $q\in\dot{E}_{\beta}^{H_{\beta}}$ and $q\Vdash\dot{c}\in(\dot{\dC}_{\beta}^{H_{\beta}})$ (here we regard $\dot{\dC}_{\beta}$ as a $\dP_{\beta}$-name for a $\dot{\dQ}_{\beta}^{\beta^+}$-name), i.e. $q\in\dot{E}_{\beta}^{H_{\beta}}$ and $q$ forces that $\dot{c}$ is a ${<}\,\beta^+$-supported function on $\beta^+$ that maps into the collection of closed bounded subsets of $\beta$ consisting of $V$-measurable cardinals. As we have seen before, if $(q,\dot{c})\in D_{\beta}$, then $\tau_{\beta}(q)\in\dot{E}_{\beta}^{H_{\beta}}$. Since $(q,\dot{c})\in\dQ_{\beta}^{\beta^+}$, $q$ forces that $\dot{c}$ and thus $\dot{\sigma}_{\iota_{\beta}}(\dot{c})$ is a ${<}\,\beta^+$-supported function from $\beta^+$ into the collection of closed bounded subsets of $\beta$ consisting of $V$-measurable cardinals. Therefore, $\tau_{\beta}(q)$ forces the same of $\tau_{\beta}^*(\dot{\sigma}_{\iota_{\beta}}(\dot{r}))$, showing that $(\tau_{\beta}(q),\tau_{\beta}^*(\dot{\sigma}_{\iota_{\beta}}(\dot{c})))\in \dot{E}_{\beta}^{H_{\beta}}$ as desired.

            A similar argument shows that $\pi_{\beta}$ preserves the pure extension ordering on $D_{\beta}$  and $(\Gamma_{\vec{\pi}}^{-1})^*(\dot{D}_{\beta})^{G_{\beta}}$. The more substantial claim is showing that the regular ordering is preserved. Recall that we defined $(q,\dot{c})\leq_{\dR_{\beta}}(q,\dot{c}')$ if $\pi(q,\dot{r})\leq_{\ro(\dC_{\beta})}\pi(q',\dot{r}')$. By Lemma \ref{Lemma: Easier description of Projection}, we have
            $$\pi((q,\dot{c}))=\bigwedge\{b\in\ro(\dC_{\beta})\;|\;q\Vdash\dot{c}\leq_{\ro(\dC_{\beta})}\check{b}\}$$
    
            To this end, let $\gamma_{\iota_{\beta}}$ be the map on $\dC_{\beta}$ defined as before: For $r\in\dC_{\beta}$, $\gamma_{\iota_{\beta}}(r)$ is a function with support $\iota_{\beta}[\supp(r)]$ such that $\gamma_{\iota_{\beta}}(r)(\iota_{\beta}(\gamma))=r(\gamma)$. We lift $\gamma_{\iota_{\beta}}$ to a function on $\ro(\dC_{\beta})$ (which we denote by the same symbol for simplicity) as follows: For $b\in\ro(\dC_{\beta})$, let $b=\bigwedge A$, where $A\subseteq\dC_{\beta}$ is an antichain. We let $\gamma_{\iota_{\beta}}(b):=\bigwedge\gamma_{\iota_{\beta}}[A]$. It is easy to see that $\gamma_{\iota_{\beta}}$ is an automorphism of $\ro(\dC_{\beta})$.

            \begin{subsubclaim}\label{Subsubclaim}
                For every $(q,\dot{c})\in D_{\beta}$, $$\dot{\pi}^{H_{\beta}}(\pi_{\beta}(q,\dot{c}))=\gamma_{\iota_{\beta}}(\pi(q,\dot{c})).$$
            \end{subsubclaim}

            \begin{proof}
                Note that $\pi_{\beta}(q,\dot{c})=(\tau_{\beta}(q),\tau_{\beta}^*(\dot{\sigma}_{\beta}(\dot{c})))$. Let $b\in\ro(\dC_{\beta})$ and suppose that
                $$\tau_{\beta}(q)\Vdash\tau_{\beta}^*(\dot{\sigma}_{\iota_{\beta}}(\dot{c}))\leq_{\ro(\dot{\dC}_{\beta})}\check{b}$$
                In particular, we have
                $$q\Vdash\dot{\sigma}_{\iota_\beta}(\dot{c})\leq_{\ro(\dot{\dC}_{\beta})}\check{b}$$
                since the definition of the ordering on $\ro(\dot{\dC}_{\beta})$ is absolute (it is defined simply by reverse inclusion on the order and pointwise end-extension). We claim that this implies that $q\Vdash\dot{c}\leq_{\ro(\dot{\dC}_{\beta})}\check{\gamma}_{\iota_{\beta}}(\check{b})$. As in the proof of Lemma \ref{Lemma: Easier description of Projection}, it suffices to show that $(q,\dot{c})$ forces that $\check{\gamma}_{\iota_{\beta}}(\check{b})$ belongs to the induced generic filter on $\ro(\dot{\dC}_{\beta})$.
                
                To this end, let $K*L$ be $\dQ_{\beta}^{\beta^+}*\dot{\dC}_{\beta}$-generic containing $q$. It follows that $K*\sigma_{\iota_{\beta}}[L]$ is $\dQ_{\beta}^{\beta^+}*\dot{\dC}_{\beta}$-generic containing $(q,\dot{\sigma}_{\iota_{\beta}}(\dot{c}))$. In particular in the model  $V[K*\sigma_{\iota_{\beta}}[L]]$, $b$ is in $\dot{G}_{\ro(\dC_{\beta})}^{K*\sigma_{\iota_{\beta}}[L]}$. (Recall that $\dot{G}_{\ro(\dC_{\beta})}^{K*\sigma_{\iota_{\beta}}[L]}$ is simply the upward closure of $\sigma_{\iota_{\beta}}[L]$, since $\sigma_{\iota_{\beta}}[L]\subseteq\dC_{\beta}$.) Thus, there must be $x\in L$ such that $\sigma_{\iota_{\beta}}(x)\leq b$. Now, by the definition of $\sigma_{\iota_{\beta}}$ and $\gamma_{\iota_{\beta}}$, we have $\sigma_{\iota_{\beta}}(x)=\gamma_{\iota_{\beta}}(x)$ for all $x\in\dot{\dC}_{\beta}^K$. Since $\gamma_{\iota_{\beta}}$ is lifted to $\ro(\dC_{\beta})$ and it is a self-inverse automorphism, we have $x\leq\gamma_{\iota_{\beta}}(b)$ and so $\gamma_{\iota_{\beta}}(b)\in\dot{G}_{\ro(\dC_{\beta})}^{K*L}$. As $K*L$ was arbitrary, we are done.

                Now by symmetry, it follows that
                $$q\Vdash_{\dQ_{\beta}^{\beta^+}}\dot{\sigma}_{\iota_{\beta}}(\dot{c})\leq_{\ro(\dC_{\beta})}\check{b}\text{ if and only if }q\Vdash_{\dQ_{\beta}^{\beta^+}}\dot{c}\leq_{\ro(\dC_{\beta})}\check{\gamma}_{\iota_{\beta}}(\check{b})$$
                proving the desired statement.
                \renewcommand{\qedsymbol}{$\square_{\text{Subsubclaim~\ref{Subsubclaim}}}$}
            \end{proof}

        Now we can finish the proof: Let $(q',\dot{c}')\leq_{\dot{\dR}_{\beta}}(q,\dot{c})$, i.e. $\pi(q',\dot{c}')\leq_{\ro(\dC_{\beta})}\pi(q,\dot{c})$. Therefore, we have
        $$\dot{\pi}^{H_{\beta}}(\pi_{\beta}(q',\dot{c}'))=\gamma_{\iota_{\beta}}(\pi(q',\dot{c}'))\leq_{\ro(\dC_{\beta})}\gamma_{\iota_{\beta}}(\pi(q,\dot{c}))=\dot{\pi}^{H_{\beta}}(\pi_{\beta}((q,\dot{c})))$$
        showing that $\pi_{\beta}(q',\dot{c}')\leq\pi_{\beta}(q,\dot{c})$.
        \renewcommand{\qedsymbol}{$\square_{\text{Subclaim~\ref{Subclaim: Pi Beta Isomorphism}}}$}
        \end{proof}

        The above proof also shows that $\pi_{\beta}$ has the form stipulated in~$(\beth)$.

        Finally, we must show the following:

        \begin{subclaim}
            $\pi_{\beta}\circ\Gamma_{\vec{\pi}}^*(\dot{\pi}_{\beta})^{G_{\beta}}=\id_{\Gamma_{\vec{\pi}}^*(\dot{D}_{\beta})^{G_{\beta}}}$.
        \end{subclaim}

        \begin{proof}
            As in the proof of Claim \ref{Claim: Pi Beta Self Inverse}, it suffices to show instead that $\dot{\pi}_{\beta}^{H_{\beta}}\circ\pi_{\beta}=\id_{D_{\beta}}$. By Claim \ref{Claim: Pi Beta Self Inverse}, we have that $\dot{\tau}_{\beta}^{H_{\beta}}\circ\tau_{\beta}=\id_{E_{\beta}}$.

            Second, we show that $\iota_{\beta}=\dot{\iota}_{\beta}^{H_{\beta}}$: Recalling the definition of $\iota_{\beta}$, we chose the $\prec_0$-least $X_{\beta}$ and $Y_{\beta}$ that were forced by conditions inside of $G_{\beta}$ and $H_{\beta}$ respectively to cover $\supp(\dot{r}_{\beta}^i)$ for $i<2$. It follows that, in $V[H_{\beta}]$, since the roles of $G_{\beta}$ and $H_{\beta}$ are reversed, $\dot{X}_{\beta}^{H_{\beta}}=Y_{\beta}$ and $\dot{Y}_{\beta}^{H_{\beta}}=X_{\beta}$. In particular, their union is the same. Thus, since $\iota_{\beta}$ was chosen in relation to $X_{\beta}\cup Y_{\beta}$ and the ground-model well-ordering, we have $\dot{\iota}_{\beta}^{H_{\beta}}=\iota_{\beta}$. Now it follows as in the proof of Lemma \ref{lemma: GammaandD} that $\dot{\pi}_{\beta}^{H_{\beta}}\circ\pi_{\beta}=\id_{D_{\beta}}$.
        \end{proof}

        This shows that $(\dot{\pi}_{\gamma},\dot{D}_{\gamma})_{\gamma<\beta+1}$ witnesses clauses~$(\aleph)$ and $(\beth)$ in the case $o^{\mathcal{U}}(\beta)\geq\beta^+$, which was our goal.
        \renewcommand{\qedsymbol}{$\square_{\text{Claim \ref{Claim: Verifying aleph and beth for order betaplus}}}$}
    \end{proof}
    
    Thereby, the induction has been accomplished and we obtain an automorphism sequence $(\dot\pi_\alpha,\dot D_\alpha)_{\alpha < \kappa+1}$ witnessing $(\aleph)$ and $(\beth)$.

    \smallskip

    We are now in a position to prove that the conclusion of almost homogeneity modulo equivalence holds for $\mathbb{P}_{\kappa+1}$ and $q_0,q_1\in\dP_{\kappa+1}$. Let $w_0\leq^*q_0$ be an element of $\mathcal{D}(\vec{D})$. Our goal is to show that $\Gamma_{\vec{\pi}}(w_0)$ is $\leq^*$-compatible with $q_1$. To this end, it suffices to show the following for every $\alpha\in A$:\footnote{Recall that $\supp(q_0)=\supp(q_1)$ and that $A$ was their common value.}
    
    \begin{enumerate}
        \item If $o(\alpha)<\alpha^+$, then $\Gamma_{\vec{\pi}}(w_0)(\alpha)$ is forced by $\Gamma_{\vec{\pi}\restriction\alpha}(w_0\restriction\alpha)$ to have the same stem as $q_1(\alpha)$.
        \item \label{keypoint for o(alpha) alpha+}If $o^{\mathcal{U}}(\alpha)=\alpha^+$, then if $\Gamma_{\vec{\pi}}(w_0)(\alpha)=(\langle\dot{s}_{\alpha}^0,\dot{S}_{\alpha}^0\rangle,\dot{u}_{\alpha}^0)$ and $q_1(\alpha)=(\langle\dot{t}_{\alpha}^1,\dot{T}_{\alpha}^1\rangle,\dot{r}_1)$, it is forced by $\Gamma_{\vec{\pi}\restriction\alpha}(w_0\restriction\alpha)$ that $\dot{s}_{\alpha}^0=\dot{t}_{\alpha}^1$ and that $\langle\dot{t}_{\alpha}^1,\dot{T}_{\alpha}^1\rangle$ forces that $\dot{u}_{\alpha}^0$ and $\dot{r}_{\alpha}^1$ are compatible.
    \end{enumerate}
    
    So let $H_{\alpha}$ be any $\dP_{\alpha}$-generic filter containing $\Gamma_{\vec{\pi}\restriction\alpha}(w_0\restriction\alpha)$ and let $G_{\alpha}:=\Gamma_{\vec{\pi}\restriction\alpha}^{-1}[H_{\alpha}]$. By the construction, 
    $$\Gamma_{\vec{\pi}}(w_0)(\alpha)^{H_{\alpha}}=\pi_{\alpha}(\Gamma_{\vec{\pi}\restriction\alpha}^*(w_0(\alpha))^{H_{\alpha}})=\pi_{\alpha}(w_0(\alpha)^{G_{\alpha}})$$

    Suppose that $o(\alpha)<\alpha^+$. By construction, the stems of $w_0(\alpha)^{G_{\alpha}}$ and $(q_0(\alpha))^{G_{\alpha}}=(h_0(\alpha))^{G_{\alpha}}$ are equal. Ergo, by the inductive hypothesis placed on $\pi_{\alpha}$, the stem of $\pi_{\alpha}(q_0(\alpha)^{G_{\alpha}})$ is equal to $(h_1(\alpha))^{H_{\alpha}}$. Note that $h_1(\alpha)$ was a name for the stem of $q_1(\alpha)$. Since $H_{\alpha}$ was arbitrary, it follows that $\Gamma_{\vec{\pi}\restriction\alpha}(w_0\restriction\alpha)$ forces that the stems of $\Gamma_{\vec{\pi}}(w_0)(\alpha)$ and $q_1(\alpha)$ are equal, just as desired.

    Suppose that $o^{\mathcal{U}}(\alpha)=\alpha^+$. Just like before, $(\dot{s}_{\alpha}^0)^{H_{\alpha}}=\dot{s}_1^{H_{\alpha}}$. Furthermore, by the inductive hypothesis, it is forced by $\dQ_{\alpha}^{\alpha^+}$ over $V[H_{\alpha}]$ that $\supp(\dot{u}_{\alpha}^0)=\iota_{\beta}[\supp(\dot{r}_{\alpha}^0)]$ and that this is disjoint from $\supp(\dot{r}_{\alpha}^1)$. Therefore, it is forced by $\dQ_{\alpha}^{\alpha^+}$ that $\dot{u}_{\alpha}^0$ and $\dot{r}_{\alpha}^1$ are compatible.

    This finishes the proof.
    \renewcommand{\qedsymbol}{$\square_{\text{Theorem~\ref{Theorem: Iteration is Homogeneous}}}$}
\end{proof}

\section{Preserving $\HOD$-supercompactness}\label{sec:PreservingHODsupercompacts}

In this section we conclude the proof of our main theorem. To this end, assume the hypotheses from Setup \ref{setup}, and let $\mathbb{P} = \mathbb{P}_{\kappa}*\dot{\mathbb{Q}}_\kappa$ be the forcing constructed in Section~\ref{sec: the iteration}. The main goal in this section is showing that after forcing with $\mathbb{P}$, $\kappa$ is {\rm HOD}-supercompact. 

\begin{lemma}
    $\mathbb{P} = \mathbb{P}_{\kappa}\ast \dot{\mathbb{Q}}_\kappa$ preserves cardinals.
\end{lemma}

\begin{proof}
This is immediate from Lemmas \ref{Lemma: Nonstatiterationspreservecardinals} and \ref{lemma: forcing equivalent}.    
\end{proof}

\begin{theorem}\label{Thm: HOD supercompactness}
    Assume that $\kappa$ is a supercompact cardinal and {\rm GCH} holds. Let $G*g\subseteq \mathbb{P}_{\kappa}*\dot{\mathbb{Q}}_{\kappa}$ be generic over $V$. Then $\kappa$ remains supercompact in $V[G*g]$. {Furthermore, if, in addition, $V = g{\rm HOD}$ holds in $V$, $\kappa$ is {\rm HOD}-supercompact in $V[G*g]$.}
\end{theorem}

\begin{proof}
    Let $\lambda> \kappa$ be such that for every  fine normal measure $\mathcal{V}$ on $\mathcal{P}_{\kappa}(\lambda)$, $$j_{\mathcal{V}}( \mathcal{U} )\restriction \kappa+1 = \mathcal{U}.$$ Let $\mathcal{V}$ be a fine normal measure on $\mathcal{P}_\kappa(\lambda)$ such that $j_{\mathcal{V}}(\ell)(\kappa) = \lambda$.
    \smallskip
    
    To simplify  forthcoming notations, let us denote $j := j_{\mathcal{V}}$ and $M := M_{\mathcal{V}}$.
   
    Since $j(\mathcal{U})\restriction \kappa+1 = \mathcal{U}$, we obtain that $$j({\mathbb{P}}_{\kappa}*{\dot{\mathbb{Q}}_{\kappa}})\restriction \kappa+1 = \mathbb{P}_{\kappa}*{\dot{\mathbb{Q}}_{\kappa}}.$$ We can thus naturally factor the iteration as
    $$j(\mathbb{P}_\kappa) = \mathbb{P}_\kappa \ast \dot{\mathbb{G}}_\kappa \ast j(\mathbb{P}_\kappa)/\dot{G}\ast \dot{g}.$$

    Let us exhibit a normal, fine, $\kappa$-complete ultrafilter on $\mathcal{P}^{V[G\ast g]}_\kappa(\lambda)$ in the  extension $V[G\ast g].$ We begin constructing a suitable master condition.

    \begin{claim}\label{Claim: master condition}
        There exists a $j(\mathbb{P}_{\kappa})\restriction \kappa+1$-name $\dot{m}\in M$ such that $(p,\dot{q}, \dot{m})\leq j(p,\dot{q})$ for all   $(p,\dot{q})\in G\ast g$.
    \end{claim}

    \begin{proof}
    The master condition $\dot{m}$ we will construct will have the form $(\dot{m}_0, \dot{m}_1)$, where it is forced that $ \dot{m}_0\in j(\mathbb{P}_\kappa)\setminus \kappa+1$, and $\one{}^\smallfrown \dot{m}_0$ forces that $\dot{m}_1\in j(\dot{\mathbb{Q}}_\kappa)$. 

\smallskip

   \textbf{Construction of $m_0$:} Recall that for every active forcing coordinate $\alpha\in \dom_1(\mathcal{U})\cup\{\kappa\}$, the poset $\mathbb{Q}_\alpha$ has a $\leq^*$-dense subset $\mathbb{D}_\alpha$ which is $<\alpha$-directed-closed with respect to the pure extension order  $\leq^*_{\mathbb{Q}_\alpha}$ (see Lemma~\ref{lemma: Some properties}). 
   
   Let
        $$  D= \{ p\in \mathbb{P}_{\kappa} \mid  \forall \alpha\in \supp(p), p\restriction \alpha\Vdash p(\alpha)\in \dot{\mathbb{D}}_\alpha \}. $$
        Clearly, $D\subseteq \mathbb{P}_{\kappa}$ is $\leq^*$-dense. 
        Over $M[G*g]$, denote
        $$j[G]\setminus (\kappa+1):= \{ j(p)\setminus (\kappa+1) \mid  p\in G \}.$$
        $$j[G\cap D]\setminus (\kappa+1):= \{ j(p)\setminus (\kappa+1) \mid  p\in G\cap D \}.$$
        By supercompactness, $j[G\cap D]\setminus (\kappa+1)\in M[G\ast g]$. Let us argue that this is a $\leq^*$-directed set of conditions. In fact, $j[G]\setminus (\kappa+1)$ is a $\leq^*$-directed set of conditions: Indeed, given $p_0,p_1\in G$, let $p^*\in G$ be a common extension. By the definition of extension, for all but finitely many values of $\alpha\in \supp(p_0)\cup \supp(p_1)$, 
        $$ p\restriction \alpha \Vdash p(\alpha)\leq^* p_0(\alpha), p_1(\alpha). $$
        In particular, $j(p)\setminus (\kappa+1)\leq^* j(p_0)\setminus (\kappa+1), j(p_1)\setminus (\kappa+1)$ as needed.

        \smallskip

      In $V[G\ast g]$, $| j[G\cap D]\setminus (\kappa+1) | = \lambda$. Since $j(\ell)(\kappa)= \lambda$, the forcing $j(\mathbb{P})\restriction (\kappa, \lambda]$ does not contain active stages. In particular, each poset in $j(\langle \mathbb{D}_{\alpha} \colon \alpha<\kappa \rangle)\setminus (\kappa+1)$ has a direct extension order which is $\lambda^+$-directed closed. Since $V[G\ast g]\vDash {}^\lambda M[G*g]\subseteq M[G*g]$, $\langle j(D)\setminus (\kappa+1), \leq^*\rangle$ is $\lambda^+$-directed closed in $V[G\ast g]$. It thus follows that there exists a condition $m_0\in j(\mathbb{P})\setminus (\kappa+1)$ such that, for every $p\in G\cap D$, $m_0\leq j(p)\setminus (\kappa+1)$. By density of $D$, we deduce that for every $p\in G$, $m_0\leq j(p)\setminus (\kappa+1)$.

      \smallskip

        \textbf{Construction of $\dot{m}_1$:} 
        Work in $V[G]$. The generic $g\s \dot{\mathbb{Q}}_\kappa^G$ induces a sequence  $\langle c^g_\alpha\colon \alpha<\kappa^+ \rangle$ of clubs on  $\kappa$. We construct $\dot{m}_1$ by defining a $j(\mathbb{P}_\kappa)$-name for a condition  $\dot{m}^*_1$ forced in $ \dot{\mathbb{Q}}_{j(\kappa)}$ which is projected onto $\dot{m}_1$. Denote  $\one$ the trivial condition of   $\mathbb{Q}^{j(\kappa)^+}_{j(\kappa)}$ and let $\dot{\delta}$ be the $\mathbb{Q}^{j(\kappa)^+}_{j(\kappa)}$-name for the least $M$-measurable cardinal above $\kappa$ in the generic club subset of $j(\kappa)$ that is induced by the generic for $\mathbb{Q}^{j(\kappa)^+}_{j(\kappa)}$. 
        Let $m^*_1 = (\one, \dot{x})$ where $\dot{x}$ is forced to be a function with domain $j[\kappa^+]\in ([ j(\kappa)^+ ]^{<j(\kappa)})^M$, and for every $\alpha<\kappa^+$, $$\dot{x}(j(\alpha)) = c^g_{\alpha}\cup \{ \check{\kappa} \}\cup \{\dot{\delta}\}.\footnote{Note that it was necessary to add $\delta$, as we typically do not have $\kappa$ in the generic club of  $\mathbb{Q}^{j(\kappa)^+}_{j(\kappa)}$, since $o^{\mathcal{U}}(\kappa) = \kappa^++1$. It was also necessary to add $\kappa$ since each $c^{g}_{\alpha}$ is unbounded in $\kappa$. Finally, both $\kappa, \delta$ are measurables in $M$, so adding them to $\dot{c}^g_{\alpha}$ results in valid conditions.} $$
        
        This concludes the definition of $\dot{m} = (\dot{m}_0, \dot{m}_1)$.
        
        \smallskip
        
      \textbf{Verification:}  It now remains to verify that for every $(p, \dot{q})\in G\ast g$, $$(p,\dot{q}, \dot{m}) \leq j(p,\dot{q}).$$ 
        
        Fix $(p, \dot{q})\in G\ast g$. Since $\mathbb{P}_{\kappa}$ has a nonstationary support, $\kappa\notin \supp(j(p))$, and thus, by the choice of $\dot{m}_0$, we have $(p,\dot{q},\dot{m}_0)\leq j(p)$. Thus, it remains to prove that 
        $$(p,\dot{q},\dot{m}_0)\Vdash \dot{m_1}\leq j(\dot{q}).$$
        It suffices to show that for every $(p,\dot{q})$ there are densely many conditions $(p^*, \dot{q}^*)\leq (p, \dot{q})$ such that $(p^*,\dot{q}^*, \dot{m}_0)\Vdash \dot{m}_1\leq j( \dot{q}^* )$.

        Fix $(p',\dot{q}')\leq (p, \dot{q})$. Let  $(p^*,\dot{q}^*)\leq (p', \dot{q}')$ and $A\in ([\kappa^+]^{<\kappa})^V$ be such that 
        $$ p^*\Vdash A = \dom(\dot{q}^*) .$$
        Such $A$ can be found because $\mathbb{P}_{\kappa}$ satisfies the $\kappa$-cover property, by Corollary \ref{corollary: cover property}.
        Since $|A|<\kappa$, we can assume, by further extending $p^*$, that for every $\alpha\in A$, there exists $\beta_\alpha<\kappa$ and a $\mathbb{P}_{\beta_\alpha}$-name $\dot{x}_\alpha\in V_{\kappa}$ such that 
        $$ p^*\restriction \beta_\alpha \Vdash p^*\setminus \beta_\alpha \Vdash \dot{q}^*(\alpha) = \dot{x}_\alpha.  $$

        Applying $j$, we obtain that for every $\alpha\in A$, 
        $$ (p^*, \dot{q}^*, \dot{m}_0)\Vdash j(\dot{q}^*)(j(\alpha)) = \dot{x}_\alpha.  $$
        Since this is true for every $\alpha\in A$, and since $(p^*,\dot{q}^*,\dot{m}_0)\Vdash j[A] = \dom(j(\dot{q}^*))$, the above choice of $\dot{m}_1$ ensures that 
        $$(p^*, \dot{q}^*, \dot{m}_0)\Vdash  \dot{m_1}\leq j(\dot{q}^*), $$
        as desired.
    \end{proof}

    Now, having constructed a master condition $\dot{m}$, consider the poset
    $$\mathbb{R}:=\langle j(\mathbb{P}_\kappa*\dot{\mathbb{Q}}_\kappa)\setminus (\kappa+1), \leq^*\rangle.$$

    \begin{claim}\label{Claim: Generic for direct extension order}
        There exists a set $H\in V[G\ast g]$ such that $H\subseteq \mathbb{R}$ is generic over $M[G\ast g]$ and $\dot{m}_{G\ast g}\in H$.
    \end{claim}
    
    \begin{proof}
    Work in $V[G*g]$. The poset $\mathbb{R}$ is $\lambda^+$-strategically closed, since it has a dense subset which is $\lambda^+$-directed closed (that is, $j(D\ast \mathbb{D}_\kappa)\setminus (\kappa+1)$, where $D$ is as in the proof of Claim \ref{Claim: master condition}). It is not hard to verify that there exists in $V[G\ast g]$ an enumeration $\langle E_\alpha\mid \alpha<\lambda^+ \rangle$ consisting of all dense subsets of $\mathbb{R}$ in $V[G\ast g]$. Now by a standard argument, construct $H$ by playing below $\dot{m}_{G\ast g}$ a game of length $\lambda^+$ in which one of the players chooses conditions in the dense sets $E_\alpha$, and the other makes sure that there are lower bounds at limit stages. The moves in the game are the conditions which generate the generic $H$ inside the extension $V[G\ast g]$.
    \end{proof}

    We are finally prepared to define a supercompact measure $\mathcal{V}^*\in V[G\ast g]$ which extends $\mathcal{V}$. For every $X\subseteq (\mathcal{P}_{\kappa}(\lambda))^{V[G\ast g]}$, let  $X\in \mathcal{V}^*$ if and only if there is a $\mathbb{P}_\kappa\ast \mathbb{Q}_\kappa$-name $\dot{X}$ for $X$, a condition $(p,\dot{q})\in G$ and a $\mathbb{P}_\kappa*\mathbb{Q}_\kappa$-name $\dot{r}$ with $(p,\dot{q})\Vdash \dot{r}\in \dot{H}$, such that
    $$ (p,\dot{q}, \dot{r})\Vdash \check{j[\lambda]}\in j(\dot{X}). $$

    \begin{claim}
        $\mathcal{V}^*$ is a fine, normal measure on $(\mathcal{P}_{\kappa}(\lambda))^{V[G\ast g]}$ which extends the measure  $\mathcal{V}$.
    \end{claim}

    \begin{proof}
        It's straightforward to verify that $\mathcal{V}^*$ is well-defined, upward closed and extends $\mathcal{V}$. Those properties imply that $\mathcal{V}^*$ is fine. The fact that $\mathcal{V}^*$ is closed under finite intersections follows from the fact that any pair of conditions in $H$ have a common extension in $H$. Thus, it remains to prove that $\mathcal{V}^*$ is a $\kappa$-complete normal ultrafilter.

        Assume that $\langle \dot{X}_i \mid i<\gamma \rangle$ for some $\gamma<\kappa$ is forced by some condition $(p,\dot{q})\in G\ast g$ to be a partition of $\dot{\mathcal{P}}_{\kappa}(\lambda)$. We thus obtain that, over $M$, $(\dot{p},\dot{q}, \dot{m})$ forces that $\langle j(\dot{X}_i)\mid i<\gamma \rangle$ is a partition of $j( \dot{\mathcal{P}}_{\kappa}(\lambda)  )$. Since $j(\mathbb{P}_\kappa \ast \mathbb{Q}_{\kappa})\setminus (\kappa+1)$ satisfies the Prikry property and $\mathbb{R}$ is $\lambda^+$-strategically closed in $M[G*g]$, the set of conditions in $\mathbb{R}$ which decide the value of $i$ for which $j[\lambda]\in j(\dot{X}_i)$ is $\leq^*$-dense below $\dot{m}$. Since $H$ is $\leq^*$-generic for $\mathbb{R}$, 
            $$(p,\dot{q})\Vdash \text{ there are }r\in \dot{H} \text{ and } i<\gamma \text{ such that } r\Vdash j[\lambda]\in j(\dot{X}_i).$$
        By extending $(p,\dot{q})$ inside $G\ast g$, we can assume that $(p,\dot{q})$ decide the value of $\dot{r},i$. Therefore $(p,\dot{q},\dot{r})\Vdash j[\lambda]\in j(\dot{X}_i)$, namely $\dot{X}_i\in \mathcal{V}^*$, as desired.  

        Finally, let us argue that $\mathcal{V}^*$ is normal. Let $\dot{f}$ be a $\mathbb{P}_\kappa*\mathbb{Q}_\kappa$-name for a regressive function on $\dot{\mathcal{P}}_{\kappa}(\lambda)$, forced by some condition $(p,\dot{q})\in G\ast g$. It thus follows that over $M$,
        $(p,\dot{q},\dot{m})$ forces that $j(\dot{f})(j[\lambda])\in j[\lambda]$. Since $j(\mathbb{P}_\kappa \ast \mathbb{Q}_{\kappa})\setminus (\kappa+1)$ satisfies the Prikry property, the set of conditions $r\in \mathbb R$ which decide the unique value $\alpha<\lambda$ such that $j(\dot{f})(j[\lambda])= j(\alpha)$ is a dense subset of $\mathbb{R}$ below $\dot{m}$ which belongs to $M[G\ast g]$, since $\mathbb{R}$ is $\lambda^+$-strategically closed in $M[G*g]$ and $j\restriction \lambda\in M$. Thus,
            $$ (p,\dot{q})\Vdash \text{ there are } r\in \dot{H} \text{ and }  \alpha<\lambda \text{ such that } r\Vdash j(\dot{f})(j[\lambda]) = j(\alpha). $$
        By extending $(p,\dot{q})$ inside $G\ast g$, we can assume that they decide the value of $\dot{r}, \alpha$. Thus $(p,\dot{q})$ witnesses that $\{ x\in (\mathcal{P}_{\kappa}(\lambda))^{V[G\ast g]} \mid f(x) = \alpha \}\in \mathcal{V}^*$.
    \end{proof}

This proves that $\kappa$ remains supercompact in $V[G\ast g]$.

\smallskip

Finally, let us show that $\kappa$ remains ${\rm HOD}$-supercompact in $V[G\ast g]$ provided $V=\mathrm{gHOD}$ holds:

\begin{claim}\label{lemma:HOD supercompactness}
   If $V = {\rm{gHOD}}$ holds then  $V[G\ast g]\models \kappa$ is {\rm{HOD}}-supercompact.
\end{claim}

\begin{proof}[Proof of claim]
    Fix $\lambda> \kappa$ a strong limit cardinal and $\mathcal{V}\in V$ a supercompactness measure on $\mathcal{P}_\kappa(\lambda)$ such that $j_{\mathcal{V}}(\ell)(\kappa) = \lambda$ and $j_{\mathcal{V}}(\mathcal{U})\restriction \kappa+1 = \mathcal{U}$. Let $\mathcal{V}^*\in V[G\ast g]$ be a supercompactness measure on $(\mathcal{P}_\kappa(\lambda))^{V[G\ast g]}$ extending $\mathcal{V}$, constructed as above. We argue that 
    $$j_{\mathcal{V}^*}( {\text{HOD}}^{V[G\ast g]} )\cap (V_{\lambda})^{V[G\ast g]} = (\text{HOD})^{V[G\ast g]}\cap (V_\lambda)^{V[G\ast g]}.$$  
    The forcing $\mathbb{P}_\kappa*\mathbb{Q}_{\kappa}$ is homogeneous modulo equivalence (see Theorem~\ref{Theorem: Iteration is Homogeneous}). Furthermore, it is ordinal-definable, because $V={\rm{gHOD}}$. Therefore, by Lemma~\ref{lemma:HODincluded}, it follows that  $\text{HOD}^{V[G\ast g]} = V$. So, the right hand side  is $V_\lambda$.
    
    Next, assuming that $j_{\mathcal{V}^*}\colon V[G\ast g]\to M^*$ is the ultrapower embedding associated with $\mathcal{V}^*$ over $V[G\ast g]$, let $M = j_{\mathcal{V}^*}(V)$ be an inner model of $M^*$ such that for $G^* = j_{\mathcal{V}^*}(G*g)$, $M^* = M[G^*]$. By elementarity, $M$ satisfies $V = {\rm{gHOD}}$. Again, because $j_{\mathcal{V}^*}(\mathbb{P}_\kappa\ast \dot{\mathbb{Q}}_\kappa)$ is homogeneous modulo equivalence, the left hand side is equal to $M\cap (V_\lambda)^{V[G\ast g]}$. Thus, it suffices to argue that $(V_\lambda)^M = V_\lambda$.

    \smallskip

    Define an embedding $k\colon M_{\mathcal{V}}\to M$ via $k( [f]_{\mathcal{V}} ) = [f]_{\mathcal{V}^*}$ for every function  $f\colon (\mathcal{P}_\kappa(\lambda))^V\to V$ in $V$. It's not hard to check that $k$ is elementary and  $$j_{\mathcal{V}^*} = k\circ j_{\mathcal{V}}.$$ 
    
    We argue that $\crit(k) > \lambda $. Let $\mu \leq \lambda$.  Then in $V$, 
    $$ \mu= [x\mapsto \otp(x\cap \mu)]_{\mathcal{V}} $$
    and thus 
    $$k(\mu) = [x\mapsto \otp(x\cap \mu)]_{\mathcal{V}^*} = \otp\left(
    j_{\mathcal{V}^*}[ \lambda ]\cap j_{\mathcal{V}^*}(\mu)\right) = \mu.$$
    
    Finally, since $\crit(k) > \lambda $, 
    $$V_\lambda = (V_\lambda)^{M_{\mathcal{V}}} = k\left((V_\lambda)^{M_{\mathcal{V}}}\right) = (V_\lambda)^{M}.$$ 
    This yields the desired equalities. \end{proof}

This concludes the proof of Theorem \ref{Thm: HOD supercompactness}.
\end{proof}

\begin{theorem}\label{thm:finally}
    Assume the existence of a supercompact cardinal. Then, there is a cardinal-preserving generic extension   where:
    \begin{enumerate}
        \item The $\HOD$ hypothesis holds.
        \item $\kappa$ is $\HOD$-supercompact.
        \item There is a club $C\s \kappa$ consisting of measurable cardinals in $\HOD$.
    \end{enumerate}
\end{theorem}
\begin{proof}
Without loss of generality, we assume that $V=\mathrm{gHOD}$ and the GCH holds (see Fact~\ref{fact:BrookeTaylor}). Working over this model let $\mathcal{U}$ be a coherent sequence of measures as in Lemma~\ref{Lemma: Construction of the coherent sequence} and force with the associated  iteration $\mathbb{P}_{\kappa+1}$. This is an ordinal definable iteration (because $\mathcal{U}$ is $\OD$), and it is almost homogeneous modulo equivalence (by Theorem~\ref{Theorem: Iteration is Homogeneous}). As a result, $\HOD^{V[G\ast g]}= V$ for every generic $G\ast g\s \mathbb{P}_{\kappa+1}$.  By Claim~\ref{lemma:HOD supercompactness}, $\kappa$ is $\HOD$-supercompact in $V[G\ast g]$. In addition, $\mathbb{Q}_\kappa$ shoots a club through the set $\{\alpha < \kappa\mid \text{$\alpha$ is $V$-measurable}\}$, ergo (3) above holds. Finally, $(1)$ above holds in that $\mathbb{P}_{\kappa+1}$ is a set-sized forcing and the $\HOD$ hypothesis holds in the ground model. 
\end{proof}
\begin{remark}
Due to the flexibility in the construction of $\mathcal{U}$ described in Remark~\ref{rem:flexibility} we can arrange a situation where $\kappa$ is $\HOD$-supercompact and, say, $C\s \kappa$ consists of supercompact cardinals in $\HOD$.
\end{remark}
A question that has eluded our methods is the following:
\begin{question}
Assume large cardinals. Is the conclusion of Theorem~\ref{thm:finally} consistent with $``\theta^{+\HOD}<\theta^+$ for all $\theta\in C$''?  
\end{question}
By results of \cite{GOP} the above configuration is known to be consistent  if $\kappa$ is just assumed  supercompact. The basic issue is the lack of tools to control lifts of normal, fine, $\theta$-complete ultrafilters on $\mathcal{P}_\theta(\theta^+)$ which in turn do not enable the homogeneity arguments given in \S\ref{sec: homogeneity}.

\section{On successors of singular cardinals in $\HOD$}\label{sec:SuccessorsOfSingulars}
In this section we want to address an issue raised by Woodin regarding successors of singular cardinals in $\HOD$. In his \emph{Suitable Extenders Model I} paper, on \cite[p.299]{WOEM1}, Woodin claims that if $\gamma$ is a singular strong limit cardinal and $\gamma^{+\HOD}<\gamma^+$ then either $\cf(\gamma^{+\HOD})=\omega$ or $\cf(\gamma^{+\HOD})=\cf(\gamma)$, and immediately after he conjectures that $\cf(\gamma^{+\HOD})=\omega$ should hold. The same issue appears referenced in Goldberg's \cite[Question~3.8]{GoldbergHODconjecture}. 

In the next theorem we use \emph{Supercompact Radin forcing}  to establish a negative answer to this conjecture. We assume familiarity with Supercompact Radin forcing as exposed in  \cite{CumFriGol}.
\begin{theorem}
    Suppose that $\kappa$ is a $\kappa^{+}$-supercompact cardinal and the axiom $\mathrm{V}=\mathrm{gHOD}$ holds. Then, there is a generic extension where $\kappa$ is a strong limit cardinal,  $\cf(\kappa^{+\HOD})=\omega_1$ and $\kappa^{+\HOD}\neq \kappa^+$. Moreover, in this model  $\cf^{\HOD}(\kappa)=\omega_1$.
\end{theorem}
\begin{proof}
    Let $j\colon V\rightarrow M$ be an elementary embedding such that  $\crit(j)=\kappa$, $j(\kappa)>\kappa^+$ and $M^{\kappa^+}\s M$. For each $\alpha<\omega_1$ define by induction a sequence $u:=\langle u(\alpha)\mid \alpha<\omega_1\rangle$ as follows: $u(0):=j[\kappa^+]$ and $u(\alpha):=\{X\s V_{\kappa}\mid u\restriction \alpha\in j(X)\}.$  It follows that $u\in \mathcal{U}^{\sup}_\infty$ as defined in \cite[Definition~3.6]{CumFriGol}.

    \smallskip

    Let $ \mathbb{R}^{\sup}_u$ be the corresponding supercompact Radin forcing (see \cite[\S3]{CumFriGol}) and $G\s \mathbb{R}^{\sup}_u$ generic for it.
     The analysis of \cite[\S4]{CumFriGol} shows that there is a Radin-like extension $V[G_0]\s V[G]$, preserving cardinals of $V$, such that $\HOD^{V[G]}\s V[G_0]$. Since $V=\mathrm{gHOD}$  holds in $V$, we have
     $$V\s \HOD^{V[G]}\s V[G_0].$$
  In $V[G]$, the following hold:
$$\kappa^{+\HOD^{V[G]}}=\kappa^{+V}<\kappa^{++V}=\kappa^{+V[G]}.$$ 
  In addition  $\cf^{V[G]}(\kappa^{+\HOD})=\cf^{V[G]}(\kappa^{+V})=\cf^{V[G]}(\kappa)=\omega_1$.

    We are left with showing that $\cf^{\HOD^{V[G]}}(\kappa)=\omega_1$. Let $C\s \kappa$ be the generic Radin club of order-type $\omega_1$ introduced by $G$. Note that all points in $C$ are  regular in $V$. We claim that $\acc(C)$ is  ordinal definable in $V[G]$. Indeed, we claim that 
$$\acc(C)=\{\rho<\kappa\mid V[G]\models ``\mathcal{G}_\kappa\models ``\cf(\rho)=\rho"\,\wedge\,\cf(\rho)<\rho"\}$$
   where  $\mathcal{G}$ is a predicate denoting the ground of $V[G]$  and $\mathcal{G}_\kappa$ the predicate for the members of $\mathcal{G}$ of rank ${<}\kappa$. 
Since $V=\mathrm{gHOD}$ holds in $V$ the satisfability relation $\models$ is ordinal definable,  the right-hand-side is in $\HOD$.
     
The left-to-right inclusion is obvious as $\otp(C)=\omega_1$. For the converse, suppose that $\rho$ is a $V$-regular cardinal that changes its cofinality. Let $\alpha<\omega_1$ be the first ordinal such that $\sigma_\alpha\leq \rho<\sigma_{\alpha+1}$. Let $p\in G$ be a condition mentioning measure sequences $u_\alpha, u_{\alpha+1}\in \mathcal{U}^{\sup}_\infty$ such that $\kappa_{u_\alpha}=\sigma_\alpha$ and  $\kappa_{u_{\alpha+1}}=\sigma_{\alpha+1}$. By the \emph{Factor Lemma} of Radin forcing, $\mathbb{R}_u^{\sup}/p$ is isomorphic to $\mathbb{R}^{\sup}_{u_\alpha}/p_0\times \mathbb{R}^{\sup}_{u}/p_1$ where $p_0$ is the truncation of $p$ up to the occurence of $u_\alpha$ and $p_1$ denotes the rest of $p$. The latter forcing does not change cofinalities below $\sigma_{\alpha+1}$ so it has to be $\mathbb{R}^{\sup}_{u_\alpha}/p_0$  the responsible for changing the cofinality of $\rho$. However, this poset is $\sigma_\alpha^+$-cc so it must be that $\rho=\sigma_\alpha$. If $\sigma_\alpha$ were to be a successor point of $C$ then it is easy to show (again, using the Factor Lemma) that $\sigma_\alpha$ does not change its cofinality, ergo $\rho$ is a member of $\acc(C)$, as needed.
\end{proof}
\begin{remark}
Following up  \cite[Question~3.8]{GoldbergHODconjecture} Goldberg asks if Clause~(5) from his \cite[Theorem~3.6(5)]{GoldbergHODconjecture} can be ruled out whenever $\kappa$ is a $\HOD$-singular cardinal; more precisely, he asks if there can be a $V$-cardinal $\kappa$ that is  $\HOD$-singular and such that $\omega <\cf(\kappa^{+\HOD})<|\kappa^{+\HOD}|^V$. The above theorem thus answers Goldberg's  question in the affirmative.
\end{remark}

\section*{Acknowledgements and AI disclosure}
No use of AI was made during the completion of this project. Poveda presented the results of this paper at the Simons Semester \emph{``Gödel's Program''} held at IMPAN. The authors are very grateful to the participants  for their comments and constructive feedback. Poveda
was supported by a 2026 Leonardo Grant for Scientific Research and Cultural Creation from the BBVA Foundation \includegraphics[width=0.1\linewidth]{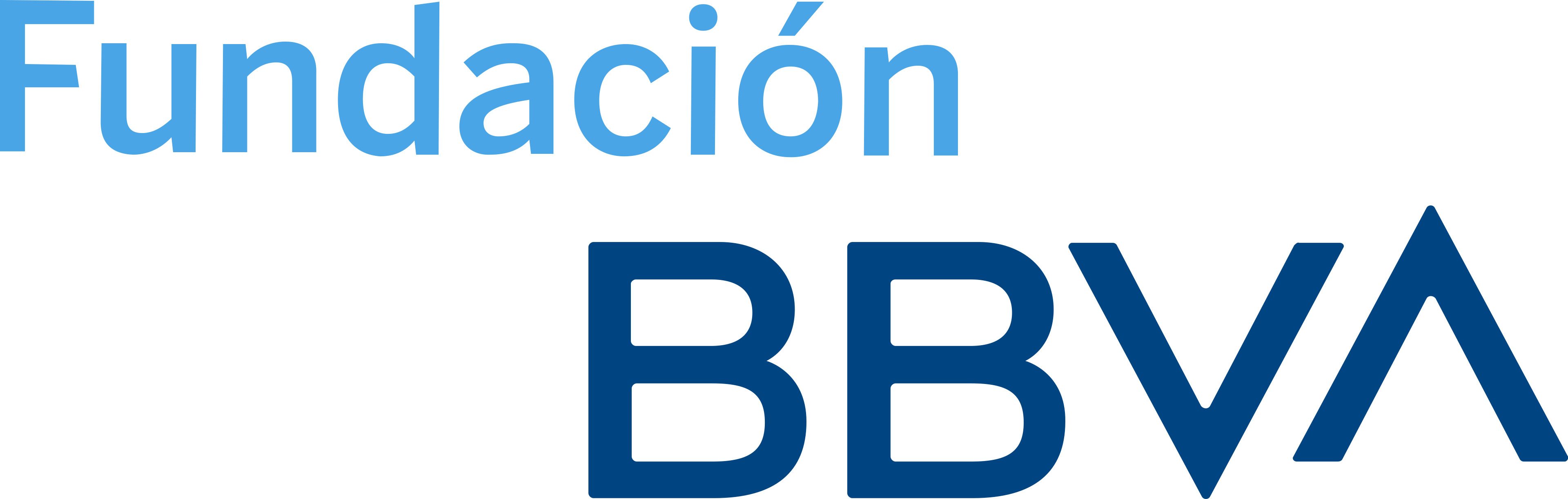} \footnote{The BBVA Foundation accepts no responsibility for the opinions, statements and contents included in the project and/or the results thereof, which are entirely the responsibility of the authors.}, by project  PID2023-147428NB-I00 from the Spanish Government and by the Alexander von Humboldt Foundation\includegraphics[width=0.08\linewidth]{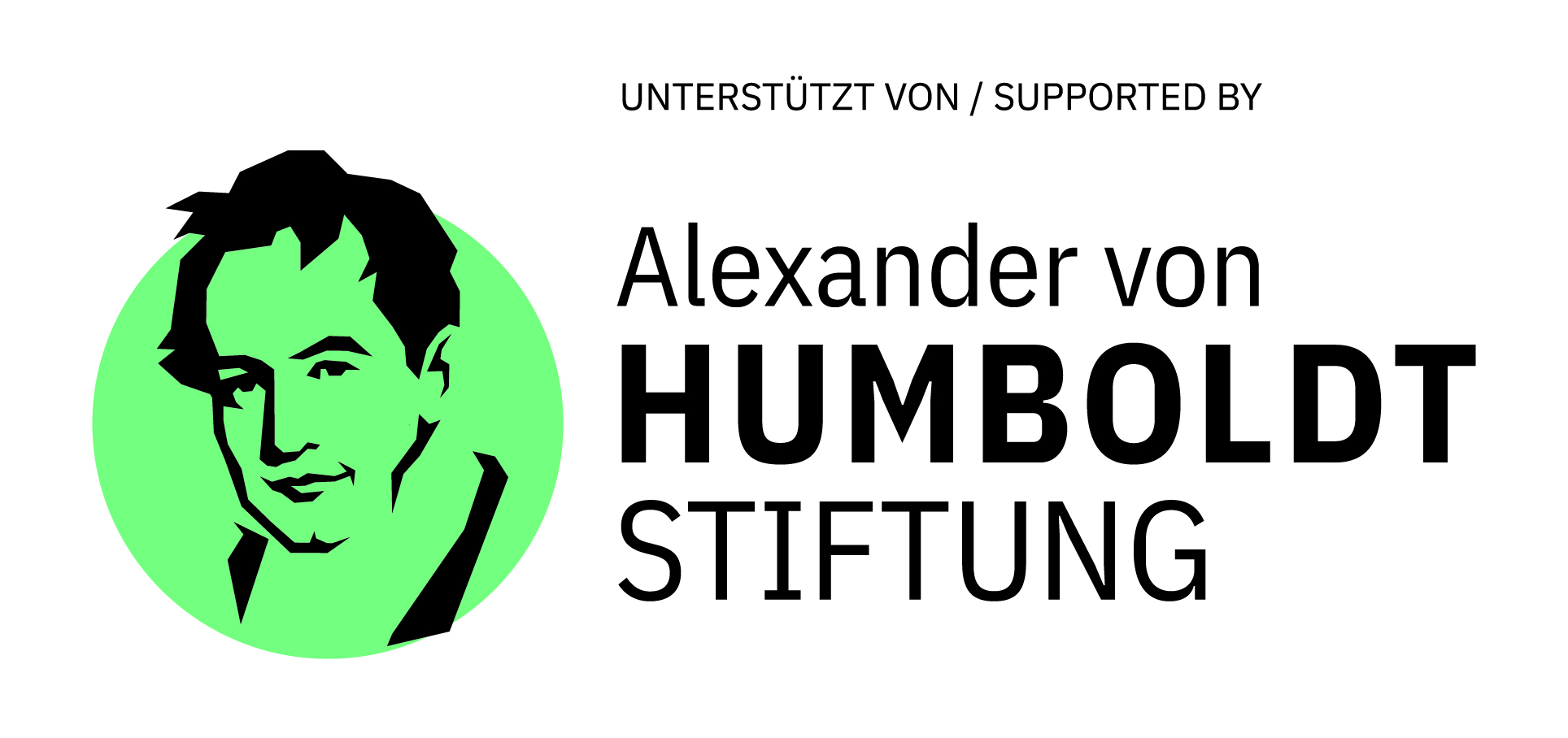}. Part of this work was developed while Jakob visited Poveda in Universität Hamburg and CUNEF Universidad. Both authors are grateful to the BBVA foundation and the AvH foundation for sponsoring these visits. 

\bibliographystyle{alpha} 
\bibliography{biblio}

\newcommand{\etalchar}[1]{$^{#1}$}
\begin{thebibliography}{BCH{\etalchar{+}}26}

\bibitem[AFF21]{MoreonHOD}
Arthur~W Apter, Shoshana Friedman, and Gunter Fuchs.
\newblock More on {HOD}-supercompactness.
\newblock {\em Annals of Pure and Applied Logic}, 172(3):102901, 2021.

\bibitem[AS83]{AbrahamShelahForcingClub}
Uri Abraham and Saharon Shelah.
\newblock Forcing closed unbounded sets.
\newblock {\em The Journal of Symbolic Logic}, 48(3):643--657, 1983.

\bibitem[Bag12]{Bag}
Joan Bagaria.
\newblock ${C^{(n)}}$-cardinals.
\newblock {\em Archive for Mathematical Logic}, 51(3):213--240, 2012.

\bibitem[BCH{\etalchar{+}}26]{CompactnessMagidor}
Tom Benhamou, James Cummings, Yair Hayut, Gabriel Goldberg, and Alejandro
  Poveda.
\newblock Compactness phenomena in hod and the optimality of magidor's covering
  theorem.
\newblock {\em arXiv preprint arXiv:2608.24190}, 2026.

\bibitem[BG24]{BagGol}
Joan Bagaria and Gabriel Goldberg.
\newblock Reflecting measures.
\newblock {\em Advances in Mathematics}, 443:109586, 2024.

\bibitem[BLS25]{Nairian}
Douglas Blue, Paul~B Larson, and Grigor Sargsyan.
\newblock Nairian models.
\newblock {\em arXiv preprint arXiv:2501.18958}, 2025.

\bibitem[BNH23]{BenHay}
Omer Ben-Neria and Yair Hayut.
\newblock On $\omega $ -strongly measurable cardinals.
\newblock {\em Forum of Mathematics, Sigma}, 11:e19, 2023.

\bibitem[BNU17]{BenUng}
Omer Ben-Neria and Spencer Unger.
\newblock Homogeneous changes in cofinalities with applications to {HOD}.
\newblock {\em Journal of Mathematical Logic}, 17(02):1750007, 2017.

\bibitem[BP23]{BagPov}
Joan Bagaria and Alejandro Poveda.
\newblock More on the preservation of large cardinals under class forcing.
\newblock {\em The Journal of Symbolic Logic}, 88(1):290--323, 2023.

\bibitem[BT09]{Broo2}
Andrew~D. Brooke-Taylor.
\newblock Large cardinals and definable well-orders on the universe.
\newblock {\em Journal of Symbolic Logic}, 74(2):641--654, 2009.

\bibitem[CFG15]{CumFriGol}
James Cummings, Sy~David Friedman, and Mohammad Golshani.
\newblock Collapsing the cardinals of {H}{O}{D}.
\newblock {\em Journal of Mathe\-matical Logic}, 15(02):1550007, 2015.

\bibitem[CIP25]{CortesPoveda}
Manuel Cort{\'e}s-Izurdiaga and Alejandro Poveda.
\newblock Almost free modules, perfect decomposition and enochs’ conjecture:
  M. cort{\'e}s-izurdiaga and a. poveda.
\newblock {\em Mediterranean Journal of Mathematics}, 22(2):51, 2025.

\bibitem[DF08]{DF}
Natasha Dobrinen and Sy-David Friedman.
\newblock Homogeneous iteration and measure one covering relative to {HOD}.
\newblock {\em Archive for Mathematical logic}, 47(7):711--718, 2008.

\bibitem[DPT25]{DPT}
Vincenzo Dimonte, Alejandro Poveda, and Sebastiano Thei.
\newblock The {B}aire and perfect set properties at singulars cardinals.
\newblock {\em To appear in Israel J. of Math.}, 2025.

\bibitem[EM02]{Eklof}
Paul~C. Eklof and Alan~H. Mekler.
\newblock {\em Almost free modules: set-theoretic methods}.
\newblock Elsevier, 2002.

\bibitem[FHR15]{Geology}
Gunter Fuchs, Joel~David Hamkins, and Jonas Reitz.
\newblock Set-theoretic geology.
\newblock {\em Annals of Pure and Applied Logic}, 166(4):464--501, 2015.

\bibitem[Git86]{GitikChanging}
Moti Gitik.
\newblock Changing cofinalities and the nonstationary ideal.
\newblock {\em Israel Journal of Mathematics}, 56:280--314, 1986.

\bibitem[Git99]{GitikOnClosed}
Moti Gitik.
\newblock On closed unbounded sets consisting of former regulars.
\newblock {\em The Journal of Symbolic Logic}, 64(1):1--12, 1999.

\bibitem[GK23a]{GitikKaplan2}
Moti Gitik and Eyal Kaplan.
\newblock Non-stationary support iterations of prikry forcings and restrictions
  of ultrapower embeddings to the ground model.
\newblock {\em Annals of Pure and Applied Logic}, 174(1):103164, 2023.

\bibitem[GK23b]{GitikKaplan}
Moti Gitik and Eyal Kaplan.
\newblock On restrictions of ultrafilters from generic extensions to ground
  models.
\newblock {\em The Journal of Symbolic Logic}, 88(1):169--190, 2023.

\bibitem[GK24]{gitik2024fresh}
Moti Gitik and Eyal Kaplan.
\newblock On fresh sets in iterations of prikry type forcing notions.
\newblock {\em arXiv preprint arXiv:2403.01617}, 2024.

\bibitem[GM18]{GitikMer}
Moti Gitik and Carmi Merimovich.
\newblock Some applications of supercompact extender based forcings to {HOD}.
\newblock {\em The Journal of Symbolic Logic}, 83(2):461--476, 2018.

\bibitem[G{\"o}d39]{GodelConstructible}
Kurt G{\"o}del.
\newblock Consistency-proof for the generalized continuum-hypothesis.
\newblock {\em Proceedings of the National Academy of Sciences},
  25(4):220--224, 1939.

\bibitem[G{\"o}d46]{Godel}
Kurt G{\"o}del.
\newblock Remarks before the princeton bicentennial conference on problems in
  mathematics.
\newblock 1946.

\bibitem[Gol24a]{GoldbergHODconjecture}
Gabriel Goldberg.
\newblock The {HOD} conjecture and its failure.
\newblock 2024.

\bibitem[Gol24b]{GoldbergHODStrongCompact}
Gabriel Goldberg.
\newblock Strongly compact cardinals and ordinal definability.
\newblock {\em Journal of Mathematical Logic}, 24(01):2250010, 2024.

\bibitem[GOP24]{GOP}
Gabriel Goldberg, Jonathan Osinski, and Alejandro Poveda.
\newblock On the optimality of the {HOD} dichotomy.
\newblock {\em arXiv preprint arXiv:2411.03558}, 2024.

\bibitem[GP25]{GolPov}
Gabriel Goldberg and Alejandro Poveda.
\newblock Compactness phenomena in {HOD}.
\newblock 13:e118, 2025.

\bibitem[HK16]{HayKag}
Yair Hayut and Asaf Karagila.
\newblock Restrictions on forcings that change cofinalities.
\newblock {\em Archive for Mathematical Logic}, 55(3):373--384, 2016.

\bibitem[HP26a]{HP2}
Yair Hayut and Alejandro Poveda.
\newblock The directedness of the rudin-keisler order at measurable cardinals.
\newblock {\em arXiv preprint arXiv:2601.10614}, 2026.

\bibitem[HP26b]{HP}
Yair Hayut and Alejandro Poveda.
\newblock The gluing property.
\newblock {\em Journal of Mathematical Logic}, 26(01):2450030, 2026.

\bibitem[Kan09]{Kan}
Akihiro Kanamori.
\newblock {\em The higher infinite}.
\newblock Springer Monographs in Mathematics. Springer-Verlag, Berlin, second
  edition, 2009.
\newblock Large cardinals in Set Theory from their beginnings, Paperback
  reprint of the 2003 edition.

\bibitem[Koe11]{Koel}
Peter Koellner.
\newblock Independence and large cardinals.
\newblock In Edward~N. Zalta, editor, {\em The Stanford Encyclopedia of
  Philosophy}. Metaphysics Research Lab, Stanford University, summer 2011
  edition, 2011.

\bibitem[Kun14]{Kun}
Kenneth Kunen.
\newblock {\em Set theory: an introduction to independence proofs}, volume 102.
\newblock Elsevier, 2014.

\bibitem[McA71]{McAloon}
Kenneth McAloon.
\newblock Consistency results about ordinal definability.
\newblock {\em Annals of Mathematical Logic}, 2(4):449--467, 1971.

\bibitem[Mit84]{MitIter}
William~J Mitchell.
\newblock The core model for sequences of measures. {I}.
\newblock In {\em Mathematical Proceedings of the Cambridge Philosophical
  Society}, volume~95, pages 229--260. Cambridge University Press, 1984.

\bibitem[Mit01]{Mitchellokappa}
William~J. Mitchell.
\newblock A measurable cardinal with a closed unbounded set of inaccessibles
  from {$o(\kappa)=\kappa$}.
\newblock {\em Transactions of the American Mathematical Society},
  353(12):4863--4897, 2001.

\bibitem[MS96]{MagShe}
Menachem Magidor and Saharon Shelah.
\newblock The tree property at successors of singular cardinals.
\newblock {\em Archive for Mathematical Logic}, 35(5-6):385--404, 1996.

\bibitem[Myh71]{Myhill}
John Myhill.
\newblock Ordinal definability john myhill and dana scott.
\newblock {\em Axiomatic Set Theory, Part 1}, 1:271, 1971.

\bibitem[Pov24]{PovOmega}
Alejandro Poveda.
\newblock Two results on extendible cardinals.
\newblock {\em Proceedings of the American Mathematical Society},
  152(06):2669--2681, 2024.

\bibitem[PT25]{PT}
Alejandro Poveda and Sebastiano Thei.
\newblock Combinatorics in higher solovay models.
\newblock {\em arXiv preprint arXiv:2509.18991}, 2025.

\bibitem[Sar08]{Sar08}
Grigor Sargsyan.
\newblock On hod-supercompactness.
\newblock {\em Archive for Mathematical Logic}, 47(7-8):765--768, 2008.

\bibitem[Sol70]{Solmod}
Robert~M. Solovay.
\newblock A model of {S}et {T}heory in which every set of reals is {L}ebesgue
  measurable.
\newblock {\em Annals of Mathematics}, pages 1--56, 1970.

\bibitem[Ste09]{Steel}
John~R Steel.
\newblock An outline of inner model theory.
\newblock In {\em Handbook of set theory}, pages 1595--1684. Springer, 2009.

\bibitem[VH72]{Vopenka}
P.~Vop{\v{e}}nka and P.~H{\'a}jek.
\newblock {\em The Theory of Semisets}.
\newblock Studies in Logic. North-Holland Publishing Company, 1972.

\bibitem[Woo10]{WOEM1}
W.~Hugh Woodin.
\newblock Suitable {E}xtender {M}odels {I}.
\newblock {\em Journal of Mathematical Logic}, 10(01n02):101--339, 2010.

\bibitem[Woo17]{midrasha}
W~Hugh Woodin.
\newblock In search of {U}ltimate-{L} the 19th midrasha mathematicae lectures.
\newblock {\em Bulletin of symbolic logic}, 23(1):1--109, 2017.

\bibitem[Zem02]{Zeeman}
Martin Zeman.
\newblock {\em Inner models and large cardinals}, volume~5.
\newblock Walter de Gruyter, 2002.

\end{thebibliography}
\end{document}